\ifx\shlhetal\undefinedctrlseq\let\shlhetal\relax\fi
\documentclass{article} 

\usepackage{amsbsy}    
\usepackage{amsfonts}  
\usepackage{amssymb}   

\newcommand{\blambda}{{\boldsymbol \lambda}}
\newcommand{\bmu}{{\boldsymbol \mu}}
\newcommand{\Col}{{\rm Col}}
\newcommand{\cl}{{\rm cl}}

\newcommand{\tim}{{\rm times}}

\newcommand{\cf}{{\rm cf}\/}
\newcommand{\hcf}{{\rm hcf}\/}
\newcommand{\dom}{{\rm Dom\/}}

\newcommand{\reg}{{\rm Reg}}
\newcommand{\rk}{{\rm rk}}
\newcommand{\fil}{{\rm Fil}}
\newcommand{\ord}{{\rm Ord}}
\newcommand{\eub}{{\it eub}}
\newcommand{\lub}{{\it lub}}
\newcommand{\rang}{{\rm Rang\/}}
\newcommand{\mod}{{\rm \ mod\:}}
\newcommand{\hpp}{{\it hpp}}
\newcommand{\epp}{{\it epp}}
\newcommand{\ehpp}{{\it ehpp}}

\newcommand{\spp}{{\it spp}}
\newcommand{\shpp}{{\it shpp}}
\newcommand{\pcf}{{\rm pcf}} 
\newcommand{\tlim}{{\rm tlim}}
\newcommand{\ts}{{\bf Ts}}
\newcommand{\tw}{{\bf Tw}}
\newcommand{\1}{{\bf 1}}
\newcommand{\pp}{{\it pp}}

\newcommand{\cP}{{\cal P}}
\newcommand{\bd}{{\rm bd}}
\newcommand{\bT}{{\bf T}}
\renewcommand{\i}{\iota}

\newcommand{\forces}{\Vdash}
\newcommand{\nforces}{\nVdash}
\newcommand{\otp}{{\rm otp}}

\newcommand{\A}{{\rm A}^*}

\newcount\skewfactor
\def\mathunderaccent#1#2 {\let\theaccent#1\skewfactor#2
\mathpalette\putaccentunder}
\def\putaccentunder#1#2{\oalign{$#1#2$\crcr\hidewidth
\vbox to.2ex{\hbox{$#1\skew\skewfactor\theaccent{}$}\vss}\hidewidth}}
\def\name{\mathunderaccent\tilde-3 }

\renewcommand{\v}{{\bf V}}
\newcommand{\gc}{{\mathfrak C}}
\newcommand{\QED}{\vrule width 6pt height 6pt depth 0pt 
\vspace{0.1in}}
\newcommand{\qed}{\QED} 

\newcommand{\Proof}{\noindent{\sc Proof:} \hspace{0.25in}} 

\newtheorem{theorem}{Theorem}[section] 
\newtheorem{lemma}[theorem]{Lemma} 
\newtheorem{notation}[theorem]{Notation}

\newtheorem{remark}[theorem]{Remark} 
\newtheorem{definition}[theorem]{Definition}
\newtheorem{claim}[theorem]{Claim} 
\newtheorem{observation}[theorem]{Observation} 

\newtheorem{convention}[theorem]{Convention}
\newtheorem{conclusion}[theorem]{Conclusion}

\title{SET THEORY WITHOUT CHOICE: NOT EVERYTHING ON COFINALITY IS POSSIBLE} 
\author{
{\bf Saharon Shelah}\thanks{Research supported by  ``The Israel
Science Foundation'' administered by The Israel Academy of Sciences
and Humanities. Publication no 497, done 11/92 -- 1/93.}\\
Institute of Mathematics\\
The Hebrew University of Jerusalem\\
Jerusalem, Israel\\
and\\
Department of Mathematics\\
Rutgers University,\\
New Brunswick NJ, USA
}

\date{
\today
}

\begin{document}

\maketitle 

\begin{abstract}
We prove in ZF+DC, e.g. that: if $\mu=|H(\mu)|$ and
$\mu>\cf(\mu)>\aleph_0$ \underline{then} $\mu ^+$ is regular but non
measurable. This is in contrast with the results on measurability for
$\mu=\aleph_\omega$ due to Apter and Magidor [ApMg].
\end{abstract}

\vfill
\eject

\centerline{{\sc Annotated Content}}

\par \noindent
{\bf \S0 Introduction}

\begin{quotation}
[In addition to presenting the results and history, we gave some basic
definitions and notation.]
\end{quotation}

\par \noindent
{\bf \S1 Exact upper bound}

\begin{quotation}
[We define some variants of least upper bound ($\lub$, $\eub$) in
$({}^{(\A)}\ord, <_D)$ for $D$ a filter on $\A$. We consider $<_D$-increasing
sequence indexed by ordinals $\delta$ or indexed by sufficiently directed
partial orders $I$, of members of ${}^{(\A)}\ord$ or of members of
${}^{(\A)}\ord/D$ (and cases in the middle). We give sufficient conditions for
existence involving large cofinality (of $\delta$ or $I$) and some amount of
choice.  Mostly we look at what the $ZFC$ proof gives without choice.  Note in
particular \ref{claim1:6}, which assumes only $DC$ ($ZF$ is not mentioned of
course), the filter is $\aleph_1$-complete and cofinality of $\delta$ large
and we find an $\eub$ only after extending the filter.]
\end{quotation}

\par \noindent
{\bf \S2 $\hpp$}

\begin{quotation}
[We look at various ways to measure the size of the set $\prod\limits_{a\in
\A} f(a)/D$, like supremum length of $<_D$-increasing sequence in
$\prod\limits_{a\in \A} f(a)$ (called $\ehpp_D$), or in $\prod\limits_{a\in
\A} f(a)/D$ (called $\hpp_D$), or we can demand them to be cofinal (getting
$\epp$ or $\pp$); when we let $D$ vary on $\Gamma$ we write
$\ehpp_\Gamma$ etc.  
So existence of $<_D-\eub$ give downward results on $\pp$.]
\end{quotation}

\par \noindent
{\bf \S3 Nice family of filters}

\begin{quotation}
[In this paper we can say little on products of countably many; we can say
something when we deal with $\aleph_1$-complete filters. So as in [Sh-g], V,
we deal with family $E$ of filters on $\A$ which is nice.  Hence suitable
ranks from function $f$ from $\A$ to ordinal and filter $D\in E$ are well
defined (i.e. the values are ordinals not infinity).  The basic properties of
those ranks are done here.

We then define some measures for the size of $\prod\limits_{a\in \A} f(a)/D$
(i.e. ${\it Tw}$, ${\it Ts}$, $T$), looking at subsets of 
$\prod\limits_{a\in \A}f(a)$ or of $\prod\limits_{a\in \A} f(a)/D$ which are
pairwise $\neq_D$. In conclusion \ref{conclusion3:??} we, under reasonable
assumption, prove that some such measures and $\sup\{\rk^2_D(f): D\in E\}$ are
equal. In \ref{claim3:12} we have a parallel of \ref{claim1:6}: sufficient
condition for the existence of $\eub$ when we allow to increase the filter. We
end defining normal filters and a generalization.

The basic point is that for every $f\in {}^{(\A)}\ord$ (if $E$ is nice) for
some $D\in E$ we have $\rk^2_D(f)< \rk^3_D(f)$ and in this case the rank
determine $f/D$, and order on rank the order among such $f$; so we can
represent $\prod\limits_{a\in \A} f^*(a)/D$ as the union of $\leq |E|$ well
ordered sets.]
\end{quotation}

\par \noindent
{\bf \S4 Investigating strong limit singular}

\begin{quotation}
[We deal with $\otimes_{\alpha, R}$, which means that we can regard $R\cap
\alpha$ as a substitute of the family of regulars (not just individually) i.e.
we can find $\langle e(i): i<\alpha\rangle$, $e(i)$ is an unbounded subset of
$i$ of order type which belongs to $R$. We give the basic properties in
\ref{observation4:2}, then move up from $\mu$ to $\rk^2_D(\mu)$ with
appropriate choice of $R$. With this we have a parallel of Galvin--Hajnal
theorem (\ref{obs6:1}). Here comes the main theorem (\ref{lemma4:5}) assuming
$DC$, $\mu$ singular of uncountable cofinality (and $E$ a set of filters on
$\cf(\mu)$ which is nice), if $\mu=|{\cal H}(\mu)|$ (kind of strong limit)
\underline{then} 
set theory ``behave nicely'' up to $2^\mu$: $2^\mu$ is an aleph, there
is no measurable $\leq 2^\mu$ and $\mu^+$ is regular.

We end defining some natural ideals in this context of $\otimes_{\mu,R}$.]
\end{quotation}

\par \noindent
{\bf \S5 The successor of a singular of uncountable cofinality}

\begin{quotation}
[Here we prove our second main theorem: if $\langle \lambda_i:i<\kappa\rangle$
is increasing continuous ($E$ a nice family of filters on $\kappa$) then for
at most one $\lambda$, for stationarily many $i<\kappa$, the cardinal
$\lambda^+_i$ has cofinality $\lambda$.]
\end{quotation}

\par \noindent
{\bf \S6 Nice $E$ exists}

\begin{quotation}
[The theorems so far have as a major assumption the existence of a nice $E$.
Using inner models we show that as in the situation with choice, if there is
no nice $E$ then the universe is similar enough to some inner model to answer
our questions on the exponentiation and cofinality.]
\end{quotation}

\setcounter{section}{-1}
\section{Introduction} 

Originally I disliked choiceless set theory but ([Sh 176]) discovering (the
first modern asymmetry of measure/category: a difference in consistency
strength and)
\begin{description}
\item[\mbox{$[ZF+DC]$}]\ \\
if there is a set of $\aleph_1$ reals then there is a Lebesgue non-measurable
set  
\end{description}
I have softened. Recently Gitik suggested to me to generalize the \pcf\ theory
to the set theory without choice, or more exactly with limited choice. E.g.:
is there a restriction on $\langle \cf(\aleph_n):n<\omega\rangle$? By Gitik
[Gi] $ZF+$ ``every aleph has cofinality $\aleph_0$'' is consistent.

It is known that if $ZF + DC + AD$ then there is a very specific pattern of
cofinality, but we have no flexibility. So we do not know e.g. if
$$
\mbox{``}ZF + DC + (\forall \delta)[ \cf(\delta)\leq \aleph_{\alpha}]+
(\forall\beta<\alpha)(\aleph_{\beta+1} \mbox{ is regular})\mbox{''} 
$$ 
is consistent for $\alpha =1$, or $\alpha =2$ etc. The general question is
what are the restrictions on the cofinality function.

Apter repeated the question above and told me of Apter and Magidor
[ApMg] in which 
the consistency of $ZF+DC_{\aleph_\omega}+|{\cal
H}(\aleph_{\omega})|=\aleph_{\omega} 
+\mbox{``}\aleph_{\omega+1}$ is measurable'' was proved (for an earlier weaker
result see Kafkoulis [Kf ]) and was quite confident that a parallel theorem
can be similarly proved for $\aleph_{\omega_1}$. 

My feeling was that while the inner model theory and the descriptive set
theory are not destroyed by dropping AC, modern combinatorial set theory says
essentially nothing in this case, and it will be nice to have such a theory.

Our results may form a modest step in this direction. The main result is
stated in the abstract.
\begin{theorem}[$ZF+DC$]
\label{A}
If $\mu$ is singular of uncountable cofinality and is strong limit in the
sense that ${\cal H}(\mu)$ has cardinality $\mu$ \underline{then} $\mu^+$ is
regular and non measurable.
\end{theorem}

Note that this work stress the difference between ``bounds for cardinal
arithmetic for singulars with uncountable cofinality'' and the same for
countable cofinality.

Another theorem (see section~\ref{par5}) says 
\begin{theorem}
\label{B}OB
If $\langle \mu_i:i\leq\kappa\rangle$ is an increasing continuous sequence of
alephs $>\kappa$, then for at most one $\lambda$, $\{i:\cf(\mu_i^+)=\lambda\}$
is a stationary subset of $\kappa$ (see \ref{definition3:?}~$(3)$). 
\end{theorem}

We were motivated by a parallel question in $ZFC$, asked by Magidor
and to which we do not know the answer: can $\aleph_{\omega_1}$ be
strong limit and for $\ell\in \{1, 2\}$ 
there is a stationary set $S_\ell\subseteq\omega_1$ such that
$\prod\limits_{\delta\in S_1} \aleph_{\delta+2}/J^{\bd}_S$ has true cofinality
$\aleph_{\omega_1+\ell}$?

It was known that ``there are two successive singulars'' has large
consistency strength.

We do not succeed to solve Woodin's problem (show consistency of ``$ZF+
DC_\omega+$ every aleph has cardinality $\aleph_0$ or $\aleph_1$''), and
Specker's problem (show consistency of ``$ZF+$ every $\cP (\kappa)$ is the
union of countably many sets of cardinality $\leq\kappa$''). For me the real
problems are:
\begin{description}
\item[(a)] ($ZF$) Is there an aleph $\kappa$ (i.e. suitable definition) such
that $DC_\kappa$ implies that the class of regular alephs is unbounded?
\item[(b)] For which $\kappa<\lambda$, does $DC_{<\kappa}$ implies $\lambda$
is regular? 
\item[(c)] ($ZF$) Is there $\kappa$ such that $DC_\kappa$ implies that for
every $\mu$ for a class of alephs $\lambda$, the set ${\cal P}(\lambda)$ is
not the union of $\leq\mu$ sets each of cardinality $\lambda$ and
$\cf(\lambda)\geq\kappa$ ?
\end{description}

We try to assume only standard knowledge in set theory but we start by
discussing what part of pcf theory survives (see [Sh-g] and history there) so
knowledge of it will be helpful; in particular section 6 imitates [Sh-g, V,
\S1] so it also uses a theorem of Dodd and Jensen [DJ]. On set theory without
full choice see [J1], on cardinal arithmetic and its history [Sh-g].  Lately
Apter and Magidor got consistency results complementary to ours and we
intend to return to the problems here in [Sh 589, \S5].

\begin{definition}
A cardinality is called an aleph if there is an ordinal with that cardinality,
or just the cardinality of a set which can be well ordered and then it is
identified with the least ordinal of that cardinality.
\end{definition}

\begin{notation}
\label{not0:2}
{\em
$\alpha ,\beta,\gamma,\epsilon,\zeta,\xi,i,j$ are ordinals; $\delta$ a limit
ordinal; $\lambda,\mu,\kappa,\chi,\theta ,\sigma$ are cardinals (not
necessarily alephs),  

$\reg$ is the class of regular alephs (see Definition \ref{1.1}(7)),

$D$ a filter on the set $\dom(D)$,

$A=\emptyset \mod D$ means $\dom(D) \setminus A \in D$,

$D^+=\{A:A \subseteq \dom(D)$, and $A \neq \emptyset \mod D\}$,

$D+A=\{X \subseteq \dom(D): X \cup (\dom(D) \setminus A) \in D\}$.

For a filter $D$, $\theta_D =\min \{\kappa$ : there is no function $h:\dom(D)
\longrightarrow\kappa$ such that for every $ i<\kappa$ we have $h^{-1}({i})
\neq\emptyset\mod D\}$. 

$\cal I , \cal J$ denote ideals on the set $\dom(\cal I), \dom(\cal J)$
respectively; definitions given for filters $D$ apply also to the dual ideal
$\{\dom(D) \setminus A: A \in D\}$ but ${\cal I}+A=\{B\subseteq\dom(\cal I):
B\setminus A\in {\cal I}\}$.  

$I, J$ are directed partial orders or just index sets.

A cone of a partial order $I$ is a subset of the form $\{a\in I:I\models
a_0\leq a\}$ for some $a_0\in I$. 

For a set $A$ let $\theta(A)$ be 
$$
\sup\{\alpha:\alpha\mbox{ is an ordinal and there is a function from $A$ onto
} \alpha\mbox{ or }\alpha \mbox{ is finite}\}. 
$$
Note that if $|A|$ is an aleph then $\theta(A)=|A|^+$. Also $\theta(A)$ is an
aleph. If $g:A\longrightarrow\delta$ has unbounded range then
$\cf(\delta)\leq\theta(A)$ (and $|\rang(y)|$ is an aleph $\leq
\theta(A)$) 
see Definition \ref{1.1}(7).

Let $\theta^-(A)$ be 
$$
\min\{\alpha: \alpha\mbox{ an ordinal and there is no one to one function from
}\alpha\mbox{ into }A\}; 
$$ 
clearly $\theta^-(A^*)$ is an aleph and $\theta^-(A)\leq \theta(A)$.

For a directed partial order $I$, ${\cal J}^{\rm bd}_I$ is the ideal of
bounded subsets of $I$ and ${\cal D}^{\rm bd}_{I}$ is the dual filter (usually
$I$ is a regular aleph).

If $f,g :\A\longrightarrow\ord$ and $D$ is a filter on $\A$ then 
$$
f\leq_D g\quad\mbox{ means }\quad\{a\in\A:f(a)\leq g(a)\}\in D,
$$ 
similarly for other relations ($<_D$,$=_D$, $\neq_D$). Note: $\leq_D$ is a
quasi order (as maybe $f_0\leq_D f_1\leq_D f_0$, $f_0\neq f_1$) but $\langle
{}^{(\A)}\ord /D, <_D \rangle$ is a partial order and $\neq_D$ is not the
negation of $=_D$. } 
\end{notation}

\begin{definition}
Let $|A|\leq^*|B|$ mean that {\em there is a function from $B$ onto $A$} (so
$|A|\leq |B| \Rightarrow |A|\leq^*|B|$ {\em but not necessarily the
converse}). (Note: for well ordered sets $\leq$, $<^*$ are equal, in fact
``$B$ is well ordered'' is enough.)
\end{definition}

\begin{definition}\
\begin{enumerate}
\item $AC_{\blambda,{\bar\bmu}}$ is the axiom of choice for every family
$\{A_t:t\in I\}$ of sets, $|A_t|\leq\mu_t$, where $\bar
\bmu=\langle\bmu_t:t\in I\rangle$, $I$ a set of cardinality $\leq \lambda$.
If $\bmu_t=\bmu$ for all $t\in I$ then we write $AC_{\blambda,\bmu}$. We may
write $AC_{I,\bmu}$, $AC_{I,\bar A}$ instead, similarly below. 

\item If in part (1), $I$ is a well ordering (e.g. an ordinal) then
$DC_{I,\bar \bmu}$ is the dependent choice version.

$DC_{I,\bmu}$ is $DC_{I,{\bar \bmu}}$ where $(\forall{t\in I})
\bmu_t=\bmu.$ 

\item Let $AC_\blambda$ mean $\forall\bmu AC_{\alpha,\bmu}$; and
$AC_{<\blambda}\equiv(\forall{\bmu<\blambda}) \ AC_\bmu$; and
$AC_{\blambda,<\bmu}$, $AC_{<\blambda,\bmu}$, $AC_{<\lambda,<\mu}$,
$AC_{\blambda, \bar \bmu}$, $AC_{<\blambda, \bar \bmu}$, $AC_{<\blambda, <\bar
\bmu}$, have the natural meanings.
\item $DC_{\alpha}$ means $\forall\bmu \ DC_{\alpha,\bmu}$, and
$DC_{<\alpha}$, $DC_{\alpha,<\bmu}$, $DC_{<\alpha,\bmu}$ and
$DC_{<\alpha,<\bmu}$ have the natural meanings. $DC$ is $DC_{\aleph_0}$. 
\end{enumerate}
\end{definition}

\section{Exact upper bound}

\begin{definition}\ 
\label{1.1}
\begin{enumerate}
\item  A partial order $I$ is $(\leq \lambda)$-directed (or
$(<\lambda)$-directed) \underline{if} for every $A\subseteq I$ of cardinality
$\leq \lambda$ (of cardinality $<\lambda$) there is an upper bound. If
$\lambda=2$ we omit it. Note that $2$-directed is equivalent to
$(<\aleph_0)$-directed.   
\item We say that $J$ is cofinal in $I$ if $J\subseteq I$ and $(\forall s \in
I)(\exists t \in J) [I\vDash s \leq t]$. 
\item I is endless if $(\forall s \in I)(\exists t \in I) [I\vDash s < t]$.
\item We say $\blambda \geq\cf(I)$ \underline{if} there is a cofinal
$J\subseteq I$ of cardinality $\leq\lambda$ and $\blambda=\cf(I)$ if
$\blambda$ is the smallest such $\lambda$ i.e. $(\forall\blambda_1)(\lambda_1
\geq \cf(I) \rightarrow \blambda_1\geq \blambda)$ (it does not necessarily
exist, but is unique by the Cantor-Bernstein theorem).
\item Let $F$ be a set of functions from $\A$ to ordinals, $D$ a filter on
$\A$.  We say $g:\A\rightarrow \ord$ is a $<_D -\eub$ (exact upper bound) of
$F$ if:
\begin{description}
\item[(i)\ \ ] \ \ $f\in F\Rightarrow f \leq_D g$
\item[(ii)\ ]\ \ If $h:\A\rightarrow\ord$, $h<\max\{g,1\}$\\
\underline{then} for some $f\in F$ we have $h<_D\max\{f,1\}$ (where
$f^*=\max\{f, 1\}$ means $f^*$ is a function, $\dom(f^*)=\dom(f)$ and for
every $x\in \dom(f^*)$ we have $f^*(x)=\max\{f(x), 1\}$). 
\end{description}
\item Let $I$ be ${(\leq \blambda)}^*$-directed if in $(1)$ we ask $|A|
{\leq}^*\lambda$. We define similarly ``${\cf}^*(I)\leq \blambda$'', and also
``$I$ is ${(< \blambda)}^*$- directed'', ``${\cf}^*(I)<\blambda$'',
``${\cf}^*(I)=\blambda$'' (e.g. ${\cf}^*(I)=\blambda$ means $(\forall
\blambda_1)[\blambda_1 \geq \cf^*({\cal I})\rightarrow \blambda\leq^*
\blambda_1]$), so actually we should say ${\cf}^*(I)=\blambda/{\equiv^*}$ where
$\blambda_1{\equiv^*} \blambda_2$ iff $\blambda_1 \leq^* \blambda_2 \leq^*
\blambda_1$. Instead ``$(\leq \blambda)^*$-'' we may write ``$(\leq^* \blambda
)$-''. 
\item The ordinal $\delta$ is regular \underline{if} $\delta=\cf(\delta)$
where  
$$
\cf(\alpha)=\min\{\otp(A): A\subseteq\alpha\mbox{ is unbounded}\}.
$$ 
Clearly if $I$ is well ordered \underline{then} $\cf(I)$ is a regular ordinal,
and a regular ordinal is a cardinal.
\item We say $A$ is unbounded in $I$ if $A\subseteq I$ and for no $t\in
I\setminus A$ do we have $(\forall {s\in A})(s\leq_I t)$.
\item $\cf_- (I)\leq\lambda$ \underline{if} there is an unbounded $A\subseteq
I$, $|A|\leq\lambda$; similarly $\cf_-^* (I)<\lambda$, $\cf_-^* (I)=\lambda$.
\item $\hcf_- (I)\leq\lambda$ (the hereditary $\cf_- (I)$ is $\leq
\lambda$) \underline{if} for 
every $J\subseteq I$, $\cf_- (J)\leq\lambda$; similarly $\hcf_-^*
(I)<\lambda$, $\hcf_-^* (I)=\lambda$.
\end{enumerate}
\end{definition}

\begin{definition}\ 
\label{def1:7}
Let $I$ be a $(<\aleph_0)$-directed (equivalently a directed) partial order
(often $I$ will be a limit ordinal $\delta$ with its standard order, then we
write just $\delta$). 
\begin{enumerate}
\item $\bar F=\langle F_t: t\in I\rangle$ is $<_D$-increasing (or
$\leq_D$-increasing) \underline{if}: 

(a) $\bigcup\limits_{t\in I} F_t \subseteq {}^{\dom(D)} \ord$ and
$F_t\neq\emptyset$,

(b) if $f_e \in F_{t_{e}}$, (for $e<2$) and $t_0 <_I t_1$ then $f_0 / D <f_1/
D$ (or $f_0/D\leq f_1/D$).
\item ${\bar F} /D$ is smooth (or $\bar F$ is $D$-smooth) \underline{if}:

(a) $\bigcup\limits_{t\in I} F_t \subseteq {}^{\dom (D)} \ord$ and
$F_t\neq\emptyset$,

(c) $f_1, f_2 \in F_t\Rightarrow f_1 / D= f_2 /D$.

\item Let $\bar F=\langle F_t : t\in I\rangle$, we say $\bar F /D$ is
semi-smooth (or $\bar F$ is semi-$D$-smooth) if: 

(a) $\bigcup\limits_{t\in I} F_t \subseteq {}^{\dom (D)} \ord$ and
$F_t\neq\emptyset$,

$(c)'$ if $I\models \mbox{``}t_0<t_1< t_2< t_3$'' and $f_l\in F_{t_l}$ for
$l<4$ \underline{then} 
$$
\{a\in\dom (D): f_1 (a)< f_2 (a)\}\subseteq \{a\in\dom (D): f_0
(a)< f_3 (a)\}\mod D,
$$
\item $\bar F$ is almost $D$-smooth (or ${\bar F}/D$ is almost smooth)
is defined similarly but 

$(c)''$ if $I\models t_0<t_1$ and $f_0,{f'}_0\in F_{t_0}$, and $f_1,{f'}_1\in 
F_{t_1}$ \underline{then} 
$$
\{a\in\dom (D): f_0 (a)< f_1 (a)\}=\{a\in\dom (D): {f'}_0
(a)< {f'}_1 (a)\}\mod D,
$$
\item $f\in {}^{\dom(D)}\ord$ is a $\lub$ (least upper bound) of $\bar F/D$
(or $<_D-\lub$ of $\bar F$) if:

(a) $f \in \bigcup\limits_{t\in I} F_t \Rightarrow f\leq_D g$,

(b) if $g' \in {}^{\dom (D)} \ord$ satisfies (a) too then $g\leq_D g'$.
\item $g \in {}^{\dom(D)} \ord$ is  the $eub$ of $\bar F/ D$ (or $<_D-\eub$ of
$\bar F$) if  

(a) $f\in \bigcup\limits_{t\in I} F_t \Rightarrow f <_D g$,

(b) if $f^* <_D \max \{g,\1_{\A}\}$, then $f^* <_D \max\{f,\1_{\A}\}$; note
that if for some $g\in F$ we have $g\neq_D 0_{A^*}$ then this is equivalent
to: if $f^* <_D g$ then 
$$
(\exists f\in \bigcup\limits_{t\in I} F_t)(f^* <_D f)
$$
(note: if $\neg (g\neq_D {\bf 0}_{\A})$ this holds vacuously).
\end{enumerate}
\end{definition}

\begin{observation}
\label{obs1:2}
\begin{enumerate}
\item We have $\cf (0)=0$, $\cf (\alpha +1)=1$. For a limit ordinal $\delta$,
$\cf(\delta)=\cf^*(\delta)$ is regular and infinite; each regular ordinal is
an aleph.  For a linear order $I$, ${\cf}^*(I) \leq \lambda$ iff $I$ is not
$(\leq \lambda)^*$-directed. For a limit ordinal $\delta$ we have:
$\cf(\delta)\leq^* |A|$ iff $\cf(\delta)<\theta (A)$.
\item If $\bar F$ is $<_D$-increasing \underline{then} it is semi-$D$-smooth.
\item If $\bar F$ is $<_D$-increasing \underline{then} it is almost
$D$-smooth.
\item If $I$ is a partial order, $|I|$ an aleph, then there is $J\subseteq
I$, linearly and well ordered by $<_I$ with no upper bound in $I\setminus J$.
\item We are assuming that each $F_t$ is a set; if we consider classes,
intersect them with ${\cal H}(\chi)$ for $\chi$ large enough.
\item If $\bar Y=\langle Y_\alpha: \alpha<\delta\rangle$ is a sequence of
subsets of $A$ and $\alpha<\beta<\delta$ implies $Y_\alpha\subseteq Y_\beta$
and $\cf(\delta)\geq \theta(A)$ \underline{then} $\bar Y$ is eventually
constant i.e.  
$$
(\exists \alpha<\delta)(\forall \beta)(\alpha \leq\beta <\delta \rightarrow
Y_\beta =Y_\alpha)
$$
\item If ${\bar F}=\langle F_{\alpha}: \alpha <\delta\rangle$ is
$\leq_D$-increasing, $A^*=\dom(D)$ and $\theta({\cal P}(\A)/D)\leq
\cf(\delta)$, \underline{then} we can define, in a uniform way, from
$\bar F$ a club 
$C$ of $\delta$ and $Y/D\in {\cal P}(\A)/D$ such that:
\begin{description}
\item[(a)] ${\bar F}\restriction C$ is $<_{D+Y}$--increasing,
\item[(b)] $f_1, f_2 \in\bigcup\limits_{t\in C} F_t \Rightarrow f_1
=_{D+(\A\setminus Y)} f_2$.
\end{description}
\item If $\langle F_\alpha: \alpha<\delta\rangle$ is $<_D$-increasing, $f$ is
a $\leq_D-\lub$ of $\bigcup\limits_{\alpha<\delta} F_\alpha$\\ 
then $\{a\in \A: f(a) \mbox{ a limit ordinal}\}\in D$.
\end{enumerate}
\end{observation}
\medskip

\par \noindent
{\sc Remark:} Not $\theta({\cal P}(\A))\geq \theta({\cal P}(\A)/D)$ as
$|{\cal P}(\A)| {}^*\geq |{\cal P}(\A)/D|$.
\medskip

\Proof (1)--(6), (8) Check.

\par \noindent
(7) For $\alpha<\beta< \delta$ let
$$
\begin{array}{ll}
E^1_{\alpha, \beta}=\{Y/D: & Y\subseteq \A\mbox{ and: } Y=\emptyset\
\mod\ D\mbox{\ \underline{or}\ }Y\in D^+\\ 
\ & \mbox{ and for some }f_1\in F_\alpha\mbox{ and }f_2\in F_\beta\mbox{ we
have  }f_1<_{D+Y} f_2\} 
\end{array}
$$
$$
\begin{array}{ll}
E^0_{\alpha, \beta}=\{Y/D: & Y\subseteq \A, Y=\emptyset \ \mod\ D\mbox{\
\underline{ or }\ } Y\in D^+\\ 
\ & \mbox{ and for every }f_1\in F_\alpha, f_2\in F_\beta\mbox{ we have }f_1
<_{D+Y} f_2\} 
\end{array}
$$
Clearly
\begin{description}
\item[$(*)_1$] $E^1_{\alpha, \beta}\subseteq {\cal P}(\A)/D$, $\emptyset/D\in
E^0_{\alpha, \beta}\subseteq E^1_{\alpha, \beta}$, and $\alpha_1\leq \alpha_2
< \beta_2 \leq \beta_1 \Rightarrow E^\ell_{\alpha_2, \beta_2}\subseteq
E^\ell_{\alpha_1, \beta_1}$ 
\end{description}
and
\begin{description}
\item[$(*)_2$] if $\alpha_1<\alpha_2 < \beta_2< \beta_1$ then $E^1_{\alpha_2,
\beta_2}\subseteq E^0_{\alpha_1, \beta_1}$. 
\end{description}
As $\cf(\delta)\geq \theta({\cal P}(\A)/ D)$, for $\ell\in \{0, 1\}$ and
$\alpha< \delta$ the sequence $\langle E^\ell_{\alpha, \beta}:\beta\in
[\alpha, \delta)\rangle$ is eventually constant (by \ref{obs1:2}(6)), so let
$\gamma^\ell_\alpha\in (\alpha, \delta)$ be minimal such that
$\gamma^\ell_\alpha\leq 
\gamma<\delta \Rightarrow E^\ell_{\alpha, \gamma} = E^\ell_{\alpha,
\gamma^\ell_\alpha}$. 

Let $E^\ell_\alpha= E^\ell_{\alpha, \gamma^\ell_\alpha}$. Clearly
\begin{description}
\item[$(*)_3$] $E^0_\alpha\subseteq E^1_\alpha$,
\item[$(*)_4$] $\ell<2 \ \&\ \alpha< \beta \Rightarrow E^\ell_\alpha\subseteq
E^\ell_\beta\subseteq {\cal P}(\A)/D$, 
\item[$(*)_5$] $\alpha< \beta \Rightarrow E^1_\alpha\subseteq E^0_\beta
\subseteq {\cal P}(\A)/D$.
\end{description}
Applying again \ref{obs1:2}(6), by $(*)_4$, we find that for some $\alpha(*)<
\delta$ we have $\alpha(*)\leq \alpha< \delta \Rightarrow E^0_\alpha =
E^0_{\alpha(*)}$ so by $(*)_3+(*)_5$ we get $\alpha(*)\leq \alpha<\delta
\Rightarrow E^0_\alpha = E^1_\alpha = E^1_{\alpha(*)}$.

Choose by induction on $\varepsilon$ an ordinal $\beta_\varepsilon<\delta$:
for $\varepsilon=0$ let $\beta_\varepsilon= \alpha(*)$, for $\varepsilon=
\zeta+1$ let $\beta_\varepsilon=\max\{\gamma^0_{\beta_\varepsilon},
\gamma^1_{\beta_\varepsilon}\}< \delta$ and for $\varepsilon$ limit 
$\beta_\varepsilon= \bigcup\limits_{\zeta< \varepsilon} \beta_\zeta$.
So for some limit ordinal $\varepsilon(*)$ we have: $\beta_\varepsilon$ is
defined iff $\varepsilon< \varepsilon(*)$, and $\langle\beta_\varepsilon: 
\varepsilon< \varepsilon(*)\rangle$ is strictly increasing with limit $\delta$.

Choose $f^*\in F_{\beta_0}$, $f^*\in F_{\beta_1}$ and let $Y^*=\{x\in A^*:
f^*(x)< f^{**}(x)\}$. Clearly $Y^*\in E^1_{\beta_0, \beta_1}=E^1_{\beta_0} =
E^1_{\alpha(*)} = E^0_{\alpha(*)}$ hence  
$$
f^\prime\in F_{\beta_0}\ \&\ f^{\prime\prime}\in F_{\beta_1} \Rightarrow
f^\prime <_{D+Y^*} f^{\prime\prime}.
$$ 
So clearly $Y^*/D$ does not depend on the choice of $f^*$, $f^{**}$, i.e.
$f^\prime\in F_{\beta_0}\ \&\ f^{\prime\prime}\in F_{\beta_1}\Rightarrow
\{x\in \A: f^\prime(x)<f^{\prime\prime}(x)\}= Y^*\ \mod\ D$.

In fact we can replace $(\beta_0, \beta_1)$ by any $(\beta_{\varepsilon_0},
\beta_{\varepsilon_1})$ with $\varepsilon_0<\varepsilon_1< \varepsilon(*)$. So
clearly the conclusion holds with $Y=\A\setminus Y^*$.
\hfill\qed$_{\ref{obs1:2}}$

\begin{claim}\
\label{claim1:3}
\begin{enumerate}
\item  If $\langle f_\beta:\beta <\delta\rangle$ is such that:
$f_\beta:\A\rightarrow \ord$ and\footnote{see \protect\ref{obs6:1} below} 
$$
\delta\geq\theta^-(\bigcup\limits_{\alpha<\theta(\A)}{\cal P}(\A\times\alpha))
$$
\underline{then}\  for some $\beta<\gamma<\delta$ we have $f_\beta\leq
f_\gamma$.
\item  $[AC_{\lambda,|I|}]$  If $I$ is $(\leq \lambda)^*$-directed, $H$ an
increasing function from $I$ to $J,\  |J| \leq \lambda$ ($I,J$ partial
orders), \underline{then}\  $H$ is constant on a cone of $I$.
\item If $I$ is $(<\theta(J))$-directed, $|I|$ an aleph, $H$ an increasing
function from $I$ to $J$ and $\theta(J)\leq\lambda$ \underline{then}\ $H$ is
constant on a cone of $I$.
\item  $[AC_{\lambda,|I|}, AC_{A,B}]$ If $I\ $ is $(\leq^* \lambda)$-directed, $H$
an increasing function from $I$ to $J={}^AB,\ B$ partially ordered set,
${\hcf}^*_- (B)\leq\lambda$ ($J$ ordered by: $ f\leq g \iff (\forall a\in A)(
f(a) \leq g(a))$), and $|A|\leq \lambda$ then $H$ is constant on a cone.
\item In part (2), instead of $|J|\leq\lambda$, actually $\hcf^*_-
(J)\leq\lambda$ suffices (see Definition \ref{1.1}(10)). 
\item $[AC_{\aleph _0}]$ If $D$ is an $\aleph_1$-complete filter on $\A$ and
${\bar F}=\langle F_{\alpha}: \alpha <\delta\rangle$ is almost $D$-smooth, and
$\delta\geq \theta^-(\bigcup\{{}^\alpha({\cal P}(\A)/D): \alpha<\theta({\cal
P}/\A)/D\mbox{ or } \alpha<\omega\})$ 
\underline{then}\ for some $\beta<\gamma<\delta$, $f_{\beta}\in F_{\beta}$,
$f_{\gamma}\in F_{\gamma}$ we have $f_{\beta}\leq_D f_{\gamma}$. 
\end{enumerate}
\end{claim}

\Proof {\em 1)}\ \ \ Assume not. For each $\beta < \delta$ choose by induction
on $i$, $\alpha(\beta, i)=\alpha_{\beta,i}$ as the first ordinal $\alpha$ such
that: 
\begin{description}
\item[$(*)$] $\alpha<\beta$ and for every $j<i$, $\alpha_{\beta,j}<\alpha$ and
for $a\in\A$ 
$$
[f_{\alpha(\beta,j)}(a)>f_{\beta}(a)\quad\iff\quad f_{\alpha(\beta,j)}(a)>
f_{\alpha}(a)].
$$
\end{description}
So $\alpha_{\beta, i}$ is strictly increasing with $i$ and is $<\beta$, hence 
$\alpha_{\beta,i}$ is defined iff $i<i_{\beta}$ for some $i_{\beta}\leq
\beta$. The sequence $\langle \langle \alpha_{\beta,i}:i<i_{\beta}\rangle:
\beta<\delta \rangle$ exists. Now for $\beta<\delta$, $i<j < i_{\beta}$ let
$u_{\beta,i,j}=:\{a \in\A: f_{\alpha(\beta,i)}(a)> f_{\alpha(\beta,j)}(a)\}$.
For each $a \in\A$ and $\beta<\delta$ we have: 
\begin{quote}
\noindent $v_{\beta}(a)=:\{i<i_{\beta}:f_{\alpha(\beta,i)}(a)>f_{\beta}(a)\}$
is a subset of $i_{\beta}$ and  
$$
\langle f_{\alpha (\beta,i)} (a): i\in v_{\beta}(a)\rangle
$$ 
is a strictly decreasing sequence of ordinals; hence $v_{\beta} (a)$ is finite.
\end{quote}
Now $i_{\beta}=\bigcup\limits_{a \in\A} v_{\beta}(a)$ (as if $j\in i_{\beta}
\setminus \bigcup\limits_{a \in\A} v_{\beta}(a)$ then $\alpha_{\beta,j},
\beta$ are as required). So we have a function $v_{\beta}$ from $\A$ onto
$\{v_{\beta}(a):a \in\A \}$, a set of finite subsets of $i_{\beta}$ whose
union is $i_{\beta}$. Now this set is well ordered of cardinality
$|i_{\beta}|$ (or both are finite), so $i_{\beta}<\theta(\A)$. 

Let $A_{\beta,i}=:\{a\in \A: i\in v_{\beta}(a)\}= \{a\in \A: f_{\alpha(\beta,
i)}(a)> f_\beta(a)\}$; also for no $\beta_1\neq\beta_2<\delta$ do we have
$\langle \alpha_{\beta_1,i} ,A_{\beta_1,i}: i< i_{\beta_1} \rangle = \langle
\alpha_{\beta_2,i}, A_{\beta_2, i}:i<i_{\beta_2}\rangle$: as by symmetry we
may assume $\beta_1<\beta_2$; now $\beta_1$ is a good candidate for being
$\alpha_{\beta_2,i_{\beta_2}}$ so $\alpha_{\beta_2 ,i_{\beta_2}}$ is well
defined and $\leq \beta_1$, contradicting the definition of $i_{\beta_2}$.
Similarly 
$$
\beta_1< \beta_2< \delta, \quad \forall
i<i_{\beta_2}[i<\alpha_{\beta_1}\ \&\ \alpha_{\beta_2, i}=
\alpha_{\beta_1 ,i}
\ \& \ A_{\beta_2,i}= A_{\beta_1, i}]
$$ 
is impossible.
Clearly if $j<i_{\gamma}$, $\beta=\alpha_{\gamma,j}$ then
$\langle\alpha_{\beta,i}: i< i_{\beta}\rangle=\langle\alpha_{\gamma,i} :
i< j\rangle$. 

Now for every $\beta<\delta$ let us define $c_{\beta}=\{\alpha_{\beta,i}:
i<i_{\beta}\}$. So $|c_{\beta}| < \theta(\A)$. Let $X_\beta=\{(i, a):
i<i_\beta\mbox{ and }a\in A_{\beta, i}\}$ so $X_\beta\subseteq i_\beta\times
\A$ and $i_\beta< \theta(\A)$. It is also clear that for $\beta_1\neq\beta_2$ 
($<\delta$), $X_{\beta_1}\neq X_{\beta_2}$. (Just let $f$ be a one-to-one
order preserving function from $c_{\beta_1}$ onto $c_{\beta_2}$, let
$\alpha\in c_{\beta_1}$ be minimal such that $f(\alpha)\neq\alpha$ and we get
a contradiction to the previous paragraph.)

So there is a one-to-one function from $\delta$ into $\bigcup\{{\cal P}(\A
\times\alpha): \alpha < \theta(\A) \}$. But $\delta\geq\theta^-(\bigcup
\limits_{\alpha<\theta(\A)}{\cal P}(\A\times\alpha))$. So we are done. 

{\em 2)}\ \ \ Let for $t\in J$, $I_t =\{s\in I: H(s)=t\}$. So $\{ I_t: t\in
J , I_t\neq\emptyset\}$ is an indexed family of $\leq |J|$ nonempty subsets of
$I$, $|I_t|\leq |I|$, so by assumption there is a choice function $F$
for this family. Let
$A=\{F(I_t): t\in J, I_t\neq\emptyset \}$, so $|A|{\leq}^* |J|$, so as $I$ is
$(\leq\blambda)^*$-directed and $|J|\leq\blambda$ there is $s(*)\in I$
such that $(\forall t\in J)(F(I_t){\leq}_I s(*))$. 

Thus $H(s(*))$ is the largest element in the range of $H$. As $H$ is
increasing it must be constant on the cone $\{s: s(*)\leq_I s\}$. 

{\em 3)}\ \ \ Like part $(2)$ (with $F(I_t)$ being the first member of $I_t$
for some fixed well ordering of $I$).

{\em 4)}\ \ \ By the proof of part $(5)$ below it suffices to prove ${\cf}^*_-
({}^A B\restriction Y, \leq)\leq\lambda$; where $Y=\rang(H)\subseteq
{}^A B$, clearly $Y$ is $(\leq^*\lambda)$-directed (proved as in the
proof of (2)). Also $AC_{\lambda, |Y|}$ holds as $|Y|\leq^* |I|$. Let
for $(a,b)\in A\times B$, 
$$
Y_{(a,b)}=\{h: h\in Y \mbox{ and } B\models b = h(a)\}.
$$ 
Let for each $a\in A$, $B_a =\{b\in B: \mbox{ for some }h\in Y\mbox{ we have
}h(a)=b\}\subseteq B$. 

\underline{Case 1}\ \ for some $a\in A$, $B_a$ is a non empty subset
of $A$ with no last element. $B_a$ has an unbounded subset $B^*$,
$|B^*|\leq^*\lambda$ 
as $Y$ is $(\leq^*\lambda)$-directed, to get contradiction
we need a choice function for $\{Y_{(a,b)}: b\in B^*
\}$; by $AC_{\lambda, |Y|}$ one exists.

\underline{Case 2}\ \ not case $1$.\\
Let $h^*\in {}^A B$ be $h(a)=$ the maximal element of $B_a$. If $h^*$
is in $Y$ then it is the maximal member of $Y$; so
$\cf^*_-({}^AB\restriction Y, \leq)=1\leq \lambda$. If not to get
contradiction it suffices to have a choice
function for $\{Y_{a,h^*(a)} : a\in A\}$ which again exists.  

{\em 5)}\ \ \ Like the proof of part $(2)$. Let $J_1=\rang (H)$ and
$J_2\subseteq J_1$ be unbounded in $J_1$, $|J_2|\leq\lambda$ and choose $F$ as
a choice function for $\{I_t: t\in J_2, I_t\neq\emptyset\}$.

{\em 6)}\ \ \ Let for $\beta<\gamma<\delta$, $Y_{\beta,\gamma}\in {\cal
P}(\A)/D$ be such that : for every $f_{\beta}\in F_{\beta}$ and $f_{\gamma}\in
F_{\gamma}$ we have $Y_{\beta,\gamma} =:\{a\in\A :
f_{\beta}(a)>f_{\gamma}(a)\}/D$.  We can assume toward contradiction that the
conclusion fails, so $\beta<\gamma<\delta \Rightarrow Y_{\beta, \gamma}\neq
\emptyset/D$.  For each $\beta<\delta$ we define by induction on $i$ an
ordinal $\alpha_{\beta,i}=\alpha (\beta,i)$ as the first $\alpha$ such that:
$$
(\forall j<i)(\alpha_{\beta,j}<\alpha<\beta)\mbox{  and  }(\forall
j<i)(Y_{\alpha_{\beta,i},\alpha}= Y_{\alpha_{\beta,i},\beta}).
$$
So clearly $\alpha_{\beta, i}$ increase with $i$ and is ${<}\beta$, hence
for some $i_{\beta}\leq \beta$ we have $\alpha_{\beta,i}$ is defined iff 
$i<i_{\beta}$, and $\langle\langle\alpha_{\beta,i}: i<i_{\beta}\rangle
: \beta<\delta\rangle$ exists. Clearly for $i<j<i_{\beta}$,
$Y_{\alpha_{\beta,i}, \alpha_{\beta,j}}= Y_{\alpha_{\beta,i}, \beta}$.
If for some $\beta$ and $Y$ the set $\{i<i_{\beta}:
Y_{\alpha_{\beta,i},\beta}=Y\}$ is infinite, let $j_n$ the $n$-th member.
By $AC_{\aleph_0}$ we can find $\langle f_{j_n}: n<\omega\rangle$ with
$f_{j_n}\in F_{j_n}$, so
$Y_{\alpha_{\beta,j_n},\alpha_{\beta,j_{n+1}}}=Y$, so 
$$
A_n=\{a\in\A :
f_{\alpha_{\beta,j_n}}(a)> f_{\alpha_{\beta,j_{n+1}}}(a)\}\in
Y_{\alpha_{\beta,j_n},\alpha_{\beta,j_{n+1}}} = Y.
$$ 
Hence
$\bigcap\limits_{n<\omega} A_n\neq\emptyset$ (as $D+A_0$ is a proper
$\aleph_1$-complete filter remembering $Y_{\beta, \gamma} \neq
\emptyset/D$) and we get a strictly decreasing sequence 
of length $\omega$ of ordinals $\langle f_{\alpha_{\beta,j_n}}(a):
n<\omega\rangle$ for $a\in\bigcap\limits_{n<\omega} A_n$, a contradiction. So
$i_\beta <\theta({\cal P}(\A)/D)$. 
(In fact let $b_\beta=\{i<i_\beta: Y_{\alpha_{\beta, i},
\beta}\notin\{Y_{\alpha_{\beta, j}, \beta}: j<i\}\}$, then
$|b_\beta|=|i_\beta|$ or both are finite, and clearly $\otp(b_\beta)<
\theta^-({\cal P}(\A)/D)$, but $\theta^-({\cal P}(\A)/D)$ is an aleph
or finite number, hence $i_n<\omega\ \vee\ i_n<\theta^-({\cal
P}(\A)/D)$). So
$\beta\mapsto \langle Y_{\alpha_{\beta,i}}: i<i_\beta\rangle$ is a
one-to-one function from $\delta$ to $\bigcup\limits_{\alpha <\theta({\cal
P}(\A)/D)} {}^\alpha({\cal P}(\A)/D)$ [why? because if $\beta_1\neq
\beta_2$ have the same image, without loss of generality
$\beta_1<\beta_2$ and $\beta_1$ 
could have served as $\alpha_{\beta_1, i_{\beta_2}}$, contradiction so
the mapping is really one-to-one]. \hfill$\QED_{\ref{claim1:3}}$

\begin{claim}\
 [$DC_{\kappa}$, $\kappa$ an aleph $>|{\cal P}(\A)|$ and: $AC_{{\cal
P}(\A)}$ or $|I|$ is an aleph]  
\label{claim1:4}

Assume: 
\begin{description}
\item[(i)]  $D$ is a filter on the set $\A$
\item[(ii)] $I$ is a $(\leq |{\cal P}(\A)|)^*$-directed
partial order \\
(e.g. an ordinal $\delta, \cf(\delta) \geq\theta({\cal P}(\A))$)
\item[(iii)] $F=\{f_t:t\in I\}$ is $\leq_D$-increasing,
$f_t:\A\rightarrow \ord$ 
\end{description}
\underline{Then} $F$ has a $\leq_D-\eub $
\end{claim}    

\Proof Let $\beta^*=\sup \bigcup\limits_{t\in I}\rang(f_t)$, it is an ordinal, and
let $X=$ ${}^{\A} (\beta^*)$. By $DC_\kappa$ as $|{\cal P}(\A)| <\kappa$, so
without loss of generality $\kappa=|{\cal P}(\A)|^+$, $|{\cal P}(\A)|$
an aleph. 

We can try to choose by induction on $\alpha <\kappa$ functions $g_\alpha
\in X$ such that: 
\begin{description}
\item[(a)]  $(\forall \beta<\alpha)(\neg g_\beta \leq_D
g_\alpha)$
\item[(b)]  $(\forall t\in I)(f_t\leq_D g_\alpha)$
\item[(c)]  $g_0(a)= \sup\{f_t(a) :t\in I\}$
\end{description}

By \ref{claim1:3}(1) 
the construction cannot continue for $\kappa$ steps; so as $DC_{\kappa}$ holds,
for some $\alpha=\alpha(*)<\kappa$ we have $\langle g_\alpha:
\alpha<\alpha(*)\rangle$ and we cannot choose $g_{\alpha(*)}$.
By clause (c) clearly
$\alpha(*)>0$.  Let 
$B_a=:\{g_\beta(a):\beta<\alpha(*)\}$ so $B_a$ is a set of ordinals,
$|B_a|\leq |\alpha(*)|$ (as $|\alpha(*)|$ is an aleph) and $\langle
B_a: a\in \A\rangle$ exists.

Define $H:I\rightarrow \prod\limits_{a\in \A}B_a$ as follows:
\[H(t)(a)=\min(B_a\backslash f_t (a))\ \ \ \mbox{(well defined as
$g_0(a) \in B_a$).}\]
Clearly $|\rang(H)| \leq |{\cal P}(\A)|^{|\A|} =|{\cal P}(\A)|$ (as
$|\A|$ is an infinite aleph hence $|\A\times\A|=|\A|$) and $I \vDash
s\leq t\Rightarrow 
H(g)\leq H(s)\ \mod\ D$.  As $I$ is $(\leq |{\cal
P}(\A)|)^*$-directed, by \ref{claim1:3}(2) if $AC_{{\cal P}(\A)}$ and
\ref{claim1:3}(3) if $|I|$ is an aleph we know that
$H$ is constant on some cone $\{t\in I: s^*\leq
t\}$.  Now consider $H (s^*)$ as a candidate for
$g_{\alpha(*)}$: it satisfies clause (b), clause (c) is
irrelevant so clause (a) necessarily fails, i.e. for some $\beta
<\alpha^*$, $g_\beta{\leq}_D H(s^*)$ so necessarily
$\alpha(*)=\beta+1$, and so $g_\beta$ is as required.
\hfill$\QED_{\ref{claim1:4}}$  

\par \noindent
{\sc Discussion\ \ref{claim1:4}A}: In \ref{claim1:4} the demand
$DC_{\kappa}$ is quite strong, 
implying ${\cal P}(\A)$ is well ordered. Clearly we need slightly less
than $DC_{\kappa}$ - only $DC_{\alpha (*)}$ but $\alpha (*)$ is not
given {\em a priori} so what we need is more than $DC_{|{\cal P}(\A)|}$.
Restricting ourselves to
$\aleph_1$- complete filters we shall do better (see \ref{claim1:5}).

\begin{claim}
\label{claim1:4B}
\begin{description}
\item[(1)] \mbox{$[AC_{|{\cal P}(\A)|,|I|}]$\ \ \ }
If (i), (ii), (iii) of \ref{claim1:4}\  hold, and $g$ is a $\leq_D$-$\lub$ of
$F$ \underline{then}\ 
\begin{description}
\item[(a)] $g$ is a $\leq_D$-$\eub$ of $F$
\end{description}
\item[(1A)] $[AC_{{\cal P}(\A\times\A),I} +AC_{\A}]$ Assume in
addition 
\begin{description}
\item[$(ii)^+$] I is $(\leq |{\cal P}(\A\times\A)|)^*$-directed partial order
\end{description}
\underline{Then} 
\begin{description}
\item[(b)] for some $A\subseteq A^*$ and $t\in I$ we have:
\begin{description}  
\item[$(\alpha)$] $A\in D^+\ \&\ t\leq_I s\Rightarrow f_s\restriction
A=g\restriction A\mod D$, 
\item[$(\beta)$] $a\in A\iff [\cf (g(a))\geq\theta({\cal P}(\A))\
\vee\ g(a)=f_t(a)]$,
\item[$(\gamma)$] if 
$A \not\in D^+$, without loss of generality $A=\emptyset$, (changing $g$) 
\end{description}
\end{description}
\item[(2)] If (i), (ii), (iii) of \ref{claim1:4}\  hold, $|I|$ an aleph, $g$ a
$\leq_D$-$\lub$ of $F$ \underline{then} (a) above holds and if
$AC_{\A}$ also (b) above holds.
\item[(3)] Assume:
\begin{description}
\item[(i)] $D$ is a filter on $\A$
\item[(ii)]  $I$ is a partial order, $\theta({\cal P}(\A)/D)$-directed
\item[(iii)] $\bar F=\langle F_t: t\in I\rangle$ is $<_D$-increasing
\item[(iv)] $|I|$ is an aleph
\item[(v)] $g$ a $<_D$-$\lub$ of $F$ 
\end{description}
\underline{Then} (a) above holds. If $AC_{\A}$ then also (b) holds.\\
Even waving $AC_{\A}$ (but assuming (i)-(v)) we have 
\begin{description}
\item[(b)$'$] for some $A\subseteq {\A}$ and $t(*)$
\begin{description}
\item[$(\alpha)$]  if $A\in D^+$ and $t(*)\leq t\in I$ \mbox{ and }
$f\in F_t$ then\\ 
$f\restriction A=g\restriction A\mod D$
\item[$(\beta)$]  assume $(ii)^+$; for no $C\subseteq \ord$ such that 
$|C|<\theta({\cal P}(A))$ and $\{a\in {\A}\setminus A : g(a)=\sup
(C\cap g(a))\}\neq\emptyset\mod D$. 
\end{description}
\end{description}
\end{description}
\end{claim}

\par \noindent
{\sc Remark:} In \ref{claim1:4B}(1A) instead $\A\times \A$ we can use
$\A\times \zeta$ for every $\zeta<\theta({\cal P}(\A)/D)$
\medskip

\Proof
{\em 1)} (a) Let $f\in {}^{\A}\ord$, $f<\max\{g, \1_{\A}\}$. We define
a function $H : I\longrightarrow {\cal P}(A)$ by $H(t)=\{a\in\A:
f(a)\geq \max\{f_t (a), \1\}\}$. By \ref{claim1:3}(2) we know that $H$ is
constant on a 
cone of $I$, more exactly $H'$, $H'(t)=H(t)/D$ is increasing and is
constant by the proof of \ref{claim1:3}(2). Let $t\in I$ be in the cone. If
$H(t)=\emptyset \mod D$ we are done, otherwise define $f^\ast$ to be
equal to $g$ on $A^\ast\setminus H(t)$ and equal to $f$ on $H(t)$.
Now $f^\ast$ contradicts the fact that $g$ is $a<_D -\lub$. Now check.

{\em 1 A)} (b) for $A\subseteq\A$ let $I_A =\{t\in I: A=\{a\in\A: f_t
(a)=g(a)\}\}$ if not empty, and $I$ otherwise. Let us apply $AC_{{\cal
P}(\A), I}$ to $\{I_A: A\in {\cal P}(\A)\}$ getting $\{t_A: A\in {\cal
P}(A)\}$. Let $t(*)\in I$ be an upper bound of $\{t_A: A\in {\cal
P}(A)\}$ (exists by assumption $(ii)$). Let 
$$
A_0=: \{a\in A^*: f_{t(*)}(a)\neq g(a)\}
$$
and 
$$
A_1=\{a\in A_0: g(a)\mbox{ is a limit ordinal}\}
$$
and lastly 
$$
A'=:\{a\in A_1: \cf(g(a))<\theta({\cal P}(\A))\}.
$$
Now clearly $t(*)\leq_I s\in I \Rightarrow f_s< g\ \mod\ (D+A_0)$ and
$A_0=A_1\ \mod\ D$. If $A'=\emptyset\ \mod\ D$
we are done, so
assume $A'\in D^+$. By $AC_{\A}$ we can find $\langle C_a:
a\in\A\rangle$, such that for $a\in A'$ we have: ${C}_a$ an unbounded
subset of $g(a)$ of order type 
$\cf[g(a)]$. Let for $a\in {\A\setminus A'}$, $C_a=\{0\}$, and for
$h\in\prod\limits_{a\in\A} C_a$ let $I^*_h =:\{t\in I: h<_D \max\{f_t,
\1_{\A}\}\}$.
Now the assumption of \ref{claim1:4B}(1A) implies that of
\ref{claim1:4B}(1), hence clause (a) holds, and hence $I^*_h$ is not empty. 
As we have $AC_{\A}$ clearly $|\A|$ is an aleph, and $\theta({\cal
P}(\A))$ is an aleph with cofinality  $>|\A|$, hence $\{\otp(C_a):
a\in \A\}$ is a bounded subset of $\theta({\cal P}(\A))$, so we can
find a function $h^*$ from ${\cal P}(\A)$ onto
$\zeta^*=\sup\{\otp(C_a): a\in A^*\}$. Let $h_a$ be the unique
one-to-one order preserving function from $\otp(C_a)$ onto $C_a$, so 
$\langle h_a: a\in \A\rangle$ exists. For $h\in \prod\limits_{a\in \A}C_a$
let 
$$
X_h=:\{(a, b): {}^a\in \A\mbox{ and }h^{-1}_a\circ
h(a)=h^*(b)\}\subseteq \A\times \A.
$$
Clearly the mapping $h\mapsto X_h$ is one-to-one from
$\prod\limits_{a\in \A} C_a$ into ${\cal P}(\A\times \A)$. Hence by
$AC_{{\cal P}(\A), I}$ 
$$
|\prod\limits_{a\in\A}
C_a|\leq |{\cal P}(\A\times\A)|=|{\cal P}(\A)|,
$$ 
so we can find a function $G$,
$\dom(G)=\prod\limits_{a\in\A} C_a$, $G(h)\in I^*_h$. By $(ii)$ we get that
$\rang(G)$ has an $\leq_I$-upper bound $t(**)$, and we can get a
contradiction. 

{\em 2)} Follows by (3) (with $F_t =:\{ f_t\}$)

{\em 3)} (a) Let $<^*$ be a well ordering of $I$. Let $f\in {}^{(\A)}
\ord$, $f<\max\{g,\1_{\A}\}$. For every $A\in D^+$ let $t_{A/D}$ be the
$<^*$-first $t\in I$ such that for some $f_t \in F_t$ we have
$\{a\in\A : f_t (a)\geq f(a)\} = A\mod D$. Now $A\mapsto t_{A/D}$ is
a mapping (maybe partial) from ${\cal P}(\A)/D$ into $I$ and $|I|$ is
an aleph, so 
$$
|\{t_{A/D}: A\in D^+\}|<\theta({\cal P}(\A)/D)\ (\leq\theta({\cal
P}(\A))).
$$
But $I$ is $\theta({\cal P}(\A)/D)$-directed so there is $t(*)\in I$
such that 
$$
(\forall A\in D^+)(t_{A/D}\leq_I t(*).
$$ 
Any $f_{t(*)} \in
F_{t(*)}$ satisfies $f<\max\{f_{t(*)}, \1_{\A}\}$.

(b) + (b)$^\prime$
For each $A\in D^+$, let $t_{A/D}$ be the $<^*$-first $t\in I$ such
that 
$$
A\subseteq \{a\in \A: f_t(a)=g(a)\} \mod D
$$ 
if there is one.
So clearly 
$$
\{t_{A/D}: A\in D^+\mbox{ and }t_{A/D}\mbox{ is well defined}\}
$$ 
is a set of $<^*$-order type
$<\theta({\cal P}(\A)/D)$, hence there is $t(*)\in I$ such that 
$$
A\in D^+ \ \&\ t_{A/D}\mbox{ well defined }\Rightarrow t_{A/D}\leq_I
t(*).
$$ 
Let 
$$
A_0=:\{ a\in \A: f_{t(*)}(a)\neq g(a)\}
$$
$$
A_1=:\{ a\in A_0: f(a)\mbox{ a limit ordinal}\}
$$
and
$$
A'=:\{ a\in A_1: \cf(g(a))<\theta({\cal P}(\A))\},
$$
as in the proof of \ref{claim1:4B}(1A) we have for $A=A_0$, clause
$(\alpha)$ of (b) and of (b)$'$ holds, and $A_1=A_0\ \mod\ D$, and if
$A'=\emptyset\ \mod\ D$ we are done. As for clause (b)$(\beta)$
by $AC_{\A}$, if it fails then $(\beta)$ of $(b)'$ fails, so it suffices
to prove $(b)'$. Now clause $(\alpha)$ holds, and we prove $(\beta)$
as in \ref{claim1:4B}\ (1A). \hfill$\QED_{\ref{claim1:4B}}$

\begin{claim}
\label{claim1:5}
$[DC_{{\aleph_0}} +AC_{{\cal P}(A^*)}]$\ \ \ Assume:
\begin{description}
\item[(i)] $D$ is an $\aleph_1$-complete filter on a set
$A^*$,
\item[(ii)] $I$ is a $(\leq |{\cal P}(A^*)|)^*$-directed partial
order, 
\item[(iii)] $F=\{f_t:t\in I\}$ is $\leq_D$-increasing,
$f_t:A^*\rightarrow \ord$. 
\end{description}
\underline{Then} $F$ has a $\leq_D-\eub$
\end{claim}

\noindent{\sc Remark \ref{claim1:5}A}:\ \ \ \ \
\begin{enumerate}
\item Given ${\bar F}=\langle F_t: t\in I\rangle$ 
is $<_D$-increasing (i.e. $[I\models s<t\ \&\ f\in I_s\ \&\ g\in I_t
\Rightarrow f <_D g]$) 
we can use
$J=(\bigcup\limits_{t\in I} F_t,<_D)$ as the partial order instead $I$ with
$f_g = g$ for $g\in\bigcup\limits_{t\in I} F_t$, so claim \ref{claim1:5}
applies to ${\bar F}=\langle F_t: t\in I\rangle$ too; similarly for
\ref{claim1:4} . Note that if $I$ is $(\leq\lambda)^*$-directed so is
$J$. Also \ref{claim1:4B} $(1)$ (in the proof replace $H$ by $H':
H'(t)=H(t)/D$), concerning \ref{claim1:4B}\ $(1A)$ check, and lastly for
\ref{claim1:4B}\ $(2)$ see \ref{claim1:4B}\ $(3)$)
\item If we want to demand only $AC_{{\cal P}(\A),\lambda}$, then
$\lambda=\prod\limits_{\alpha\in\A}\beta_a$ is large enough. 
\item Note $AC_{{\cal P}(\A)}$ can be replaced by $AC_{{\cal
P}(\A)/D}$. 
\end{enumerate}

\Proof By \ref{claim1:4B}$(1)$ it suffices to find a $\leq_D$-lub.
Let for $a\in 
\A$, 
$$
\beta_a=:\sup\{f_t(a): t\in I\}+1 \in \ord.
$$ 
For every $A\in D^+$, there is no decreasing $\omega$-sequence in
$(\prod\limits_{a\in A^*} \beta_a, <_{D+A} )$; hence by $DC_{\aleph_0}$
there is a function $g \in \prod\limits_{a \in A^*}\beta_a$ satisfying:
\begin{description}
\item[$(i)_A$]\ \ $(\forall t\in I) [f_t \leq_{D+A} g]$
\item[$(ii)_A$]\ \ if $g'$ satisfies (i) then $\neg(g'<_{{D+A}}g)$.
\end{description}
By $AC_{{\cal P}(A^*)}$ we can find  $\langle g_A: A \in
D^+\rangle$ with $g_A \in \prod\limits_{a\in A} \beta_a$ that satisfies
$(i)_A+(ii)_A$.  
Let $B_a=\{g_A(a): A\in D^+\} \cup \{\beta_a\}$ so $|B_a|\leq
|D^+|\leq|{\cal P}(A^*)|$ 

Let $H: I\longrightarrow B^*=\prod\{B_a:a\in A^*\}$ be
$H(t)(a)=\min\big(B_a\backslash f_t (a)\big)$. Clearly $H$ is an order
preserving mapping from $(I, \leq_I)$ to $(B^* /D,\leq_D)$. Also
$$
|B^*/D|\leq^* |B^*|\leq^*\prod\limits_{a\in\A} |D^+|=|D^+|^{|\A|}\leq |{\cal
P}(\A)|^{|\A|}.
$$ 
But as $AC_{{\cal P}(\A)}$ holds  clearly $\A$ is well ordered, hence
$|\A|\times |\A|=|\A|$, hence $|{\cal P}(\A)|^{|\A|}=|{\cal P}(\A)|$.

By \ref{claim1:3}(2) we know that $H$ is constant on a cone, say
$\{t:s(*)\leq_I t\}$.  

Assume $g^*=: H(s(*))$ is not a $\leq_D$-\lub\ , so some $g'$ exemplifies
it.  Let $A=\{a\in A^*: g'(a)<g^*(a)\}$, 
so $A\in D^+$ and $g_A$ is well defined.

By the choice of $g_A$ and $g'$ clearly $\neg(g'<_{D+A} g_A)$ hence
\[A_1=\{a\in A^*:g'(a)\geq g_A (a)\}\neq\emptyset\mod D+A\]

So $A_2=A_1\cap A\in D^+$; and $g_A\leq_{D+A_2} g'$.
But by the choice of $g^*= H(s(*))$ (as $(\forall a)(g_A(a)\in B_a)$)
clearly $g^*\leq_{D+A} g_A$. Together we get $g^*\leq_{D+A_2} g'$, but
this contradicts the definition of $A$ as $A_2\subseteq A$.
\hfill$\QED_{\ref{claim1:5}}$    

\begin{claim}
\label{claim1:6}
$[DC_{\aleph_0}]$\ \ \ Consider the conditions:
\begin{description}
\item[(i)] $D$ is an $\aleph_1$-complete filter on $A^*$,
\item[(ii)] $I$ is a $(\leq |{\cal P}(\A)/D|)^*$-directed partial order,
\item[(iii)] $\bar{F}=\langle F_t : t\in I\rangle$ is
${\leq}_D$-increasing, so $f:A^*\rightarrow\ord$ for
$f\in\bigcup\limits_{t\in I} F_t$,
\item[(iv)] $g:A^*\rightarrow \ord$ and $(\forall t\in I)(\forall f\in
F_t)(f\leq_D g)$, 
\item[(v)] no $g' <_D g$ satisfies (iv).
\end{description}
\underline{Then} the following implications hold:
\begin{description}
\item[$(1)$] If (i), (iii) \underline{then} some $g$ satisfies (iv).
\item[$(2)$] If (i), (iii), (iv) \underline{then} some $g^*\leq g$
satisfies clauses (iv), (v); if for $g$ clauses (i), (iii), (iv),
$\neg$(v) hold, we can ask also $g^* <_D g$. 
\item[$(3)$] $[AC_{{\cal P}(\A)/D,I}]$\ \ \ If (i)-(v), and $\bar F$ is
$D$--smooth \underline{then} there is $D^*$, such that 

\ \ \ \ (a) $D^*$ an $\aleph_1$-complete filter extending
$D$,

\ \ \ \ (b) $g$ is $<_{D^*}-\lub$ of $\bar{F}$.
\end{description}
\item[$(4)$] In (3) we can omit smoothness. 
\item[$(5)$] If (i)-(v) and $|I|$ is an aleph then for some $D^\ast$ clauses
(a), (b) of (3) hold (i.e. we can drop $AC$ in part (3)).
\end{claim}

\Proof  $(1)$, $(2)$ are contained in the proof of \ref{claim1:5}

$(3)$  Let $D^*=\{A\subseteq A^*:$ for $D+(A^*\backslash A)$,
 clause (v) fails$\}$. Clearly $D^*$ is $\aleph_1$-complete filter and
$\emptyset\notin D^*$. 
Clearly $\bar F$ is $\leq_{D^*}$-increasing and $g$ is a $\leq_{D^*}$-upper
bound of $\bar F$. If $g$ is not a $\leq_{D^*}$-$\lub$, there is $g'\in
{}^{(\A)} \ord$, a $\leq_{D^*}$-upper bound of $\bar F$ such that $\neg
[g\leq_{D^*} g']$, so 
$$
A'=:\{a\in\A: g'(a)<g(a)\}\neq\emptyset\mod
D^*.
$$ 
For each $t\in I$, for every $f_t \in F_t$ let 
$$
Y_t=:\{a\in\A:
f_t(a)>g'(a)\}/D\in D^*/D
$$ 
(well defined as $\bar F$ is
$D$-smooth), and $t\mapsto Y_t$ is an 
increasing function from $I$ to $(D^*/D,\leq_D)$, by
\ref{claim1:3}(2) we get a contradiction.

$(4)$ Use \ref{obs?}(3) below to regain smoothness.

$(5)$ Left to the reader (just use \ref{claim1:3}(3) instead
\ref{claim1:3}(2)). \hfill$\QED_{\ref{claim1:6}}$ 
\medskip 

\noindent We have used
\medskip

\begin{observation}
\label{obs?}
\begin{enumerate}
\item \ \mbox{[$AC_{{\cal P}(\A)/D,I}$ or $|I|$ is an aleph]}\ \ \ \ \ 
If $\bar F=\langle F_t : t\in I\rangle$ is
$<_D$-increasing (where $f_t :\A =:\dom(D)\longrightarrow\ord$), 
$I$ is $(\leq|{\cal P}(A)/D|)^*$-directed \underline{then}\ for every
$g: \A\longrightarrow\ord$ for some $t^\ast\in I$ and $A\subseteq
\A$, for all $f\in\bigcup\limits_{t^*\leq t\in I}F_t$
\[\{a\in\A: f(a)>g(a)\}=A\mod D.\]
Similarly for $f(a)=g(a)$, $f(a)<g(a)$.
\item In (1) in the case $|I|$ is an aleph, ``$I$ is $(<\theta({\cal
P}(A)/D))^*$- directed'' suffices.
\item \mbox{$[AC_{\aleph_0}]$} Assume $\bar F=\langle F_t: t\in
I\rangle$ is $<_D$-increasing.
Let $I$ be a partial order
$(\leq\aleph_0)^*$- directed, 
$$I_0=\{\langle t_n: n<\omega\rangle :
t_n\in I \mbox{  and  } t_n\leq_I t_{n+1} \mbox{  for  } n<\omega\},$$
$${\bar t}^0\leq_{I_0} {\bar t}^1\iff [(\forall n<\omega)(\exists
m<\omega)(t_n^0\leq_I t_m^1)],$$
$${\bar t}^0 <_I {\bar t}^1\mbox{ if }{\bar
t}^0\leq_I {\bar t}^1 \& \neg ({\bar t}^1\leq_I {\bar t}^0).$$
Let $F^*_{\bar t} =\{\sup_{n<\omega} f_n : \langle f_n :
n<\omega\rangle\in\prod\limits_{n<\omega} F_{t_n}\}$ and ${\bar F}^*=\langle
F^*_{\bar t} : {\bar t}\in I\rangle$ where of course
$(\sup\limits_{n<\omega} f_n)(a)=\sup\limits_{n<\omega} f_n(a)$.
\underline{Then}
\begin{description}
\item[(a)] ${\bar F}^*$ is $\leq_D$-increasing and smooth
\item[(b)] $g$ is a $\leq_D$-$\lub$ of $\bar F$ iff $g$ is a
$\leq_D$-$\lub$ of ${\bar F}^*$ 
\item[(c)] $g$ is a $\leq_D$-$\eub$ of $\bar F$ iff $g$ is a
$\leq_D$-$\eub$ of ${\bar F}^*$
\item[(d)] $I$, $J$ are equi-directed, what means: there is an
embedding $H$ of $I$ into $J$ (a partial order) such that its range is
a cofinal subset of $J$ [use $H(t)=\langle t: n<\omega\rangle$].
\end{description}
\item \mbox{$[AC_{\aleph_0}]$} Assume $I=\delta$ where $\delta$ is a
limit ordinal of cofinality $>\aleph_0$, and let $\bar F=\langle
F_\alpha: \alpha<\delta\rangle$ be $<_D$-increasing (where each
$F_\alpha\subseteq {}^{(\A)}\ord$ is not empty).
Let $J=\{\alpha<\delta: \cf(\alpha)=\aleph_0\}$, and for $\alpha\in J$
let 
$$
\begin{array}{ll}
F^*=:\{\sup\limits_{n} f_n: & \mbox{ for some strictly increasing
sequence }\langle \alpha_n: n<\omega\rangle\\
\ & \mbox{ of ordinals }<\alpha\mbox{ we have
}\alpha=\bigcup\limits_{n<\omega} \alpha_n, \mbox{ and }f_n\in
F_{\alpha_n}\}.
\end{array}
$$
Then $\bar F^* = \langle F^*_\alpha: \alpha\in J\rangle$ satisfies
clauses (a)-(d) from part (3), and of course $J$ is well ordering.
\end{enumerate}
\end{observation}

\Proof Check.

\section{$hpp$}

\begin{definition}\ 
\label{definition2:1}
\begin{enumerate}
\item Let $\Gamma$ be a set of filters on
$A^*=A^*_\Gamma=\dom(\Gamma)$. For an ordinal $\alpha$,  we let
$\hpp_\Gamma(\alpha)$ be  
the supremum of the ordinals $\beta+1$ for which there
is a witness $(\bar F,D)$, which means:
\begin{description}
\item[(i)] $D\in \Gamma$,
\item[(ii)] $F=\langle F_\gamma: \gamma<\beta\rangle$, with
$F_\gamma\neq \emptyset$,
\item[(iii)] $\bigcup\limits_{\gamma<\beta} F_\gamma\subseteq {}^{(\A)}\alpha$,
\item[(iv)] $\bar F$ is $<_D$-increasing.
\end{description}
\item $\ehpp_\Gamma(\alpha)$ is defined similarly, but each $F_\gamma$
a singleton; the definition of $\shpp_\Gamma(\alpha)$ is similar too,
but $\bar F$ is smooth.
We can
replace $\alpha$ by $f\in {}^{(\A)}\ord$ in all these cases here (i.e.
$F\subseteq \prod\limits_{a\in\A}f(a)$). If $F_\alpha=\{ f_\alpha\}$
we may write $F=\langle f_\alpha:\alpha<\beta\rangle$ instead of $\bar{F}$.
\item If $\Gamma=\{D\}$ then we write $D$ instead of $\Gamma$. We say
that $\Gamma$ is $\aleph_1$-complete if each $D\in\Gamma$ is
$\aleph_1$-complete, and similarly for other properties.
\item We define $\pp_\Gamma(\alpha)$ as in (1) but add to {\bf
(i)}-{\bf (iv)}
also
\begin{description}
\item[(v)] there exists $\langle\alpha_a: a\in\dom(D)\rangle$ such
that: 
\begin{description}
\item[$(\alpha)$] \ \ \ $\theta(\A)\leq\alpha_a <\alpha$,
\item[$(\beta)$] \ \ \ $\cf_D(\langle\alpha_a: a\in\A\rangle)=\alpha$
(see below), each $\alpha_a$ a limit ordinal,
\item[$(\gamma)$] 
\ \ \ for every $g\in\prod\limits_{a\in\A}\alpha_a$
there is $f\in \bigcup\limits_{\gamma<\beta} F_\gamma$ such that
$g\leq_D f$. 
\end{description}
\end{description}
\item $\epp_\Gamma(\alpha)$, $\spp_\Gamma(\alpha)$ are defined similarly.
\item For a filter $D$ with the domain $\A$,
$$
\begin{array}{l}
\cf_D(\langle\alpha_a: a\in\A\rangle ) \\
\quad =\inf\{\otp(C): C\subseteq
\bigcup\limits_{a\in\A} \alpha_a\mbox{ and }\{a\in\A:
\alpha_a=\sup(C\cap\alpha_a)\}\in D\},\\ 
\quad\lim\sup\limits_D(\langle\alpha_a : a\in\A\rangle )=\min\{\alpha
: \{a\in\A : \alpha_a\leq\alpha\}\in D\}.
\end{array}
$$
\end{enumerate}
\end{definition}

\par \noindent
{\sc Remark:} 
1) Note that in \ref{definition2:1}(4) clause $(\beta)$ (i.e.
$\cf_D(\langle \alpha_a: a\in \A\rangle)= \alpha$) is a replacement to
``$\langle \alpha_a: a\in \A\rangle$ is a sequence of limit ordinals
with $\tlim \langle \cf(\alpha_a): a\in \A\rangle =\alpha$''.

2) Note $\pp$ stands for pseudo power, $h$ for hereditary, $s$ for
smooth, $e$ for element (rather than set).
\medskip

\begin{observation}
\label{observation2:2} 
\begin{enumerate}
\item $\hpp_\Gamma(\alpha)$, $\ehpp_\Gamma(\alpha)$ and
$\shpp_\Gamma(\alpha)$ increase with $\alpha$ and $\Gamma$; 
and  $\pp_\Gamma(\alpha)$, $\epp_\Gamma(\alpha)$ and
$\spp_\Gamma(\alpha)$ increase with $\Gamma$. 
\item $\ehpp_\Gamma(\alpha)\leq \shpp_\Gamma (\alpha)
\leq\hpp_\Gamma(\alpha) \leq 
\theta({}^{A^*}\alpha)$ and $\epp_\Gamma(\alpha) \leq
\pp_\Gamma(\alpha) \leq \theta({}^{(\A)}\alpha)$, and 
$\epp_\Gamma(\alpha)\leq\ehpp_\Gamma(\alpha)$ and $\pp_\Gamma(\alpha)
\leq \hpp_\Gamma(\alpha)$.
\item ${\it x}_\Gamma(\alpha)=\sup\{x_D(\alpha): D\in\Gamma\}$ for
${\it x}\in\{\pp, \epp, \hpp, \ehpp, \spp, \shpp\}$.
\item ${\it x}_D(\alpha)=\sup\{{\it x}_D(\beta): \theta(\dom(D))
\leq\beta<\alpha\}$ for ${\it x}\in\{\hpp,\ehpp, \shpp\}$ if
$\theta(\dom(D))\leq \cf(\alpha)$.
\item If $a\neq b\Rightarrow\alpha_a\neq\alpha_b$ then
$\cf_D(\langle\alpha_a: a\in\A\rangle)\geq \min_{A\in D}|A|$.   
\end{enumerate}
\end{observation}

\par \noindent
{\sc Remark \ref{observation2:2}A:}

If $\hpp_\Gamma(\alpha)>\beta$, and $\cf(\beta)\geq \theta({\cal
P}(\dom(D))$ and relevant criterion for existence of $<_D-\eub$ holds  
for $D\in \Gamma$ \underline{then} for some $\alpha'\leq
\alpha$ we have $\pp_\Gamma(\alpha')>\beta$. See below.
\medskip

\Proof
Easy, e.g.\\
5) let $C\subseteq \bigcup\limits_{a\in \A} \alpha_a$ exemplify
$\cf_D(\langle \alpha_a: a\in \A\rangle)=\otp(C)$. Let 
$$
A'=\{a\in \A: C_a\neq a\}\quad \mbox{ where }
$$ 
$$
C_a=\{\beta\in C: \beta\leq \alpha_a\mbox{ but for no }b\in \A\mbox{
do we have }\beta\leq \alpha_b< \alpha_a\}.
$$ 
So there is a one-to-one function from $A'$ into $C$ and
$$
\otp\{\alpha_a: a\in A'\}< \otp\{\alpha_a: a\in A'\}+\omega,
$$
 so we
can finish. 

\begin{definition}\
\label{definition2:4}
\begin{enumerate}
\item $I$ has the true cofinality
$\delta$, $\delta$ an ordinal {\em if} there is a cofinal $J\subseteq I$ and
a function $h$ from $J$ onto $\delta$ such that $h(f_1)<h(f_2)\Rightarrow
f_1 <_I f_2$.
\item $I$ has strict cofinality $\delta$, $\delta$ an ordinal
{\em if} the function $h$ above is one-to-one.
\end{enumerate}
\end{definition}

\begin{claim}
\label{claim2:5}
\begin{enumerate}
\item If $I$ has the strict cofinality $\delta$ {\em then} $I$ has
the true cofinality $\delta$ and $\cf(I)\leq\delta$.
\item If $I$ has the true (strict) cofinality $\delta$ {\em then} $I$
has the true (strict) cofinality $\cf(\delta)$ (which is regular).
\item If $\bar{F}=\langle F_\alpha:\alpha<\delta\rangle$,
$\bigcup\limits_{\alpha<\delta} F_\alpha\subseteq {}^{(\A)}\ord$, $\bar{F}$ is
$<_D$-increasing and $\bar{F}$ has $<_D$-\eub\ $\langle\alpha_a:
a\in\A\rangle$ \underline{then}  $\prod\limits_{a\in\A} \alpha_a/D$ has the true
cofinality $\cf(\delta)$.
\item If $I$ has the true cofinality $\delta_1$ and $\delta_2$ then
$\cf(\delta_1)=\cf(\delta_2)$. 
\end{enumerate}
\end{claim}

\noindent{\sc Remark \ref{claim2:5}A:}\ \ \ We did not say ``$I$ has
the true cofinality $\lambda$ $\Rightarrow$ $\cf(I)\leq\lambda$''. 
But it is true that: $I$ has true cofinality $\lambda$ implies
$\cf^*(I)\leq^* \lambda$.

\begin{claim}
Let $\Gamma$, $A^*$ be as in \ref{definition2:1}(1).
\label{claim2:3} 
\begin{enumerate}
\item \ $[ |\A|$ an aleph, $DC_{|\A|^+}]$\ \ Assume $\delta <
\hpp_\Gamma(\alpha), 
\delta$ a cardinal and $\cf(\delta)\geq\theta({\cal P}(A^*))$. Then for some
$D\in\Gamma$ and $\bar{\alpha}=\langle\alpha_a:a \in A^*\rangle$,
with each $\alpha_a$ being a limit ordinal $\leq \alpha$ we have,
$\prod\limits_{a\in {A^*}}\alpha_a/ D$ has the true cofinality $\delta$.

\item \ {\em [$AC_{\A}$]} \ \ If $\prod\limits_{a\in
{A^*}} \alpha_a/ D$ has the true cofinality $\mu$ \underline{then}
$\prod\limits_{a\in {A^*}} \cf(\alpha_a)/D$ also has the true cofinality $\mu$.  
\item \ $[|A^*|$ an aleph and $DC_{|\A|^+}]$\ \ Assume $\omega\leq\alpha <
\aleph_\alpha$ and $\aleph_\alpha >\theta({\cal P}(\A))$. If
$\aleph_\gamma\leq\hpp_D(\aleph_\alpha)$ \underline{then} 
${}^{(A^*)}\alpha$ can be mapped onto $\gamma$. Similarly for
$\aleph_{\alpha_0 +\gamma}\leq\pp_D (\aleph_{\alpha_0 +\alpha})$. 
\item  \ Similar claims hold for $\pp$, $\shpp$, $\spp$, $\ehpp$ and $\epp$.
\end{enumerate}
\end{claim}

\Proof  $1)$ By the definition of $\hpp_\Gamma(\alpha)$ we can
find $D\in \Gamma$ and a $<_D$-increasing sequence $\bar F=\langle
F_\alpha: \alpha< \delta\rangle$, with $F_\alpha\subseteq
{}^{(\A)}\alpha$ non empty. \\
If $g\in {}^{(\A)}(\alpha+1)$ is a
$<_D-\eub$ of $\bar F$ then by \ref{claim2:5}(3), $\alpha_a=:g(a)$
for $a\in \A$ are as required. Now such $g$ exists by \ref{claim1:4}
which is applicable by 1.7A(1). 
 
$2)$\ \  Let ${\bar F} =\langle F_{\alpha} : \alpha <\mu\rangle$
exemplify that $\prod\limits_{a\in\A} \alpha_a$ has true cofinality $\mu$.
By $AC_{\A}$ we can find $\langle C_a : a\in\A\rangle\in\v$, $C_a$ a
club of $\alpha_a$ of order type $\cf(\alpha_a)$. For $f\in
\prod\limits_{a\in\A} \alpha_a$ let $f^{\otimes}\in\prod\limits_{a\in\A} C_a$ be
$f^{\otimes} (a)=\min (C_a\setminus f(a))$, so $f\leq f^{\otimes}
\in\prod\limits_{a\in\A} C_a$ and $f_1\leq_D f_2\Rightarrow f^{\otimes}_1\leq_D
f^{\otimes}_2$. Now apply \ref{obs1:2}(7) to $\langle\{f^{\otimes}
: f\in F_{\alpha}\}: \alpha <\mu\rangle$.

$3)$\ \ It is enough to map ${}^{(\A)} \alpha$ onto
$(\gamma -1)\setminus\alpha$. We
define a function $H$ from ${}^{(\A)} \alpha$ into
$\gamma\setminus\alpha$. For $f\in {}^{(\A)} \alpha$, if
$(\prod\limits_{a\in\A} \aleph_{f(a)}, <_D)$ has true cofinality
$\aleph_{\beta+1}\geq\aleph_\gamma$ then let $H(f)=\beta$. By
$AC_{\A}$ and \ref{claim2:3}\ (2) it is enough for every $\beta$ such
that $\aleph_\alpha <\aleph_{\beta +1} <\aleph_\gamma$ to find
$f\in\ {}^{(\A)} (\aleph_{\alpha})$ such that $(\prod\limits_{a\in\A} f(a),
<_D)$ has the true cofinality $\aleph_{\beta +1}$. Now we use part
(1) of \ref{claim2:3}.

$4)$\ \ By \ref{observation2:2} or repeating the
proofs.\hfill$\QED_{\ref{claim2:3}}$

\begin{claim}
\label{claim2:6} $[DC_{\aleph_0}+ (\forall D\in \Gamma)AC_{{\cal
P}(\dom(D))}]$\ \ \ In \ref{claim2:3}(1), if
$\Gamma$ is $\aleph_1$- complete \underline{then} the conclusion
holds.
\end{claim}

\Proof\ We can use \ref{claim1:5} 
instead of \ref{claim1:4}.

\begin{definition}
\label{definition2:7} 
$$
\begin{array}{ll}
\pcf_\Gamma \{\alpha_a: a \in
A^*\}=: \{\lambda: & \lambda \mbox{ is the true cofinality of }
\prod\alpha_a/D\\ 
\ & \mbox{ for some filter }D \mbox{ on }A^*\mbox{ which belongs to
}\Gamma\}.
\end{array}
$$
\end{definition}

\noindent{\sc Remark \ref{definition2:7}A:}\ \ We could phrase
\ref{claim2:3} for $\pcf_\Gamma$.

\begin{claim}
\label{claim2:8}
\begin{enumerate}
\item If $\lambda$ is an aleph, $\mu<\lambda<\epp_D(\mu)$ and
$\theta(\dom(D))<\lambda$ \underline{then} $\lambda$ is not measurable
\item If $\lambda$ is an aleph, $\mu<\lambda<\epp_D(\lambda)$,\\
\underline{then} there is no \mbox{$(\leq|\dom(D)|+\mu)$-complete}
uniform ultrafilter on $\lambda$.
\end{enumerate}
\end{claim}

\Proof {\em 1.}\ \ Let $(\bar F,D)$ witness
$\lambda<\epp_D(\mu)$, so $f\in F\Rightarrow \dom (f) = A^*$, where
$A^*=:\dom(D)$; so let $F=\{f_\alpha:\alpha<\lambda\}$, where
$F_\alpha=\{f_\alpha\}$. Assume $D^*$ is a $\lambda$-complete
ultrafilter on $\lambda$.  For each $a\in A^*$, we have a function
$g_a:\lambda\rightarrow\mu$, defined by $g_a(\alpha)=f_\alpha(a)$. So
$g_a$ is a one-to-one map 
$$
\mbox{ from }s'=:\{\alpha<\lambda:
(\forall\beta < \alpha)(g_a(\beta)\neq g_a(\alpha))\}\mbox{ onto }\rang
(g_a)\subseteq\mu,
$$ 
so $|\rang(g_a)|$ is an aleph $\leq \mu$ (or
finite). Hence
$|\rang(g_a)|<\lambda$. By the choice of $D^*$ for some unique
$\gamma_a$ we have $B_a=:\{\alpha<\lambda:f_\alpha(a)=\gamma_a\}\in D^*$, as
$\langle f_\alpha:\alpha<\lambda\rangle$ exists also $\langle
\gamma_a:a \in A^*\rangle$ exists as well as $\langle B_a:a\in
A^*\rangle$. Clearly $|\{\gamma_a:a\in\A\}| < \theta(\A)$. As $D^*$
is $(\leq |A^*|)$-complete, $B^*=:\bigcap\limits_{a\in A^*} B_a\in D^*$, but
$[\alpha,\beta\in B^*\Rightarrow f_\alpha=f_\beta]$, a contradiction.

{\em 2.}\ \  Same proof. \hfill$\QED_{\ref{claim2:8}}$

\noindent{\sc Remark \ref{claim2:8}A:}\ \ You can also phrase the
theorem in terms of $\lambda$-complete filters on $\lambda$ which are weakly
$\kappa$-saturated (i.e. for every
$h:\lambda\longrightarrow\lambda'<\lambda$, for some
$C\subseteq\lambda'$ of cardinality $\leq\kappa$,
$\{\alpha<\lambda:h(\alpha)\in C\}\in D$). Here instead of 
$|A^*|<\lambda$ we need: $D$ is uniform and $\neg(\lambda\leq
|{}^{(A^*)}\kappa|)$.

\begin{definition}
\label{definition2:9}
Let $E$ be a set of filters on a set $\A$.
\begin{enumerate}
\item $E$ has the $I$-\lub\ property ($I$ a directed set) if:\\
\underline{for every} $D_0\in E$ and $F=\{f_t:t\in I\}\subseteq
{}^{(A^*)}\ord$ such that $\{f_t: t\in I\}$ is 
$<_{D_0}$-increasing,\\
\underline{there is} $D$, $D_0\subseteq D\in E$ such that there is a
$<_D$-\lub\ for $F$.
\item ``$E$ has the $I$-\eub\ property'' is defined similarly.
The $I\mbox{-}^*\eub$ is defined similarly for $<_D$-increasing
$\bar{F}=\langle F_t: t\in I\rangle$. Similarly $I\mbox{-}^*\lub$.
\item If $E=\{D\}$ we shall say $D$ has this property. 
\end{enumerate}
\end{definition}

\section{Nice families of filters}
\label{nice}
\noindent{\sc Hypothesis \ref{nice}.0}: $E$ is a family of
$\aleph_1$-complete filters on $A^*=\dom(E)$, such that (for
simplicity) $D\in E\ \&\ A\in D^+\Rightarrow D+A\in E$. 
Our main interest is in nice $E$ (defined below).

\begin{definition}
\label{hypothesis3:1}
We call $E$ nice if 
$$(\forall f:A^*\rightarrow\ord)(\forall D\in
E)(\rk^2_D(f)<\infty)\leqno\oplus_E$$  
(see below). Let $\oplus_{A^*}$
mean $\oplus_E$ for some $E$, with $A^*=\dom(E)$ and
$\oplus_\kappa[E]=\oplus_E$ with $\A=\kappa$. 

Let $E_{[D_0]}=\{D: D_0\subseteq D\in E\}$.
\end{definition}

\begin{definition}\ 
\label{definition3:2}
\begin{enumerate}
\item The truth value of $\rk^2_D(f,E)\leq \alpha$ (for $\alpha$ an ordinal
$f:A^*\rightarrow \ord$, $E$ usually omitted): 
\begin{description}
 \item[$\rk^2_D (f)\leq \alpha$] \underbar{if}
for every $A\in D^+$ and $f_1<_{D+A} f$ (and $f_1:A^*\rightarrow
\ord$) there is $D_1$, $D+A\subseteq D_1\in E$ and $\beta<\alpha$ such
that $\rk^2_{D_1}(f_1)\leq \beta$. 
\end{description}
(So $f =_D {\bf 0}_{\A}$ implies $\rk^2(f)\leq 0$)
\item $\rk^2_D(f)=\alpha$\ \ if\ \ $\rk^2_D(f)\leq \alpha$ and 
$\neg[\rk^2_d(f)\leq\beta]$ for $\beta <\alpha$\\
$\rk^2_D(f)=\infty$\ \ if\ \ $\rk^2_D(f)\leq \alpha$ for no $\alpha$\\
$\rk^3_D(f)=\min\{\rk^2_{D_1}(f):D\subseteq D_1 \in E\}$\\
(really we should have written $\rk^2_D(f,E)$ etc.).
\end{enumerate}
\end{definition}

\par \noindent
{\sc Remark:} Why start with $\rk^2$? To be consistent with [Sh-g]
Ch.V.

\begin{convention}
Let $f,g$ vary on ${}^{(A^*)}\ord$ and $D\in E$ and $A,B\subseteq A^*$.
\end{convention}

\begin{claim}
\label{claim3:3}
\begin{enumerate}
\item $\rk^2_D(f)\leq\beta$, $\beta\leq\alpha$ implies
$\rk^2_D(f)\leq\alpha$; so $\rk^2_D(f)$ is well defined as an ordinal or
$\infty$ (ZF is enough for the definition).
\item If $f\leq g$ or just $f\leq_D g$ and $l=2,3$ then
$\rk^l_D(f)\leq\rk^l_D(g)$ (so $f=_Dg$ implies $\rk^l_D(f)=\rk^l_D(g)$). 
\item In \ref{definition3:2}(1) we can demand in addition $f_1\leq f$.
\item $\rk^3_D(f)\leq\rk^2_D(f)$.
\item If $D_1\subseteq D_2$ then $\rk^3_{D_1}(f)\leq\rk^3_{D_2}(f)$.
\item For every $f, D$, for some $D_1$ such that $D\subseteq D_1\in E$ we
have 
\[\rk^3_D(f)=\rk^2_{D_1}(f)=\rk^3_{D_1}(f).\]
\item $\rk^2_D(f)=\sup\{\rk^3_{D+A}(g)+1: A\in D^+\mbox{ \ and\ }\
g<_{D+A}f\}$. 
\item If $A\in D^+$ then $\rk^3_D(f)\leq\rk^3_{D+A}(f)\leq\rk^2_{D+A}(f)\leq
\rk^2_D(f)$.\\ 
\relax[Why? By parts (5), (4), (7)+(4) respectively.] 
\item If $\rk^2_D(f)=\rk^3_D(f)$ \underline{then} for every $A\in D^+$
\[\rk^2_{D+A}(f)=\rk^3_{D+A}(f)=\rk^2_D(f)=\rk^3_{D}(f).\]
\item If $f<_D g$ \underline{then} $\rk^3_D(f)<\rk^3_D(g)$ or both are
$\infty$.\\
\relax [Why? Apply part (6) to $g$, $D$ and get $D_1$. Now
$$
\begin{array}{l}
\rk^3_D(f)\leq \rk^3_{D_1}(f)< \rk^3_D(f)+1\\
\qquad\leq \sup\{\rk^3_{D_1+A}(h)+1:
h<_{D_1+A} g\mbox{ and }A\in D^+_1\}\\
\qquad =\rk^2_{D_1}(g)= \rk^3_{D_1}(g).
\end{array}
$$
(Why the inequalities? by $D\subseteq D_1$, using part (5); trivially;
as $f<_{D_1} g$ hence $f<_D g$; by definition of $\rk^2$; and by the
choice of $D_1$ respectively.)]
\item $\|f\|_D\leq\rk^3_D(f)$ ($\|.\|_D$ -- the Galvin--Hajnal rank
which is defined by: $\|f\|_D=\sup\{\|g\|_D+1 : g <_D f\}$).\\ 
\relax [Why? part (10).]
\item If for $l=1,2$, $\rk^2_D(f_l)=\rk^3_D(f_l)=\alpha_l$ and
$\alpha_1<\alpha_2$ \underline{then} $f_1<_D f_2$. \\
\relax [Why? If not then for some $A\in D^+$ we have $f_2\leq_{D+A} f_1$, so
by part (2) we have $\rk^2_{D+A} (f_2)\leq \rk^2_{D+A}(f_1)$, but by part
(1) $\rk^2_{D+A}(f_2)=\alpha_2 > \alpha_1=\rk^2_{D+A}(f_1)$,
a contradiction.]
\item If for $l=1, 2$ $\rk^2_d(f_l)=\rk^3(f_l)=\alpha_l$ and
$\alpha_1=\alpha_2$ \underline{then} $f_1=_D f_2$.\\ 
\relax [Why? If not then by symmetry for some $A\in D^+$ we have $f_1<_{D+A}
f_2$, so by part (10) we have $\rk^3_{D+A} < \rk^3_{D+A}(f_2)$ and by part (9)
we get a contradiction.] 
\item If $\rk^2_{D_0}(f)<\infty$ and $\neg(g\geq_{D_0}f)$
\underline{then} for some $D_1\in E_{[D_0]}$,
$\rk^2_{D_1}(g)<\rk^2_{D_0} (f)<\infty$. 
\item If $\{A_t: t\in I\}\subseteq D^+$ and $(\forall D_1)(D\subseteq
D_1\in E \rightarrow (\exists t\in I)(A_t\in D^+_1))$
(e.g. $I$ finite, $\bigcup\limits_{t\in I} A_t\in D$) then 
$$
\rk^2_(f)=\sup \{\rk^2_{D+A_t}(f): t\in I\}\mbox{ and }
\rk^3_D(f)=\min\{\rk^3_{D+A_t}(f): t\in I\}
$$
[Why? The inequality $\geq$ by part (8), the other by the definition
(or part (7) and part (9)).]
\end{enumerate}
\end{claim}

\begin{observation}
\label{obs3:4A}
\begin{enumerate}
\item
\begin{description}
\item[(a)] If $f=g+1$, $D\in E$ then $\rk^2_D(f)\leq \rk^2_D(g)+1$; if
in addition $\rk^3_D(g)= \rk^2_D(g)$ 
\underline{then} $\rk^2_D (f)=\rk^3_D (g)+1$.
\item[(b)] If every $f(a)$ is a limit ordinal \underline{then}
$\rk^3_D (f)$ is a limit ordinal.
\item[(c)] $f=_D {\bf 0}_{\A}$ \underline{iff} $\rk^2_D (f)=0$.
\item[(d)] $\rk^3_D(f)=0$ iff $\neg (0_{\A}<_D f)$.
\end{description}
\item If $f=g+1$ \underline{then}:
\begin{description}
\item[(a)] $\rk^2_D (f)\in\{\rk^2_D (g),\ \rk^2_D (g)+1\}$,
\item[(b)] if $\rk^2_D (f)=\rk^2_D (g)<\infty$ \underline{then}
$\rk^2_D (g)$ is a limit ordinal of cofinality $<\theta({\cal P}(\A))$
and even $<\theta({\cal P}(\A)/D)$.
\end{description}
\item If $\rang(f)$ is a set of limit ordinals,
$F\subseteq\prod\limits_{a\in\A} f(a)$ is such that 
$$
[g<f\Rightarrow (\exists
h\in F)(g\leq h)]
$$ 
\underline{then} $\rk^2_D (f)=\sup\{\rk^3_{D+A}(h):
h\in F, A\in D^+\}$.
\item If $\delta=\rk^2_D (f)$, $\cf (\delta)\geq\theta (E)$
\underline{then}\footnote{ we can use this to prove a result parallel to
\protect\ref{claim2:3}\ (3)} for some $D_1$, $D\subseteq D_1 \in E$,
$\prod\limits_{a\in\A} f(a)/D_1$ has true cofinality $\delta$
(moreover by a smooth sequence $\langle F_\alpha:
\alpha<\cf(\delta)\rangle$) 
\item If $\rk^2_D(f)=\alpha<\infty$ \underline{then} for every
$\beta<\alpha$ there are $D$, $g$ such that $D\subseteq D_1\in E$,
$g<_{D_1} f$ and $\rk^2_{D_1}(g)=\rk^3_{D_1}(g)=\beta$.
\end{enumerate}
\end{observation}
\Proof 1)(a) By \ref{claim3:3}(7), \ref{claim3:3}(2),
\ref{claim3:3}(8) and common sense respectively 
$$
\begin{array}{ll}
\rk^2_D(f) & =\sup\{\rk^3_{D+A}(g')+1: A\in D^+\mbox{ and
}g'<_{D+A}f\}\\
\ & \leq \sup\{\rk^3_{D+A}(g)+1: A\in D^+\}\\
\ & \leq \sup\{\rk^2_D(g)+1\}= \rk^2_D(g)+1
\end{array}
$$
and by \ref{claim3:3}(4) and \ref{claim3:3}(10) respectively
$$
\rk^2_D(f)\geq \rk^3_D(f)\geq \rk^3(g)+1.
$$
Together we have finished.

\par \noindent
(b) By part (c) proved below, $\rk^3_D(f)>0$, so assume toward
contradiction that $\rk^3_D(f)=\alpha+1$, so (by \ref{claim3:3}(6))
for some $D_1$ we have $D\subseteq D_1$ ($\in E$) and
$\rk^3_D(f)=\rk^2_{D_1}(f)= \rk^3_{D_1}(f)$. By \ref{claim3:3}(7) for
some $A\in D^+$ and $g<_{D_1+A} f$ we have $\rk^2_{D_1+A}(g)=\alpha$.
By \ref{claim3:3}(10), \ref{claim3:3}(10), \ref{claim3:3}(9) and the 
choice of $\alpha$ respectively
$$
\rk^3_{D_1+A}(g)< \rk^3_{D_1+A}(g+1)<
\rk^3_{D_1+A}(f)=\rk^3_D(f)=\alpha+1.
$$
So $\rk^3_{D_1+A}(g)<\alpha$ contradicting the choice of $A$, and $g$.

\par \noindent
(c) If $f=_D 0_{\A}$ then the supremum in \ref{claim3:3}(7) (or the
definition) is taken on an empty set so $\rk^2_D(f)=0$. If
$\neg (f=_D 0_{\A})$ then $A=: \{a\in \A: f(a)>0\}$ and $g\in
{}^{(\A)}\ord$, $g(a)=0$ appears in the supremum in \ref{claim3:3}(7)
(or the definition) so $\rk^2_D(f)\geq 0+1=1>0$.

\par \noindent
(d) By (c) and the definition of $\rk^3$.
\medskip

\par \noindent
2) By \ref{claim3:3}(2), \ref{obs3:4A}(a)
respectively we have
$$
\rk^2_D(g)\leq \rk^2_D(f)\leq \rk^2_D(g)+1.
$$
Thus clause (a) holds. For clause (b) assume $\rk^2_D(f)=\rk^2_D(g)$,
and call this ordinal $\alpha$. If $\alpha=0$ we get contradiction by
\ref{obs3:4A}(1)(c) as $\rk^2_D(f)=\alpha$, $f=g+1$. If
$\alpha=\beta+1$ then for some $A\in D^+$ and $g_1<_{D+A} g$ we have
$\rk^3_{D+A}(g_1)=\beta$ so by \ref{claim3:3}(7), \ref{claim3:3}(10),
and the choice of $g_1$, $A$, and the choice of $\beta$ respectively:
$$
\rk^2_D(f)\geq \rk^3_{D+A} (g)+1 > \rk^3_{D+A}(g_1)+1 =\beta+1 = \alpha;
$$  
contradiction. So $\alpha$ is a limit ordinal and $\{\rk^3_{D+A}(g):
A\in D^+\}$ is an unbounded subset of $\alpha$ hence
$\cf(\alpha)<\theta({\cal P}(\A))$ (in fact $\cf(\alpha)< \theta({\cal
P}(\A)/D)$). 
\medskip

\par \noindent
3) By \ref{claim3:3}(7) and \ref{claim3:3}(2).
\medskip

\par \noindent
4) By \ref{claim3:3}(7), and $\rk^2_D(f)$ being a limit ordinal
$$
\begin{array}{ll}
\delta=\rk^2_D(f) & =\sup\{\rk^3_{D+A}(g)+1: A\in D^+\mbox{ and
}g<_{D+A} f \}\\
\ & =\sup \{\rk^3_{D+A}(g): A\in D^+\mbox{ and }g<_{D+A}f \mbox{ and
}g\leq f\}\\
\ & =\sup\{\rk^3_{D_1}(g): \mbox{ for some }A\in D^+, g<_{D+A}f
\mbox{ and }g\leq f\\
\ &\ \qquad\qquad \mbox{ and } D+A\subseteq D_1\in E\mbox{ and }\\
\ &\ \qquad\qquad\rk^3_{D_1}(g)= \rk^2_{D_1}(g)< \delta\}\\
\ & =\sup\{\rk^3_{D_1}(g): D\subseteq D_1\in E, g<_{D_1} f\mbox{ and
}g\leq f\\
\ &\ \qquad\qquad\mbox{ and }\rk^3_{D_1}(g)=\rk^2_{D_1}(g)<\delta\}\\
\ & =\sup\limits_{D_1\in E, D_1\supseteq D} \sup\{\rk^3_{D_1}(g): 
g<_{D_1} f\mbox{ and }g\leq f\\
\ &\ \qquad\qquad\mbox{ and }\rk^3_{D_1}(g)= \rk^2_{D_1}(g)<\delta\}.
\end{array}
$$
As $\cf(\delta)\geq \theta(E)$ necessarily for some $D_1$ we have
$D\subseteq D_1\in E$ and 
$$
\delta=\sup\{\rk^3_{D_1}(g): g<_{D_1}
f\mbox{ and }g\leq f \mbox{ and }\rk^3_{D_1}(g)= \rk^2_{D_1}(g)< \delta\}.
$$
For $\alpha< \delta$ let 
$$
F_\alpha=\{g: g\leq f, g<_{D_1} f, \mbox{ and
}\alpha=\rk^3_{D_1}(g)=\rk^2_{D_1}(g)\}
$$
and $S=\{\alpha< \delta: F_\alpha\neq \emptyset\}$, it is necessarily
unbounded in $\delta$. Now by \ref{claim3:3}(12), (13), $\bar
F=\langle F_\alpha: \alpha\in S\rangle$ is $<_{D_1}$-increasing and smooth.

\par \noindent
5) Suppose $\beta<\alpha$, $f$, $D$ form a counterexample. Then we
prove by induction on $\gamma\geq \beta$ that 
\begin{description}
\item[$(*)_\gamma$] if $g\in {}^{(\A)}\ord$, and $D\subseteq D_1\in E$
and $g<_{D_1} f$ and $\rk^3_{D_1}(g)\geq \beta$ \underline{then}
$\rk^3_D(g)\geq \gamma$.
\end{description}

For $\gamma=\beta$ clearly $(*)_\gamma$ holds trivially.

For $\gamma>\beta$, we can find $D_2$ such that $D_1\subseteq D_2\in
E$ and $\rk^3_{D_1}(g)=\rk^2_{D_2}(g)= \rk^3_{D_2}(g)$, but $\beta\leq
\rk^3_{D_1}(g)$ by assumption, so as $\beta$, $\alpha$, $f$, $D$ form
a counterexample $\beta<\rk^3_{D_1}(g)$ so by \ref{claim3:3}(7)
applied to $\rk^2_{D_2}(g)$ we can find $A\in D^+_2$ and $g_1<_{D_2+A}
g$ such that $\rk^3_{D_2+A}(g_1)\geq \beta$, so by the induction
hypothesis $\rk^3_{D_1+A}(g_1)\geq \bigcup\{\gamma_1: \beta\leq
\gamma_1<\gamma\}$. Now by \ref{claim3:3}(9) we know
$\rk^3_{D_2+A}(g)= \rk^3_{D_2}(g)$. Together we get the required
conclusion in $(*)$.
So $\rk^2_{D}(f)= \sup \{\rk^3_{D+A}(g)+1: A\in D^+\mbox{ and
}g<_{D+A}f\}$ is $\alpha$ hence $>\beta$ so for some $A\in D^+$ and
$g<_{D+A} f$ we have $\rk^3_{D+A}(g)\geq \beta$ hence by $(*)$ it is
$>\alpha$, contradiction. \hfill\qed$_{\ref{obs3:4A}}$

\begin{definition}
\label{definition3:4} 
For $D$ a filter on $A^*$ and $f:A^*\rightarrow \ord$ let
\[\ts_D(f)=\{F: F\subseteq
\prod\limits_{i\in A^*} f(i)\mbox{ satisfies } f_1\in F\ \&\ f_2 \in F\ \&\
f_1\neq f_2\Rightarrow f_1\neq_D f_2\}\]
$$
\begin{array}{ll}
\tw_D(f)=\big\{(F,{\bf e}):& F\subseteq
\prod\limits_{i\in A^*} f(i),\ {\bf e} \mbox{ is an equivalence relation on $F$},\\
\ & \neg f_1\ {\bf e}\ f_2 \Rightarrow f_1 \neq_D f_2\big\},
\end{array}
$$
$$
\bT_D(f)=\{(F, {\bf e})\in {\it Tw}_D(f): {\bf e}\mbox{ is
}=_D\restriction F\},
$$
$${\it Ts}_D(f)=\sup\{|F|: F\in\ts_D(f)\},$$
$${\it Tw}_D(f)=\sup\{|F/{\bf e}|: F/{\bf e}\in\tw_D(f)\},$$ 
$$
T_D(f)=\sup\{|F/{\bf e}|: (F, {\bf e})\in \bT_D(f)\}
$$
(it is a kind of cardinality; of course the $\sup$ do not necessarily
exist). We may write $F/{\bf e}\in \tw_D(f)$ 
instead of $(F, {\bf e})\in \tw_D(f)$ and also may write $\bar F=\langle
F_t: t\in I\rangle$ for $(F, {\bf e})$ if $F=\bigcup\limits_t F_t$,
and $f {\bf e} g \leftrightarrow (\exists t)[\{f, g\}\subseteq F_t]$
and so the
$F_t$'s are pairwise disjoint.
\end{definition}

\begin{observation}
\label{obs3:7} 
\begin{enumerate}
\item $F\in\ts_D(f)\ \Rightarrow\ (F,=)\in\tw_D(f)$.
\item ${\it Ts}_D(f)\leq{\it Tw}_D(f)$.
\item If $\bar{F}=\langle F_\alpha:\alpha<\alpha^*\rangle$ is
$<_D$-increasing, $\bigcup F_\alpha\subseteq\prod\limits_{a\in\A}f^*(a)$
\underline{then} 
$\rk^3_D(f^*)\geq\alpha^\ast$; so $\hpp_D(f^*)\leq\rk^3_D(f^*)$; also
$\langle F_\alpha: \alpha<\alpha^*\rangle\in \tw_D(f)$ hence
$|\alpha^*|\in {\it Tw}_D(f)$.
\item $[AC_{\A}]$ If $f, g\in {}^{\A}\ord$ and $D$ is a filter on $\A$
and $\{a\in \A: f(a)=g(a)\}\in D$ and $X\in \{T_s, Tw, T\}$
\underline{then} $X_D(f)= X_D(g)$.
\end{enumerate}
\end{observation}

\begin{claim}
\label{claim3:6}
\begin{enumerate}
\item \mbox{[$E$ is nice]}\ \ \ \ Assume $F\subseteq {}^{(A^*)}\ord$.  We
can find $\langle(F_D,h_D): D\in E\rangle$ such that:
\begin{description}
\item[(a)] $F=\bigcup\limits_{D\in E} F_D$
\item[(b)] $h_D:F_D\rightarrow \ord$
\item[(c)] if $f_1,f_2 \in F_D$ \underline{then} $h_D(f_1)\leq
h_D(f_2)\iff f_1\leq_D f_2$ and hence 
\[h_D(f_1)=h_D(f_2)\iff f_1=_D f_2,\mbox{ so }\]
\[h_D(f_1)<h_D(f_2)\iff f_1<_D f_2\]
\item[(d)] if $F\subseteq \prod\limits_{i\in A^*}f(i)$ \underline{then} $(F_D,=_D)\in\tw(f)$.
\end{description}
\item Instead of ``$E$ is nice'' it suffices that
$F\subseteq\prod\limits_{a\in\A}f^*(a)$, $\rk^2_D(f^*)<\infty$.
\end{enumerate}
\end{claim}

\Proof
{\em 1)} Let $F_D=:\{f\in F: \rk^2_D(f)=\rk^3_D(f) (<\infty)\}$,
$h_D(f)=\rk^2_D(f)$. Now clause (a) holds by \ref{claim3:3} $(6)$,
clause (b) holds as $E$ is nice (see \ref{hypothesis3:1}).
For (c): 
it holds by \ref{claim3:3}(12), \ref{claim3:3}(13).
Lastly clause (d) follows from (c). 

{\em 2)} Similar (see \ref{claim3:3} $(14)$).\hfill$\QED_{\ref{claim3:6}}$

\begin{conclusion}
\label{conclusion3:7}
Let $D_0\in E$.
\begin{enumerate}
\item $\mbox{[ } E\mbox{ is nice ]}$  Assume $F\subseteq
\prod\limits_{a\in\A} f^*(a)$ is in $\ts_{D_0} (f^\ast)$.  
\begin{description}
\item[(a)\ ]\ If $|E|$ is an aleph \underline{then} so is $|F|$.
\item[(b)\ ]\ $F$ can be represented as $\bigcup \{F_{D}: D_0\subseteq D\in
E\}$ such that $F_D$'s cardinality is  an aleph $<\ehpp_E(f^\ast)$.
\item[(c)\ ]\ $|F|\leq^* |E|\times \ehpp_E(f^*)$ 
\end{description}
\item $\mbox{[ } E\mbox{ is nice ]}$  Assume $(F,{\bf e})$ is in
$\tw_{D_0}(f^*)$ \\
(any $(F,=_D)$ will do for  $F\subseteq {}^{(\A)}\ord$). 
\begin{description}
\item[(a)\ ]\  If $|E|$ is an aleph \underline{then} so is $|F/{\bf e}|$.
\item[(b)\ ]\  $F/{\bf e}$ can be represented as $\bigcup\{F_D/{\bf e}: D_0 
\subseteq D\in E\}$ such that $|F/{\bf e}|$ is an aleph $\leq \hpp(f^*)$.
\item[(c)\ ]\  $|F/{\bf e}|\leq^* |E|\times\hpp_E (f^*)$ 
\end{description}
\item Instead of ``$E$ is nice'', ``$\rk^2_{D_0} (f^*)< \infty$'' suffices.
\end{enumerate}
\end{conclusion}

\Proof {\em 1)}\ \ It follows from {\em 2)} as in this case ${\bf e}$ is the
equality. 

\noindent{\em 2)}\ \ For $D\in E_{[D_0]}$, let  
$$
\begin{array}{ll}
F_D=\big\{f\in F:\rk^2_D(f)=\rk^3_D(f)\mbox{ and:}& \mbox{ if } f_1\in f/{\bf
e}\mbox{ and }
\rk^2_D(f_1)=\rk^3_D(f_1)\\
\ & \mbox{ then }\rk^2_D(f_1)\geq\rk^2_D(f)\}.
\end{array}
$$
Let $h_D:F_D\rightarrow \ord$ be $h_D(f)=\rk^2_D(f)$. So (for $f_1,f_2\in F_D$): 
\[h_D(f_1)=h_D(f_2)\ \mbox{ iff }\ f_1/{\bf e}=f_2/{\bf
e}.\]
Also $F/{\bf e}=\bigcup \{F_D/{\bf e} : D\in E\}$. So clause (b)
holds, for clause (c) let $G(D, \alpha)$ be $y$ iff $y$ has the form
$f/{\bf e}$ where $f\in F_D$, $h_D(f)=\alpha$; so $G$ is a partial
function from $|E|\times \hpp_E(f^*)$ onto $F/{\bf e}$. Note that
clause (a) follows as $|X|\leq |Y|$, $|Y|$ is an aleph implies $|X|$ is an
aleph.  
(We can choose a well ordering $<$ of $E$, we can let $F'_D=\{f: f\in
F_D$ and letting $f\in F_t$, for no $f'\in f/{\bf e}$ and $D'<D$ do we
have $f'\in F_{D'}\}$. Let $\bigcup \{ h_D\restriction F'_D: D\in E\}$
so letting $h(f/{\bf e})=h_D(f)$ for $f\in F'_D$, 
clearly $h$ is one-to-one function from $F$ onto a set of ordinals.)

\noindent{\em 3)}\ \ Similar proof (remember \ref{claim3:3} $(13)$).
$\hfill\QED_{\ref{conclusion3:7}}$

\begin{claim}
\label{claim3:8}
\begin{enumerate}
\item If $F$ is in $\ts_D(f)$ and $\rk^3_D(f)=\alpha <\infty$
\underline{then} $|F|\leq^* |\alpha|\times |E|$.
\item Assume $(F,{\bf e})\in\tw_D(f)$ and $\rk^3_D(f)=\alpha <\infty$.
\underline{Then} $|F/{\bf e}|\leq^* |\alpha|\times |E|$.
\end{enumerate}
\end{claim}

\Proof Included in the proof of \ref{conclusion3:7}.

\begin{claim}
\label{claim3:9}
If $\alpha=\rk^2_{D_0}(f^*)$, $D_0 \in E$ \underline{then} we can find
$\langle 
A_D: D\in E_{[D_0]}\rangle$, $\alpha=\bigcup\{A_D:D\in E_{[D_0]}\}$, and
$\langle(F_D,h_D): D\in E_{[D_0]}\rangle$ with $(F_D,h_D)$ as in
\ref{claim3:6}, $F_D\subseteq \prod\limits_{a\in\A} f^*(a)$
 with $\rang(h_D)=A_D$ and 
$$
h_D(f)=\alpha\ \Rightarrow\
\rk^2_D(f)=\rk^3_D (f)=\alpha.
$$
\end{claim}

\Proof Because there is ``no hole in the possible
ranks''. I.e. we apply \ref{claim3:6} to $F=\prod\limits_{a\in\A}(f^*(a)+1)$,
and get $\langle (F_D ,h_D): D\in E_{[D_0]}\rangle$. Our main problem is
that some $\beta <\alpha$ is not in $\bigcup\limits_{D\in
E{[D_0]}}\rang(h_D)$. Then use \ref{obs3:4A}(5).
\hfill$\QED_{\ref{claim3:9}}$ 

\begin{conclusion}
\label{conclusion3:??}

{\em Assume $E$ is nice and for simplicity $|E|$ is an aleph, \\
$f\in
{}^{(\A)}\ord$  and $(\forall a\in\dom(E) )$\ 
$(f(a)\geq |E|)$.
\begin{enumerate} 
\item Then the cardinals 
$$\sup_{D\in E_{[D_0]}}\!|\rk^2_D(f)|,\ 
\sup_{D\in E_{[D_0]}}\!({\bf T}_D(f)),\ 
\sup_{D\in E_{[D_0]}}\!(Tw_D(f)) \mbox{ and }
\sup_{D\in E_{[D_0]}}\!\{\hpp_D(f)\}
$$ 
are equal. 
\item Assume $AC_{\A}$. If $f_1$, $f_2$ are as in the assumption and
$\{a\in \A: |f_1(a)|=|f_2(a)|\}\in D$ then
$$
\sup\limits_{D\in E_{[D_0]}} |\rk^2_D(f_1)|=\sup\limits_{D\in
E_{D_0}}|\rk^2_D(f_2)| 
$$
have the same cardinality.
\item $\rk^2_{D_0}(f^+)\leq \sup\limits_{D\in E_{[D_0]}}
|\rk^2_D(f)|^+$ where $f^+\in {}^{(\A)}\ord$ is $f^+(a)=: |f(a)|^+$
\end{enumerate}
} 
\end{conclusion}
\medskip

\Proof \ 1) 
\par \noindent
\underline{Step A}
\begin{description}
\item[$(*)$] $\sup\limits_{D\in E_{[D_0]}}(\rk^2_D(f)+1)\geq
\hpp_{E_{[D_0]}}(f)$ 
\end{description}

If $\beta< \hpp_{E_{[D_0]}}(f)$ then we can find $D\in E_{[D_0]}$ and
$\bar F=\langle F_\alpha: \alpha<\beta\rangle$ such that
$F_\alpha\subseteq \prod\limits_{a\in \A} f(a)$ non empty and $\bar F$
is $<_D$-increasing. We can prove using \ref{claim3:3}(9) by induction
on $\alpha< \beta$ that 
\begin{description}
\item[$(**)$] $g\in F_\alpha \Rightarrow \rk^3_D(f)\leq \alpha$.
\end{description}

Now by \ref{claim3:3}(7) we have $\rk^2_D(f)>\alpha$ and
$\alpha<\beta$ so $\rk^2_D(f)\geq \beta$.

\par \noindent
\underline{Step B}

$\rk^2_D(f)$ can be represented as the union of $E$ sets each of order
type $<\hpp_{E_{[D_0]}}(f)$.

By \ref{claim3:8}.

\par \noindent
\underline{Step C}

$\sup\limits_{D\in E_{[D_0]}} |\rk^2_D(f)|= \hpp_{E_{[D_0]}} (f)$.

Why? By steps A, B as $|E|$ is an aleph and $f(a)\geq |E|$ for every
$a$.

\par \noindent
\underline{Step D}

$\hpp_{E_{[D_0]}} (f) \leq {\it T}_D(f)\leq {\it Tw}_D(f)$.

By definition (true for each $D$ separately).

\par \noindent
\underline{Step E}

$\sup\limits_{D\in E_{[D_0]}} {\it T}_D(f)\leq \sup\limits_{D\in E_{[D_0]}}
\rk^2_D(f)$.

Why? By \ref{claim3:6}

\par \noindent
2) By part (1) and \ref{obs3:7}(4). 

\par \noindent
3) By part (2) and the definition of $\rk^2_{D_0}$.
\hfill\qed$_{\ref{conclusion3:??}}$
\medskip

\par \noindent
{\sc Remark \ref{conclusion3:??}A:}
If we waive ``$|E|$ is an aleph'' but still assume $f(a)\geq
\theta(E)$ we get that those three cardinals are not too far.
\medskip

\begin{claim}
\label{claim3:11}  Assume
\begin{description}
\item[(A)\ \ \ ] $\lambda$ is an aleph,
$\langle\lambda_i:i<\delta\rangle$ is a strictly increasing continuous
sequence of alephs with limit $\lambda$,
\item[(B)\ \ \ ] there is $G^*=\bigcup\limits_{i<\delta} G_i$, $G_i$ a one-to-one
function from $\prod\limits_{j<i}\lambda_j$ into some 
$\lambda^*_i<\lambda_{i+1}$ ( so if $AC_\delta$, we need just
that each $G_i$ exists, i.e. $|\prod\limits_{j<i}\lambda_j|<\lambda_{i+1}$). 
\end{description}
\underline{Then} there is $F\in\ts_{D_\delta^{bd}}
(\langle\lambda^*_i:i<\delta\rangle)$, and $|F|=\prod\limits_{i<\delta}
\lambda_i$.\\
($D^\bd_\delta$ is the filter generated by the cobounded sets of $\delta$).
\end{claim}

\noindent{\sc Remark:}\ \ \ This goes back to Galvin, Hajnal [GH], see [Sh
386, 5.2A(1)]. 

\Proof Define a function $G:\prod\limits_{i<\delta}\lambda_i\longrightarrow
\prod\limits_{i<\delta}\lambda_i^*$ by\\
$G(f)(i)=G_i(f\restriction i)$ and let
$F=\rang(G)$\hfill$\QED_{\ref{claim3:11}}$ 

\begin{claim}
\label{claim3:12}
Assume that $D_0\in E$ and $f_\alpha\in {}^{(A^*)}\ord$ for $\alpha<\alpha^*$.
\begin{enumerate}
\item $\mbox{[ } E\mbox{ is nice ]}$ There are $g\in {}^{(A^*)}\ord$ and
$D$ such that: 
\begin{description}
\item[(a)\ ] $D_0\subseteq D\in E$,
\item[(b)\ ] $\rk^2_D(g)=\rk^3_D(g)$,
\item[(c)\ ] $g$ is a $<_D$-\lub\ of $\{f_\alpha:\alpha<\alpha^*\}$,
i.e. :
\begin{description}
\item[(i)\ ] $\alpha<\alpha^*\ \Rightarrow\ f_\alpha\leq_D g$
\item[(ii)] if $g'\in{}^{(A^*)}\ord$ and $f_\alpha\leq_D g'$ for all
$\alpha<\alpha^*$ then $g\leq_D g'$.
\end{description}
\item[(d)\ ] Moreover for every $A\in D^+$, $g$ is also a $<_{D+A}$-\lub\ of
$\{f_\alpha:\alpha<\alpha^*\}$.
\end{description}
\item If $f_\alpha\leq_{D_0} f^*$ for $\alpha<\alpha^*$, $\rk^2_{D_0}
(f_{\alpha}) <\infty$ then there are $g\in {}^{(\A)}\ord$, $D$
satisfying (a)-(d) and we can add 
\begin{description}
\item[(e)\ ] $g\leq f^*$.
\end{description}
\end{enumerate}
\end{claim}

\Proof
{\em 1)}  Follows by (2), just let $f^* (a)=\sup_{\alpha <\alpha^*}
f_{\alpha}(a)$.

{\em 2)} Clearly the set 
\[ K=\{(D,g): D_0\subseteq D\in E\ \&\
(\forall\alpha<\alpha^*)(f_\alpha \leq_D g),\ \rk^2_D (g)<\infty\
\mbox{ and }\ g\leq f^* \}\] 
is not empty (because the pair $(D_0, f^*)$ is in it ). Choose
among those pairs one $(D_1,g)$ with $\rk^3_D(g)$ minimal. So for some $D$,
$D_1\subseteq D\in E$ we have $\rk^3_{D_1}(g)=\rk^2_D(g)=\rk^3_D(g)$ (see
\ref{claim3:3}(6)), so also $(D,g)\in K$. Clearly $(D,g)$
satisfies (a), (b), (e) and 
(c)(i). If (c)(ii) fails then let $g'$ exemplify it and so
\[A=\{a\in\kappa: g'(a)<g(a)\}\neq\emptyset\mod D.\]
But clearly also $(D+A,g')$ is in the family and (see
\ref{claim3:3}(10) and \ref{claim3:3}(9) respectively): 
\[\rk^3_{D+A}(g')<\rk^3_{D+A}(g)=\rk^3_D(g)\]
- a contradiction.

Clause (d) follows from \ref{claim3:3}(9) replacing $D$ by $D+A$.
\hfill$\QED_{\ref{claim3:12}}$

\begin{definition}
\begin{enumerate}
\item Let $E$ be as in $3.0$, $\kappa$ an aleph and
$\kappa=\cf(\kappa)>\aleph_0$, and $D^\oplus$ a filter on $\kappa$. We
say that $E$ is [weakly] $D^\oplus$-normal if there is a function $\i$
witnessing it which means:

$\i$ is a function from $A^*=\dom (E)$ to $\kappa$ such that 
\begin{description}
\item[(a)\ ] for every $D\in E$
\[i(D)=:\{A\subseteq\kappa: \i^{-1}(A)\in D\}\]
is a filter on $\kappa$ extending $D^\oplus$, 
\item[(b)\ ] every $D\in E$ is [weakly] $\i$-normal, which means: if 
$f:A^*\longrightarrow\kappa$ is $\i$-pressing down (i.e.
$f(a)<1+\i(a)$ for $a\in \A$) then
on some $A\in D^+$, the function $f$ is constant [bounded].
\end{description}
\item We say that $E$ is [weakly] $\kappa$-normal if this holds for
some $\kappa$ as above, 
filter $D^\oplus$ over $\kappa$ and function $\i$.
\item We say that $E$ is [weakly] $(\kappa,S)$-normal, where
$S\subseteq\kappa$, if this holds for some filter $D^\oplus$
over $\kappa$ such that $S\in D^{\oplus}$.
\item We say that $D$ is $\kappa$-normal (with $\i$ witnessing it) if
$E:=\{D+A: A\in D^+\}$ is, with $\i$ witnessing it.\\
If $\dom (E)=\kappa$, $\i$
the identity we omit it.\\ 
Omitting $D$ from ``$D$-normal'' means omitting clause (a).
\end{enumerate}
\end{definition}

\begin{remark}
{\em
\begin{enumerate}
\item Note: for a filter $D$ on $\kappa$ the normality is also defined as
the closure under diagonal intersection; this is equivalent. But it is
not enough that the diagonal intersection of clubs is a club, we need
that the diagonal intersection of sets including clubs includes a club. 
\item The club filter $D^{\rm cl}_\kappa$ on a regular
$\kappa>\aleph_0$ is not 
necessarily normal, but there is a minimal normal filter on it,
$D^\otimes_\kappa$, but possibly $\emptyset\in D^\otimes_\kappa$, see below. 
\end{enumerate}}
\end{remark}

\begin{definition}
\label{definition3:?}
\begin{enumerate}
\item For an ordinal $\delta$
\[D^{\rm cl}_\delta=\{A\subseteq\delta: A \mbox{ contains a club of }
\delta\}.\] 
\item For ${\cal P}\subseteq {\cal P}(\delta)$ we define  $D^{\rm
nor}_{\delta,\zeta}[{\cal P}]$ by induction on $\zeta$ as follows:
\begin{description}
\item[$\zeta=0$]\ \ \ \ \ \ \ \ \ $D^{\rm nor}_{\delta,\zeta}[{\cal P}]$ is
the filter of subsets of $\delta$ generated by ${\cal P}\cup D^{\rm
cl}_\delta$. 
\item[$\zeta>0$ limit]\ \  $D^{\rm nor}_{\delta,\zeta}[{\cal
P}]=\bigcup\limits_{\xi<\zeta} D^{\rm nor}_{\delta,\xi}[{\cal P}]$.
\item[$\zeta=\xi+1$]\ \ \ \  $D^{\rm nor}_{\delta,\zeta}[{\cal
P}]=\{\delta\setminus B:$ there is a function with domain $B$, regressive
(i.e. $f(\alpha)<1+\alpha$), such that for every $\beta<\delta$ the set
$\{\alpha\in \delta\setminus B:
f(\alpha)=\beta\}=\emptyset\mod D^{\rm nor}_{\delta,\xi}[{\cal P}]\}$ 
\item $D^{\rm nor}_{\delta}[{\cal P}]=\bigcup\limits_{\zeta} D^{\rm
nor}_{\delta,\zeta}[{\cal P}]$
\item Above we replace {\rm nor} by {\rm wnr} if for $\zeta=\xi+1$ we
replace $f(\alpha)=\beta$ by $f(\alpha)\leq\beta$
\item If ${\cal P}=\emptyset$ we omit it.
\end{description}
\item We call $A\subseteq\kappa$ stationary if $A\neq\emptyset\mod D^{\rm
nor}_{\kappa}$.
\end{enumerate}
\end{definition}

\begin{claim}
\label{claim3:15}
\begin{enumerate}
\item $D^{\rm nor}_{\delta,\zeta}[{\cal P}]$ increases with $\zeta$ and is
constant for $\zeta\geq\theta({\cal P}(\delta))$ (so $D^{\rm
nor}_{\delta}[{\cal P}]=D^{\rm nor}_{\delta,\theta({\cal P}(\delta))}[{\cal
P}]$), moreover it is constant for $\zeta\geq \zeta^{\rm nor}_\delta$ for
some $\zeta^{\rm nor}_\delta< \theta({\cal P}(\delta))$. Note that 
we get $D^{\rm cl}_\delta\subseteq D^{\rm nor}_{\delta,1}[{\cal
P}]$ even if we redefine $D^{\rm nor}_{\delta,0}[{\cal P}]$ as the filter
generated by ${\cal P}\cup D^{\rm bd}_{\delta}$. Clearly
$D^{\rm nor}_{\delta} [{\cal P}]$ is the minimal normal filter on
$\delta$ which includes ${\cal P}$.
\item $D^{\rm wnr}_{\delta,\zeta}[{\cal P}]$ increases with $\zeta$ and it
is constant for $\zeta\geq\theta ({\cal P}(\delta))$ (consequently $D^{\rm wnr}_\delta[{\cal
P}]=D^{\rm wnr}_{\delta,\theta({\cal P}(\delta))}[{\cal 
P}]$ and $D^{\rm wnr}_{\delta,\zeta}[{\cal P}]\subseteq D^{\rm
nor}_{\delta,\zeta}[{\cal P}]$), moreover for $\zeta \geq \zeta^{\rm
wnr}_\delta$ for some $\zeta^{\rm wnr}_\delta< \theta({\cal
P}(\delta))$. Also we get $D^{\rm cl}_\delta\subseteq 
D^{\rm wnr}_{\delta,1}$ even if we redefine $D^{\rm
wnr}_{\delta, 0}[{\cal P}]$ as the filter generated by ${\cal P}\cup
D^{\rm bd}_{\delta}$. Clearly 
$D^{\rm wnr}_{\delta}[{\cal P}]$ is the minimal weakly-normal filter on
$\delta$ which includes ${\cal P}$ 
\item In 1), 2) if $\cf(\delta)=\delta_1$, the filters $D^{\rm
wnr}_\delta$, $D^{\rm wnr}_{\delta_1}$ (and also $D^{\rm wnr}_{\delta,
\zeta}$, $D^{\rm wnr}_{\delta_1, \zeta}$)
 are essentially the same. I.e. let
$h\colon\delta_1\longrightarrow\delta$ (strictly) increases
continuously with unbounded range, \underline{then}
\begin{quotation}
if $A\subseteq\delta$, $A_1\subseteq\delta_1$, $A\cap\rang(h)=h''(A_1)$

then $A\in D^{\rm wnr}_{\delta,\zeta}[{\cal P}] \iff A_1\in D^{\rm
wnr}_{\delta_1,\zeta}[{\cal P}]$. 
\end{quotation}
\item If $\delta$ is not a regular uncountable cardinal then $D^{\rm
nor}_{\delta}={\cal P}(\delta)$. 
\item  $[AC_{\delta, {\cal P}(\delta)} + DC]$ If $\delta$ regular
uncountable, then in 1), 2) we can replace 
$\theta({\cal P}(\delta))$ by $(2^{\cf(\delta)})^+$ and we have
$\emptyset\notin D^{\rm nor}_\delta$. Moreover for every regular
$\sigma<\delta$, 
\[\{\alpha<\delta\colon \
\cf(\alpha)=\sigma\}\neq\emptyset\ \mod D^{\rm 
nor}_\delta.\]
\end{enumerate}
\end{claim} 

\Proof 1) - 4) Check.

\par \noindent
5) We use the variant of the definition starting with $D^{\bd}_\delta$.
We choose by induction on $n<\omega$, an equivalence relation $E_n$
on $\delta$ such that:
\begin{description}
\item[(i)] $E_0$ is the equality on $\delta$,
\item[(ii)] $E_{n+1}$ refine $E_n$,
\item[(iii)] the function $f_n$ is regressive where 
$$
f_n(\alpha)=\otp\{\beta: \beta E_n \alpha\mbox{ but }\beta<\alpha \
\&\ \neg \beta
E_{n+1} \alpha\mbox{ and } \beta=\min(\beta/E_{n+1})\}
$$
(so it is definable from $E_n$, $E_{n+1}$),
\item[(iv)] for each $\alpha<\delta$ and $n<\omega$
$$
\zeta^{n+1}_\alpha=: \min\{\zeta: \alpha/E_{n+1}\in D^{\rm
nor}_{\delta, \zeta}\}< \zeta^{n}_\alpha =: \min\{\zeta: \alpha/
E_n\in D^{\rm nor}_{\delta, \zeta}\}
$$
\underline{or} both are zero.
\end{description}

We can carry the induction by $DC$. 
For $n=0$ use clause (i) to define $E_0$.
For the choice of $E_{n+1}$, for
each $\alpha<\delta$, $f_n\restriction (\alpha/ E_n)$ as required exists by the
inductive definition of $D^{\rm nor}_{\delta, \zeta}$ (does not matter
if we let $D^{\rm nor}_{\delta, 0}$ be $D^{\cl}_\delta$ or
$D^{\bd}_\delta$), but we have to choose 
$$
\langle f_n\restriction
(\alpha/ E_n): \alpha<\delta, \alpha=\min(\alpha/E_n)\rangle;
$$ 
so we have to make $\leq |\delta|$ choices, each among the family of
regressive function on 
$\delta$. But as we have a pairing function on $\delta$, this is
equivalent to a choice of a subset of $\delta$, so $AC_{\delta, {\cal
P}(\delta)}$ which we assume is enough.

Now by $DC$ we can choose by induction on $n<\omega$, an ordinal
$\alpha_n<\delta$ such that 
\begin{description}
\item[$(*)$] if $\beta\leq \alpha_n$, $m<\omega$ and $\zeta^\beta_m=0$
then $\sup(\beta/E_m)<\alpha_{n+1}$ (note $\beta/E_m$ is bounded in $\delta$).
\end{description}
(Not hard to show that $\alpha_{n+1}$ exists.) Now letting
$\alpha(*)=\bigcup\limits_{n<\omega} \alpha_n< \delta$, we get easy
contradiction as $\langle \zeta^{\alpha(*)}_n : n<\omega\rangle$ is
eventually zero.
\hfill\qed$_{\ref{claim3:15}}$

\section{Investigating strong limit singular $\mu$}
\label{par4}

\begin{definition}\
\label{definition4:1}
\begin{enumerate}
\item $\otimes_{\alpha,R}$\ \ \  means: there is a function $e$
exemplifying it which means:
\begin{description}
\item[$\otimes_{\alpha,R}\mbox{[$e$]}$]\ \ \ $e$ is a function,
$\dom(e)=\{\delta:\delta<\alpha$ a limit ordinal $\}$ and for every limit
$\delta<\alpha$, $e(\delta)$ is an unbounded subset of $\delta$ such that 
it is of order type from $R$. (It follows that $R\supseteq \alpha\cap \reg$).
\end{description}
\item If $R\cap\alpha$ is the set of infinite regular cardinals $\leq
\alpha$ we omit it (then $\otp(e(\delta))=\cf(\delta)$). \\
If $R$ is the set of regular cardinals $<\alpha$ union with $\sigma$
(not $\{\sigma\}$!) we write $\sigma$ instead of $R$.
\item Let $\otimes^*_{\alpha, R}$ means $\otimes_{\alpha,
R\cup(\alpha\cap\reg)}$.  
\end{enumerate}
\end{definition}

\begin{observation}\
\label{observation4:2}
\begin{enumerate}
\item If $\otimes_{\alpha,R}$, $\sigma\in R$ is not regular and
$R'=R\setminus\{\sigma\}$ \underline{then} $\otimes_{\alpha,R'}$.
\item If there is $e$ satisfying $\otimes'_{\alpha,R}[e]$ below then
$\otimes_{\alpha,R}$ holds (for another $e$)
\begin{description}
\item[$\otimes'_{\alpha,R}\mbox{$[e]$}$]\ \ \ $e$ is a function,
$\dom(e)=\{\delta:\delta<\alpha$ a limit ordinal, $\delta\notin R\}$ and for
every limit $\delta \in \dom(e)$, $e(\delta)$ is an unbounded subset of
$\delta$, of order type $<\delta$.
\end{description}
\item Also the converse of (2) is true.
\item  If $\otimes_{{\alpha_2},{\alpha_1}\cup{R_1}}$ and
$\otimes_{{\alpha_1},{\alpha_0}\cup{R_0}}$ then
$\otimes_{{\alpha_2},{\alpha_0}\cup {R_0}\cup{R_1}}$.
\item  For any ordinal $\alpha$ letting $j_{|\alpha|}$ be 1 if
$|\alpha|$ is regular and zero otherwise we have
$\otimes_{\alpha,|\alpha|+j}$, hence if 
$\otimes_{\alpha,\beta}$ then $\otimes_{\alpha,|\beta|+j_\beta}$.
\item  If $\otimes^*_{\alpha,\zeta}$ then we can define $\langle
f_\beta:\beta\in [\zeta,\alpha]\rangle$, $f_\beta$ a one-to-one function
from $\beta$ onto $|\beta|$.
\item If $\otimes_{\alpha, R}$ then we can define $\bar{f}=\langle
f_\beta:\beta<\alpha\rangle$, $f_\beta$ a one-to-one function from $\beta$
onto $\sup(\beta\cap R)$.
\item  If $\otimes_{\alpha,\sigma}$ and  $\sigma<\lambda^+\leq\alpha$ (so
$\lambda^+$ is a successor) \underline{then} $\lambda^+$ is regular.
\item If we have: $R\subseteq\alpha$ closed, $\bar{f}=\langle
f_\beta:\beta<\alpha\rangle$, $f_\beta$ a function from $\sup(R\cap\beta)$
onto $\beta$ then $\otimes_{\alpha,R}$. 
\item If $\otimes_{\alpha , R}$ and $\beta <\gamma <\alpha$ and $[\beta,
\gamma]\cap R=\emptyset$ then $|\beta|=|\gamma|$.
\end{enumerate}
\end{observation}

\Proof {\em 1)}\ \ By part 2) it suffices to have $\otimes_{\alpha,R'}'[e]$.
Now there is $e$ satisfying $\otimes_{\alpha,R}[e]$. We define a
function $e'\supseteq
e\restriction (\alpha\setminus R)$ such that 
$\dom(e')=\alpha\setminus R'$: just choose for $e'(\sigma)$ a club of
$\sigma$ of order type $\cf(\sigma)$. 

{\em 2)}\ \ We are given $e$ such that $\otimes'_{\alpha,R}[e]$.  We
define $e'(\delta)$ for $\delta \in\alpha$ limit by induction on
$\delta$ such that $\otimes_{\alpha, R}[e']$ holds.
\begin{description}
\item[$\mbox{\underline{Case 1:}}$]\ \ \ \  $\delta\in R$\\
We let $e'(\delta)=\delta$.
\item[$\mbox{\underline{case 2:}}$]\ \ \ \  $\delta\notin R$\\
We let $\gamma_\delta=\otp(e(\delta))$ so $\gamma_\delta<\delta$. Let
$g_\delta$ be the unique order 
preserving function from $\gamma_\delta$ onto $e(\delta)$.  Necessarily
$\gamma_\delta$ is a limit ordinal (as $e(\delta)$ has no last
element), and 
hence $e'(\gamma_\delta)$ is a well defined unbounded subset of
$\gamma_\delta$ of the order type from $R$. Let 
$$e'(\delta)=:\{g_\delta(\beta):\beta\in e'(\gamma_\delta)\}.$$
\end{description}

{\em 3)}\ \ \  Straightforward.

{\em 4)}\ \ \ \ Let $e_l '$ exemplify
${\otimes '}_{\alpha_{l+1},\alpha_l\cup R_l}$ (see (3)),\  then
$e_0'\cup e_1'$ exemplifies
 
\noindent $\otimes'_{{\alpha_2},{\alpha_0\cup R_0\cup R_1}}$ and by
part 2) we can finish. 

{\em 5)}\ \ Let $f$ be a one-to-one function from $\alpha$ onto $|\alpha|$.
Define a function $e'$, $\dom(e')=\{\delta: |\alpha|+j_{|\alpha|}\leq
\delta<\alpha$, 
$\delta$ a limit ordinal$\}$  which will satisfy
$\otimes'_{\alpha,|\alpha|+j_{|\alpha|}}$ (enough by part 2)).  

If $\delta=|\alpha|$ choose $e'(\delta)$ as an unbounded subset of $|\alpha|$
of order type $\cf(|\alpha|)$. If $\delta>|\alpha|$, let
$\xi_\delta=\min\{\beta:\delta=\sup\{f(\gamma):\gamma <\beta\}\}$,
necessarily it is a limit ordinal. Now, if
$\xi_\delta<|\alpha|$ we let $e(\delta)=\{f(\gamma):
\gamma<\xi_\delta\}$. We are left with the case $\xi_\delta=|\alpha|$.
In this case define by induction on $i<|\alpha|$, the ordinal
$\beta_i=\sup\{g(\gamma):\gamma<i\}$. As  
$\xi_\delta=|\alpha|$ clearly $(\forall i<|\alpha|)(\beta_i<\delta)$, also
clearly $(\forall j<i)(\beta_i\leq\beta_j)$. Hence $\{\beta_i:i<|\alpha|\}$
has order type $\leq |\alpha|$; as $\delta=\sup\{g(\gamma):\gamma <
|\alpha|\}$, clearly $\delta=\bigcup\limits_{i<|\alpha|} \beta_i$, so 
$e'(\delta)\stackrel{{\rm def}}{=}\{\beta_i:i<|\alpha|\}$ is as required.
Hence we have finished defining $e'$ and the proof is completed.

{\em 6)}\ \ Let $\otimes^*_{\alpha,\zeta}[e]$; by \ref{observation4:2}(4),
\ref{observation4:2}(5) w.l.o.g $\zeta=|\zeta|$. We define $f_\beta$ by
the induction on $\beta$:
\begin{quotation}
If $\beta=\gamma+1$ let $f_\beta(\gamma)=0$, $f_\beta(\epsilon)
=1+f_\gamma(\epsilon)$ for $\epsilon<\gamma$. 

If $\beta$ is a limit ordinal, first define $g_\beta$:
\[g_\beta(\gamma)= (\otp(\gamma\cap e(\beta)), f_{\min(e(\beta)\backslash
(\gamma+1))}(\gamma)).\]

\end{quotation} 
So $g_\beta$ is a one-to-one function from $\beta$ into
$\otp(e(\beta))\times\sup_{\gamma<\beta} |\gamma|$.  As usual we can well order
this set by 
($<_{lx}$ is the lexicographic order):
$$(i_1,i_2) <^*_\beta (i_2,j_2)\iff (\max \{i_1, i_2\}, i_1,i_2) <_{lx}
(\max\{j_1,j_2\},j_1,j_2).$$
Let $h_\beta$ be the one-to-one order preserving function from $(\rang
(g_\beta),<^*_\beta)$ onto some ordinal $\gamma_\beta$. Let $f_\beta=h_\beta
\circ g_\beta$; check that $\gamma_\beta$ is a cardinal.

{\em 7)}\ \ Similar proof. (We use the fact that  there is a definable
function giving for any infinite ordinal $\alpha$ a one-to-one function
from $\alpha\times \alpha$ into $\alpha$: just let $\beta_\alpha\leq
\alpha$ be the maximal limit ordinal such that $(\forall
\gamma)(\gamma<\beta \rightarrow \gamma\times \gamma < \beta)$, so for
some $n$, $\beta^n> n$, and as above we can define a one-to-one
function from $\alpha$ into
$\underbrace{\beta\times\ldots\times\beta}_{n\ \tim}$ and from
it into $\beta$).

{\em 8)}\ \ Included in the proof of (6).

{\em 9)}\ \ Like the proof of part (5).

{\em 10)}\ \ By (7).\hfill$\QED_{\ref{observation4:2}}$

\begin{lemma}
\label{lemma4:3}
Assume $\otimes_{\mu,R}$ and $\mu\geq\theta(E\times E\times {\cal
P}(A))$ and
$$
\alpha^*=\sup\{\rk^2_D(\mu_{\A},E): D\in E\}<\infty,
$$ 
and $E$,
$\A$ are as in hypothesis~\ref{nice}.0 and $\mu_{\A}$ stands for the
constant function with domain $\A$ and value $\mu$. \underline{Then}
$\otimes_{{\alpha^*},{R^*}}$, where $R\subseteq R^*\subseteq R\cup
[\mu,\alpha^*)$ and there is a $\langle Y_\sigma: \sigma\in R^*\setminus
R\rangle$ (note that $\sigma\in R^*\setminus R$ is just an ordinal not
necessarily an aleph) such that:
\begin{description}
\item[(a)\ ] $Y_\sigma$ is a non empty set of pairs $(D, \bar{\delta})$ such
that  
\[D\in E,\ \ \bar{\delta}=\langle\delta_a: a\in\A\rangle,\ \ \delta_a\in
R\] 
and $\prod\limits_{a\in\A}\delta_a/D$ has the true cofinality $\sigma$ (see
Def \ref{definition2:4}(1)).
\item[(b)\ ] The $Y_\sigma$'s are pairwise disjoint. Moreover if
$(D_\ell, \bar \delta^\ell)\in Y_{\sigma_\ell}$ for $\ell=1, 2$ and
$D_1=D_2$ then $\bar \delta^1\neq_{D_1}\bar \delta^2$.
\end{description}
\end{lemma}

\noindent{\sc Remark:}\ \ \ Instead of
$\rk^2_D(\mu_{\A}, E)<\infty$ for every $D\in E$ it is enough to assume
$\rk^2_D(\mu,E)<\infty$ for some $D\in E$ with
$$
\alpha^*=\sup\{r^2_D(\mu_{\A}, D): D\in E\mbox{ and }\rk^2_D(\mu_{\A},
D)<\infty\}.
$$
\medskip

\Proof  For $\alpha <\alpha^*$ and $D\in E$ let
$$F^D_\alpha=\{f\in {}^{\A}\mu:\rk^2_D(f)=\rk^3_D(f)=\alpha\}\ \mbox{ and
}\ A_D=\{\alpha <\alpha^*:F^D_\alpha \not=\emptyset\}.$$
So (see \ref{claim3:6}):
\begin{description}
\item[(a)]  $(\alpha,D) \mapsto F^D_\alpha$ and $D\mapsto A_D$ are well
defined (so there are such functions),
\item[(b)] $\alpha^*=\bigcup\limits_{D\in E} A_D$,
\item[(c)]  if $f_1,f_2 \in F^D_\alpha$ then $f_1/D=f_2/D$,
\item[(d)]  if $f_e \in F^D_{\alpha_e}$ for $e=1,2$ and $\alpha_1<\alpha_2$
then $f_1/D<f_2/D$. 
\end{description}
By \ref{observation4:2}(2)+(1) it suffices to prove $\otimes'_{\alpha^*,R^*\cup
(\mu+1)}$. Let $e$ be such that $\otimes_{\mu, R}[e]$ holds.
For $\delta\in (\mu,\alpha^*)$, we try to define the truth value of
$\delta\in R^\ast$ and $e'(\delta)$ such
that $\otimes'_{{\alpha^*},R^*\cup (\mu+1)}[e']$ holds. We make three
tries; an easier case is when the definition gives an unbounded subset of
$\delta$ of order type $< \delta$: decide $\delta\notin R^*$ and
choose this set as $e'(\delta)$. If not, we assume we fail and
continue, and if we fail in all three of them then we decide
$\delta\in R$, and choose $Y_\delta$.
\medskip

\noindent \underline{First try:}\ \ $e_1(\delta)\stackrel{{\rm
def}}{=}\{\sup(\delta 
\cap A_D): D\in E\mbox{ and } \sup(\delta\cap A_D)<\delta\}$.\\
Clearly this set has cardinality $<\theta(E)\leq \mu\leq|\delta|$.
So the problem is that it may be bounded in $\delta$. In this case, by (b)
above, for some $D\in E$, $\delta=\sup(\delta\cap A_D)$.
\medskip

\noindent \underline{Second try:}\ \ Let $E_\delta=\{D: \delta=\sup(A_D\cap
\delta)\}$.\\
So we can assume $E_\delta\neq\emptyset$.  For each $D_0\in E_\delta$ let
$$
\begin{array}{ll}
E_\delta(D_0)=\{ D\in E:& D_0\subseteq D\mbox{ and there is }f\in
{}^{(\A)}\mu\mbox{ such that }f/D\mbox{ is }\\
\ & \leq_D-\lub\mbox{ of }\{g/D:g\in F^{D_0}_\beta,\mbox{
for some }\beta\in A_{D_0}\cap\delta\}\}. 
\end{array}
$$
For $D\in E_\delta(D_0)$ let 
\[F^{{(D_0},D)}_\delta=\{f: f/D\mbox{ is }\leq_D-\lub \mbox{ of }\{g/D:g\in
F^{D_0}_\beta\mbox{ for some }\beta \in A_{D_0}\cap\delta\}\}.\] 
For $f\in F_\delta^{({D_0},D)}$, let $B(f)=\{a\in\A:f(a)$ a limit
ordinal$\}$.  
Now $B(f)\in D$ by \ref{obs1:2}(8) because of the
assumptions $f$ is a $\leq_D-\lub$ and $A_{D_0}\cap \delta$ is
unbounded in $\delta$ and let
$$
H(f)=\{g\in {}^{\A}\mu:(\forall a\!\in\! B(f))(g(a)\!\in e(f(a)))\ \&\
(\forall a\!\in\!\A\setminus B(f))(g(a)=0)\}.
$$ 
(remember $e$ is a witness for $\otimes_{\mu, R}$). Note that
\begin{description}
\item[(e)] for $D_0\in E_{\delta}$, $E_{\delta}(D_0)$ is not empty [by
\ref{claim3:12}(1)],
\item[(f)] if $D_0\in E_{\delta}$, $D\in E_{\delta}(D_0)$ then
$F_{\delta}^{(D_0,D)}\neq\emptyset$,
\item[(g)] if $f_1, f_2\in F_{\delta}^{(D_0,D)}$ then $f_1/D= f_2/D$\\ 
(as $<_D$-\lub\ is unique $\mod D$) hence
\item[(g)$'$] $H_D(f_1)=H_D(f_2)$ where $H_D(f)=\{g/D: g\in H(f)\}$,
\item[(h)] if $D_0\in E_\delta$ and $f\in F_{\delta}^{(D_0,D)}$ and
$g\in H(f)$ \underline{then} $g<_D f$,
\item[(i)] if $f\in F_{\delta}^{(D_0,D)}$ then $|H(f)|\leq \prod\limits_{i\in
B_2(f)} e(f(i))$.
\end{description}
We now define, for $f^*\in F_{\delta}^{(D_0,D)}$ a function $h=h_{\delta,
f^*, D_0, D}$ from $H(f^*)$ to $\delta$:
\[h(g)=\min\{\alpha<\delta: \alpha\in A_{D_0}\mbox{ and there is }
f\in F^{D_0}_{\alpha}\mbox{ such that }\neg(f<_D g)\}.\]
(by the way, equivalently for every $f\in F^{D_0}_\alpha$). Now
\begin{description}
\item[(j)] $h(g)$ is well defined.
\end{description}
[Why? As otherwise $g$ exemplify $f^\ast$ is not a $<_D$-\lub\  of
$\bigcup\{F_{\alpha}^{D_0}: \alpha\in A_{D_0}\cap\delta\}$ ($g$ is a smaller
$<_{D_0}$-upper bound).]

Also
\begin{description}
\item[(k)] $\rang(h_{\delta,f^*,D_0,D})$ is unbounded below $\delta$.
\end{description}
[Why? For every $\alpha<\delta$ there is $\beta\in A_{D_0}\cap (\alpha,
\delta)$. Choose $f\in F^{D_0}_{\beta}$, and define $g\in {}^{(\A)}\ord$ by
$g(a)=\min(e(f^\ast(a))\setminus (f(a)+1))$. Now $g\in H(f^*)$ and $h(g)>\beta
>\alpha$.]  
\begin{description}
\item[(l)] $\rang(h_{\delta,f^*,D_0,D})$ does not depend on $f^*$, i.e. is the
same for all $f^*\in F^{(D_0, D)}_\alpha$, [by (g)], and we denote it by
$t_{\delta}(D_0,D)$.
\item[(m)] $h$ is a nondecreasing function from $(H(f^*),\leq_D)$ to $\delta$.
\end{description}
If $|E|$ is an aleph, choose a fixed well ordering $<^*_E$ of $E$ and if for
our $\delta$ for some pair $(D_0,D)$ we have $\otp(t_{\delta}(D_0,D))<\delta$,
choose $e(\delta)=t_{\delta}(D_0,D)$ for the first such pair; but this
assumption on $|E|$ is not really necessary: 

If for some $D_0\in E_{\delta}$, $D\in E_{\delta}(D_0)$ we consider (i.e. if
this set is O.K., we choose it; easily it is an unbounded subset of $\delta$):
\[\begin{array}{ll}\{t_{\delta}(D_0,D): & D_0\in E_{\delta}, D\in E_{\delta}
(D_0) \mbox{ and for any other such }({D'}_0, D')\\
\ & \mbox{ we have }\otp(t_{\delta}(D_0,D))\leq\otp(t_{\delta}({D'}_0,D'))\}.
\end{array}\]
This is an indexed family of unbounded subsets of $\delta$, indexed by a
subset of $E\times E$, all of the same order type, which we call $\beta^*$. As
we know $\theta(E\times E)$ is a cardinal $\leq \delta$. By observation
\ref{lemma4:3}A below it is enough to have $\beta^*<|\delta|$. 
\medskip

\noindent{\sc Observation} \ref{lemma4:3}A:\ \ \ 1) If ${\cal B}=\{B_c: c\in
C\}$, $B_c\subseteq\delta=\sup (B_c)$, $\delta$ a limit ordinal,
$\sup\limits_{c\in C}\otp(B_c)<|\delta|$ and $\theta(C)\leq\delta$
\underline{then} 
we can define (uniformly) from  $\delta$, ${\cal B}$ an unbounded subset $B$
of $\delta$ such that $\otp(B)<\delta$. 

\par \noindent
2) If in part (1) we omit the assumption $\delta=\sup(B_c)$, the conclusion
still holds provided that $\neg (\ast)$ where 
\begin{description}
\item[$(*)$] for every $\alpha^*<\delta$ we have
$\delta=\sup\bigcup\{B_c:\alpha^*< \sup (B_c) \leq \delta\}$. 
\end{description}
\medskip

\Proof For each $i<\delta=:\sup\limits_{c\in C}\otp (B_c)$ define
$$
B^*_i=\{\gamma: \mbox{ for some }c\in C\mbox{ the ordinal }\gamma\mbox{ is the
}i\mbox{-th member of }B_c\}. 
$$
So there is a partial function from $C$ onto $B^*_i$, hence $\otp(B^*_i)<
\theta(C)\leq \delta$. So if for some $i$, the set $B^*_i$ is unbounded in
$\delta$ then let $B$ be $B^*_i$ for the minimal such $i$. If there is no such
$i$ then let $\gamma^*=\sup\{\otp(B_c): c\in C\}$. By assumption $\gamma^*<
|\delta|$, let  
$$
B=: \{\sup(B^*_i): i<\gamma^*\},
$$
so $B\subseteq\delta$, and $|B|\leq^* |\gamma^*|$ hence $|B|\leq |\gamma|$
hence $\otp (B)<|\delta|\leq \delta$ and for every $\beta<\delta$ for
some $c\in 
C$ there is $\gamma\in B_c\setminus\beta$, so for $i=\otp(B_c\cap \gamma)$ we
have: $\gamma\in B^*_i$ hence $\gamma\leq \sup (B^*_i)\leq \sup (B)$,
so $\sup(B)= \delta$. 
So $B$ is an unbounded subset of $\delta$ of order type $<\delta$.

\par \noindent
2) Clearly 
$$
B=\{\sup (B_c): c\in C\mbox{ and }\sup(B_c)<\delta\}
$$
has order type $<\theta(C)\leq \delta$, so if $\delta=\sup(B)$ we are done; if
not let $\alpha^*=\sup(B)<\delta$ and $C'=\{c\in C:\delta=\sup(B)\}$,
${\cal B}'=:\{B_c: c\in C'\}$ are as in the assumption of
\ref{lemma4:3}A(1). 
\hfill\qed$_{\ref{lemma4:3}A}$
\medskip

\par
\noindent
{\sc Continuation of the proof \ref{lemma4:3}:}

\noindent \underline{Third try:}\ \  We are left with the case that every
$t_{\delta}(D_0,D)$, when well defined, has order type $\delta$. Continue with
our $f^*\in F_{\delta}^{(D_0,D)}$, $h=h_{\delta,f^*,D_0,D}$. For each $g\in
H(f^*)$, we know that there are $A\in D^+$ and $f\in F_{h(g)}^{D_0}$ such that
$g\restriction A\leq f\restriction A$. Now turn the table: for $A\in D^+$ let 
$$
H(f^*,A)=\{g\in H(f^*): \mbox{ for some (=all) }f\in F_{h(g)}^{D_0}
\mbox{ we have }g{\leq}_{A} f\}
$$ 
(clearly the choice of $f$ is immaterial : some, all are the same).  Now:\\
$(H(f^*,A), {\leq}_{D+A})$ is mapped by $h$ into $t_{\delta}(D_0,D)\subseteq
\delta$; moreover for $g'$,$g''\in H(f^\ast,A)$ we have: $h(g')<h(g'')
\Rightarrow g'<_{D+A} g''$ (otherwise the minimality of $h(g'')$ is
contradicted). So $(H(f^*,A),\leq_{D+A})$ has the true cofinality
$\otp(t_{\delta}^1 (D_0,D,A))$, where $t_{\delta}^1(D_0,D,A)=:(\rang(h
\restriction H(f^*,A)))$. Now by observation \ref{lemma4:3}A if we do
not succeed, for some
$\beta^*<\delta$ we have 
\begin{description}
\item[$(*)$] if $(D_0,D,A)$ above and $t^1_\delta(D_\delta,D,A)\not\subseteq
\beta^*$ \underline{then} $\otp(t^1_\delta(D_0,D,A))=\delta$. 
\end{description}

So we assume that for some $\beta^* <\delta$ we have $(*)$. Without loss of
generality $\beta^*$ is minimal. Consider such a triple $(D_0,D,A)$ and the
appropriate $f^*= f^*_\delta$. Notice that $H(f^*,A)$ is cofinal in $(H(f^*),
<_{D+A})$, (we use clause (d) of \ref{claim3:12}). Now consider whether the
quadruple 
$$
\bar x=(D_0,D,A,\langle\otp(e(f^*_\delta(a))): a\in\A\rangle/(D+A))
$$
was considered by some earlier $\delta^{\otimes}$, if so, choose minimal
$\delta^{\otimes}$, and we shall finish by the observation \ref{lemma4:3}B
below. To stress the dependency  on $\delta$ we may write
$H_\delta(f^*_\delta, A)$, and 
define $\Col_\alpha$ for $\alpha<\mu$ such that  $\alpha\in \dom(e)$
onto $e(\alpha)$ as the unique order preserving function from
$\otp(e(\alpha))$ onto $e(\alpha)$. 
Let us define for $g\in H(f^*_\delta, A)$, the function $\Col_g$ with
domain $A^*$, $\Col_g(a)$ being $\otp(e(f^*_\delta(a))\cap g(a))$ if
$f^*_\delta(a)$ is a limit ordinal, zero otherwise. Of course $\Col_g$ depends
on the choice of $f^*_\delta$ but $\Col_g/D$ does not, and let 
$H_\delta(f^*_\delta, A)=
\{\Col_g: g\in H(f^*_\delta, A)\}$, so
$H_\delta(A)=:(H_\delta(f^*_\delta, A)/(D+A), <_{D+A})$ has
order type $\delta$ and $H_\delta(A)$ is cofinal in $\prod\limits_{a\in \A}
\otp f^*_\delta(a)/(D+A)$, and $H$ and $H_\delta(A)$ were definable from
$\delta$, $(D_0, D, A)$ in a fixed way. 

Similarly for $H_{\delta^\otimes}(A)$. So observation \ref{lemma4:3}B below
give us a definition of an unbounded subset $Z_0(D_0,D,A,\delta,
\delta^\otimes)$ of $H_\delta(A)$ of order type $\delta^\otimes$, hence also
of $\otp(H_\delta(A),<_D)=\delta$ which we call $Z(D_0,D,A,\delta,
\delta^\otimes)$. So
$$
\{Z(D_0, D, A, \delta, \delta^\otimes): (D_0, D, A)=\bar x\restriction 3\mbox{
for some }\bar x \in {\mathfrak C}\} 
$$
is a family of unbounded subsets of $\delta$ of order type $\delta^\otimes$ so
by observation \ref{lemma4:3}A we are done. 

\par \noindent
\underline{Fourth try:}

All previous ones failed, in particular there is no $\delta^\otimes$ as above.
We put $\delta$ in $R^\ast$, and let $Y_\delta$ be the set of quadruple
$\langle D_0, D, A, \langle \otp(e(f^*_\delta(a)): a\in A^*\rangle/D\rangle$
as above (the last one is uniquely determined by the earlier ones).
\medskip

\noindent{\sc Observation \ref{lemma4:3}\ B}\ \ \  If for $l=1,2$ we have
$Y_l\subseteq {}^{(A^*)}\mu$ is cofinal in
$\prod\limits_{a\in\A}\alpha^l_a /D$ 
and $[g_1, g_2\in Y_\ell \Rightarrow (g_1=_D g_2)\vee (g_1<_D g_2)\vee
(g_2 <_D g_1)]$ and $(Y_\ell/D, <_D)$ is well ordered of order type
$\delta_\ell$ and $\delta_1<\delta_2$ and $\{a\in \A:
\alpha^1_a=\alpha^2_a\}\in D$ (and $\alpha^\ell_a>0$ for simplicity)
\underline{then} from the parameters $D$, $Y_\ell/d$, $\langle
\alpha^\ell_a: a\in \A\rangle/D$ for $\ell=1,2$ we can (uniformly)
define one of the following:
\begin{description}
\item[(a)] $Y\subseteq Y_1/D$ cofinal in $\prod\limits_{a\in \A}
\alpha^1_a/ D$ and a function $h$ from $Y$ into $\delta_2$,
$<_D$-increasing such that $\delta_2=\sup \rang(h\restriction Y)$.
\item[(b)] cofinal $X\subseteq \delta_2$ of order type $<\theta({\cal P}(\A))$.
\end{description}
\medskip

\Proof For $g_1\in Y_1$ let 
$$
\begin{array}{ll}
K(g_1)=: \{ g_1\in Y: & g_1\leq_D g_2\mbox{ and there is no }g'_2\in
Y_2 \mbox{ such that }\\
\ & \qquad \qquad g'_2<_D g_2\mbox{ and } g_1\leq_D g'_2\},
\end{array}
$$
$$
K_D(g_1)=: \{g_2/D: g_2\in K(g_1)\},
$$
$$
K'_D(g_1/D)= \bigcup \{K_D(g'_1): g'_1\in Y_1\mbox{ and }g'_1=_D
g_1\}.
$$

${\bf k}(g_1/D)$ is (if exists) the $\leq_D$-minimal $g_2/D\in Y_2$
such that 
$$
g'_2/D\in K'_D(g_1/D) \Rightarrow g'_2/D\leq_D g_2/D.
$$
Now $K(g_1)$ is a subset of $Y_2$, $K_D(g_1)$ is a subset of $Y_2/D$
and $g'_1=_D g^{\prime\prime}_1 \Rightarrow K_D(g'_1)=
K_D(g^{\prime\prime}_1)$ hence $K_D(g_1)= K'_D(g_1/D)$.

\par \noindent
\underline{Case 1:} for some $g_1/D\in Y_1/D$, the set
$K^{\prime}_D(g_1/D)$ is 
unbounded in $Y_2/D$.

We can choose the $<_D$-minimal such $g_1/D$, so $K'_D(g_1/D)$ is a
well defined cofinal subset of $Y_2/D$ (remember $Y_2/D$ is $<_D$-well
ordered) and easily $\otp(K'_D(g_1/D), <_D)< \theta(\{A/D: A\in
D^+\})\leq \theta({\cal P}(\A))$.

\par \noindent
\underline{Case 2:} not case 1.

So ${\bf k}$ is a  well defined function from $Y_1/D$ into $Y_2/D$,
and easily $g'_1/D\leq_D g^{\prime\prime}_1/D \Rightarrow {\bf
k}(g'_1/D)\leq_D {\bf k}(g^{\prime\prime}_1/D)$.

Let 
$$
\begin{array}{ll}
Y=: \{g_1/D: & g_1/D\in Y_1/D\mbox{ and }\\
\ & [g'_1/D\in Y_1/D\ \&\ g'_1/D
<_D g_1/D\Rightarrow {\bf k}(g'_1/D)<_D {\bf k}(g_1/D)]\},
\end{array}
$$
and let $h(g_1/D)=\otp(\{g_2/D\in Y_2: g_2/D<_D {\bf k}(g_1/D)\})$, so
$h: Y\rightarrow \delta_2$. Now $Y$, $h$ are as required. 
\hfill\qed$_{\ref{lemma4:3}B}$\\
\null\hfill\qed$_{\ref{lemma4:3}}$

\begin{claim}
\label{claim4:4} 
$[DC+\oplus_\kappa[E]]$\ \ \ \ Assume $\mu>\A=\cf(\mu)>\aleph_0$, $\mu>
\theta(E)+\kappa^+$, $\otimes_{\mu,R}$ and
$\lambda=|\rk^2_D(\mu,E)|$. \underline{Then} $\lambda<\aleph_\gamma$, for some
$\gamma<\theta(E\times E\times {\cal P}(\A)\times |R|^{|\A|})$. 
\end{claim}

\Proof By \ref{lemma4:3} and \ref{observation4:2}(10).
\hfill$\QED_{\ref{claim4:4}}$  

\begin{observation}
\label{obs6:1}
\begin{enumerate}
\item $\theta(\bigcup\limits_{x\in X} A_x)\leq (\theta(X)+\sup_{x\in
X}\theta(A_x))^+$. 
\item If $\lambda$ is regular, $\lambda\geq\theta(X)$,
$\lambda\geq\theta (A_x)$ (for $x\in X$) \underline{then}
$\theta (\bigcup\limits_{x\in X} A_x)\leq\lambda$. So e.g. $\theta(A\times
A)=\theta(A)$ if $\theta (A)$ is regular and $\theta(A\times B)\leq
(\theta(A)+\theta (B))^+$, and $\theta (A\times B)=\max
\{\theta (A), \theta (B)\}$ (or all three are finite) when the later
is regular. 
\end{enumerate}
\end{observation}

\Proof {\em 1)}\ \ As $A\leq^* B\Rightarrow \theta(A)\leq \theta(B)$
and $\bigcup\limits_{x\in X}A_x\leq^* \bigcup\limits_{x\in
X}(\{x\}\times A_x)$ clearly 
without loss of generality $\langle A_x: x\in X\rangle$ is pairwise disjoint.
Let $A=\bigcup\limits_{x\in X} A_x$ and $f: A\stackrel{\rm onto}
{\longrightarrow}\alpha$, so $\gamma_x\stackrel{\rm
def}{=}\otp(\rang(f\restriction A_x)) 
<\theta(A_x)\leq\sup_{x} \theta(A_x)$ and call the latter $\lambda$. So
$\alpha=\bigcup\{B_{\beta}: \beta<\lambda\}$ where 
$$
\begin{array}{rr}
B_{\beta}=:&\big\{\alpha:\mbox{ for some } x\in X\mbox{ we have }\alpha\in
\rang(f\restriction A_x)\mbox{ and}\\
\ &\beta=\otp(\alpha\cap\rang(f\restriction A_x))\big\}.
\end{array}
$$    
Let 
$$
\begin{array}{ll}
A^{\beta}=:\{ a\in \bigcup\limits_{x\in X} A_x : & \mbox{ for the
}x\in X\mbox{ such that } a\in A_x\mbox{ we have:}\\
\ & \beta\mbox{ is the order type of }\{f(b): b\in A_x\mbox{ and
}f(b)<f(a)\}\}. 
\end{array}
$$
Let $E$ be the relation on $A$ defined by $a\ E\ b\iff\bigvee\limits_{x\in X}
\{a,b\}\subseteq A_x$ so $E$ is an equivalence relation, and the $A_x$ are the
equivalence classes. 

Now $A^\beta\cap A_x$ has at most one element, so $f\restriction A^\beta$
respects $E$, so $f$  induces a function from $A^\beta/E$ onto $B_\beta$. So 
$$
|B_\beta|<\theta(A^\beta/E)=\theta(\{A_x: A^\beta\cap
A_x\neq\emptyset\})\leq \theta (X).
$$
Hence $\otp(B_{\beta}) < \theta(X)$. So $|\alpha|\leq \sum\limits_{\beta<\lambda}
|B_{\beta}|\leq \lambda\times\theta(X)$ and 
$$
\theta(\bigcup\limits_{x\in A}
A_x)\leq (\lambda+\theta(X))^+=[\sup_{x\in X}
\theta(A_x)+\theta(X)]^+.
$$

{\em 2)}\ \ We repeat the proof above. Clearly $\gamma_x<\theta(A_x)\leq
\lambda$, so $x\mapsto\gamma_x$ is a function from $X$ to $\lambda$, so as
$\lambda$ is regular, necessarily $\gamma^\ast=(\sup\limits_{x\in X}\gamma_x)
<\lambda$. So $B_\beta$ is defined for $\beta <\gamma^\ast$ only and again
$|B_\beta|<\theta (X)\leq\lambda$. Hence $\otp (B_\beta) <\lambda$ and hence by
``$\gamma^\ast<\lambda$, $\lambda$ regular'' we have $\beta^\ast=\sup
\limits_{\beta < \gamma^\ast}\otp (B_\beta) <\lambda$ and we finish as above.\\
The ``e.g.'' follows (with $X=A$, $A_x=A\times\{x\}$, $\lambda=\theta(A)$).
\hfill$\QED_{\ref{obs6:1}}$

\begin{theorem}
\label{lemma4:5}  
$[DC]$\ \ \ Assume
\begin{description}
\item[(a)] $\mu >\cf(\mu)=\kappa>\aleph_0$,
\item[(b)] $|\bigcup\limits_{\alpha<\mu}{\cal P(\alpha)}|=\mu$,
\item[(c)] $\oplus_\kappa(E)$. (I.e. $E$ is a nice family of filters
on $\kappa$.)
\end{description}
\underline{Then}
\begin{description}
\item[($\alpha$)]  $2^\mu$ is an aleph,
\item[($\beta$)]  $\mu^+$ is regular, and $\otimes_{2^\mu,\{\delta<
2^\mu:cf(\delta)>\mu\}}$; if $\mu=\aleph_\gamma$ then $\otimes_{2^\mu,R}$,
where $R\cap\mu=\reg\cap \mu$ and $|R\setminus\mu|\leq\theta(\gamma^\kappa
\times|E|^2)$,
\item[($\gamma$)] $\mu<\lambda\leq 2^\mu\Rightarrow \lambda$ not measurable,
\end{description}
\end{theorem}

\noindent{\sc Remark \ref{lemma4:5}A:}  Instead of (b), $\otimes_\mu$ 
suffices, if $2^\mu$ is replaced by $\mu^\kappa$ and $|E|$ is well ordered. 
\medskip

\Proof By the proof of \ref{claim3:11} (and clause (b) of the assumption)
there is $F$ exemplifying $2^\mu\in {\ts}_{{D^{\rm bd}_\kappa}}(\mu)$.

Let 
$$
A_D=:\{\alpha: F^D_\alpha\neq\emptyset\}\quad
\mbox{ and }\quad F^D_\alpha=:\{f\in
F:\rk^2_D(f)=\rk^3_D(f)=\alpha\}.
$$ 
Note that $F^D_\alpha$ is a singleton
denoted by $f^D_\alpha$, so $F=\bigcup\limits_{D\in E} F_D$ where
$F_D=\{f^D_\alpha:\alpha\in A_D\}$, so the $F_D$ are not neccessarily
pairwise disjoint. By clause (b) we can prove that as $\kappa<\mu$ also
$|{\cal P}(\kappa)|<\mu$ hence $|{\cal P}({\cal P}(\kappa))|<\mu$ and hence
$|E|<\mu$ so $|E|$ is an aleph. As $|E|$ is an aleph we can well order $E$
(say by $<_E$) and hence can well order $F$:
\[f <^*g\mbox{\ \ \ iff\ \ \ }\alpha(f)<\alpha(g)\mbox{\ \ or\ \
}\alpha(f)=\alpha(g)\mbox{\ \ and\ \ }D(f)<_E D(g)\] 
where $\alpha(f)$ is the unique $\alpha$ such that for some $D$ we have
$f=f^D_{\alpha(f)}$ and $D(f)$ is the $<_E$-first $D$ which is suitable. As
$F$ is well ordered, $|F|$ is an aleph, but $F$ was chosen by
\ref{claim3:11} such that $|F|=2^\mu$, so we have proved clause ($\alpha$).

If $\lambda\in(\mu,2^\mu]$, there is $F'\subseteq F$ of cardinality $\lambda$
so as in \ref{claim2:8}, $\lambda$ is not measurable so clause $(\gamma)$
holds.  

Now we shall deal with clause $(\beta)$.\\ 
By assumption (b) clearly $\otimes_{\mu}$, so \ref{lemma4:3} applies and we
get $\otimes_{2^\mu, R^*}$ for $R^*$ as there. Now  suppose that for
$\langle\sigma_i : i<\kappa\rangle\in {}^{\kappa}(\reg\cap \mu)$, and $D\in 
E$, $\prod\limits_{i<\kappa}\sigma_i /D$ has the true cofinality $\sigma$. If
for some $A\in D^+$, $\sup_{i\in A}\sigma_i<\mu$, $\prod\sigma_i/(D+A)$ is
well ordered by (b), so $\sigma<\mu$. So assume $\tlim_D \sigma_i =\mu$
(i.e. $\mu=\lim\sup\limits_D \langle \sigma_i: i<\kappa\rangle=\lim
\inf\limits_D\langle \sigma_i: i<\kappa\rangle$),
then
$\prod\limits_i\sigma_i /D$ is $\mu^+$-directed so $\cf(\sigma)>\mu$, hence
$\cf(\sigma)\geq \mu^+$. We conclude $R^*\cap\mu^+\subseteq R$, hence, by
\ref{lemma4:3}, $\otimes_{\mu^+,R}$, so, by \ref{observation4:2}(8), $\mu^+$ is
regular. This gives the first phrase in clause $(\beta)$, the second is
straightforward.\hfill$\QED_{\ref{lemma4:5}}$ 
\medskip

\noindent{\sc Discussion}\ \ref{lemma4:5}B:
\begin{enumerate}
\item If $|E|+|{\cal P}(\kappa)|$ is an aleph and the situation is as in
\ref{lemma4:5}, then we can choose
$A\subseteq\theta(2^{\mu})$ which codes the relevant instances of
$\otimes_{\mu^+,R}$, $\otimes_{\lambda,R^*}$ (or
$\otimes_{2^{\mu},R^*}$) and then work in $L[A]$ and apply theorems on
cardinal arithmetic (see in [Sh-g]) as in \ref{claim4:4} (so we can ask
on weakly inaccessible etc.), but we have to translate them back to
\v, with $A$ ensuring enough absoluteness.
\item If $|E|+|{\cal P}(\kappa)|$ is not an aleph, we can force this
situation not collapsing much, see \S\ref{par5}.
\end{enumerate}

\begin{definition}
\label{def4:7}
For an ordinal $\delta$ let 
$$ID^1_{\delta}=\{A\subseteq\delta:
\otimes_{\delta,{\delta\setminus A}}\mbox{ holds }\},$$
$$ID^1_{\delta,\alpha}=\{R\subseteq\delta: R\setminus\alpha\in
ID^1_{\delta}\},$$ 
$$ID^2_{\delta,R}=\{A\subseteq\delta: \mbox{ there is a
function } e \mbox{ with domain } A \mbox{ such that}$$
$$\ \ \ \ \ \ \ \ \ \ \ \ \ \ \ \mbox{the requirement in }
\otimes_{\delta, R}[e]\mbox{ holds for } \alpha\in A\},$$
$$ID^2_{\delta,R,\alpha}=ID^2_{\delta,R\cup \alpha}.$$
Omitting $R$ means $\reg\cap\delta$.
\end{definition}
  
\begin{claim}
\label{claim4:8}
\begin{enumerate}
\item $ID^1_{\delta,0}=ID^1_{\delta}\mbox{  and  }\delta\in\
ID^2_{\delta,R}\iff\otimes_{\delta,R}$.
\item $ID^1_{\delta,\alpha}$, $ID^2_{\delta}$ are ideals of subsets of
$\delta$.
\item Assume $\alpha_1 <\alpha_2 <\alpha_3\in\ord$ and $ID^2_{\delta_{l+1},R,
\alpha_l}$ for $l=1,2$. Then $ID^2_{\alpha_3,R,\alpha_1}$.
\item $\delta\in ID^2_{\delta}$ if $AC_{\delta,<|\delta|}$.
\item \mbox{[$AC_{|\alpha|}$]} $ID^2_{\delta,R,\alpha}$ is
$|\alpha|^+$-complete (see \ref{lemma4:3}A, we need to choose the witnesses
$e_1$).
\item $A\in ID^1_{\delta} \iff A\in ID^2_{\delta,{\delta\setminus A}}$
\item \mbox{[$AC_{|\alpha|}$]} $ID^1_{\delta,\alpha}$ is
$|\alpha|^+$-complete.  
\end{enumerate}
\end{claim}

We think that those ideals are very interesting.

\section{The successor of a singular of uncountable cofinality}
\label{par5}
\begin{claim}
\label{claim5:1} 
Assume $\cf(\mu)=\kappa>\aleph_0$, $\mu=\bigcup\limits_{i<\kappa}\mu_i$ where $\mu_i$
increasing continuous, $\mu$ an aleph ($\mu_i$ may be merely an ordinal).
If $f^*\in {}^\kappa\ord$, $f^*(i)=|\mu_i|^+$ \underline{then}
$\|f^*\|_{D_\kappa^{\rm bd}}\geq \mu^+$ (hence $\rk^2_D(f^*)
\geq \mu^+$ when $D\in E$ extend $D^{\rm bd}_\kappa$).
\end{claim}

\Proof Let $\alpha<\mu^+$ and we shall prove
\begin{description}
\item[$(*)$] there is $\langle f_\beta:\beta\leq \alpha\rangle$, $f_\beta \in
\prod\limits_{i<\kappa} f^*(i)$,  such that $\beta<\gamma <\alpha\Rightarrow
f_\beta<_{D^{\rm bd}_\kappa}f_\gamma$ 
\end{description}
Let $g$ be a one-to-one function from $\alpha$ into $\mu$.  For every
$\beta\leq\alpha$ let 
$$
f_\beta(i)=:\otp\{\gamma <|\mu_i|^+ :g(\gamma)<\beta\}.
$$
\hfill$\QED_{\ref{claim5:1}}$

\begin{observation}
\label{observation5:2}  
[$DC+AC_\kappa+\oplus_\kappa[E]$, $E$ normal, $\kappa$ an aleph]\ \ \ Assume
\begin{description}
\item[(a)\ ] $\langle\mu_i:i\leq\kappa\rangle$ is strictly increasing
continuous sequence of alephs,
\item[(b)\ ] $f\in\prod\limits_{i<\kappa}(\mu_i^++1)$, $D_0\in E$ and
$\rk^2_{D_0}(f,E)=\rk^3_{D_0}(f,E)=\mu^+$,
\item[(c)\ ] ${\it Tw}_D(\mu_i)<\mu$ for $D\in E$ and $i<\kappa$.
\end{description}
\underline{Then} $\{i:f(i)=\mu^+_i\} \in D$.  
\end{observation}

\Proof Otherwise without loss of generality $(\forall i)(0<f(i)<\mu_i^+)$ hence $(\forall
i)(\cf[f(i)]<\mu_i)$. By $AC_\kappa$ there is $\langle g_i:i<\kappa\rangle$,
$g_i$ is a one-to-one function from $f(i)$ into $\mu_i$. By the normality and
the possibility to replace $D$ by $D+A$ (for any $A\in D^+$), for every
$f_1\in\prod\limits_{i<\kappa} f(i)$ and $D_1$, $D_0\subseteq D_1\in E$, for
some $A\in D^+_1$, and $j<\kappa$ we have $(\forall i\in B)(g_i\circ f_1
(i)<\mu_j)$, so we conclude 
$$
\mu^+=\bigcup \{W_{D,j,A}: D_0\subseteq D\in E,j<\kappa, B\in D\},
$$ 
$\langle F_{D,j,A,\alpha}: \alpha\in W_{D,j,A}\rangle$ is $<_D$-increasing and
$D$-smooth, where
$$
\begin{array}{rr}
F_{D,j,A,\alpha}= &\big\{f_1:f_1\in\prod\limits_{i<\kappa}f(i), [i\in
A\Rightarrow g_i\circ f_1(i)<\mu_j]\quad\mbox{ and}\\
\ &\rk^2_D (f_1,E)=\rk^3_D (f_1,E)\big\},
\end{array}
$$ 
$$
W_{D, j, A}=: \{\alpha: F_{D, j, A, \alpha}\neq \emptyset\}.
$$
As ${\it Tw}_D(\mu_j)\leq\mu$, clearly $|W_{D,j,A}|\leq\mu$. So 
$$
\mu^+=\bigcup\{ W_{D,j,A}:D_0\subseteq D\in E,j<\kappa,A\subseteq\kappa\}
$$ 
with each set of cardinality $<\mu$, as $\theta(E)<\mu$, $\theta({\cal
P}(\kappa))<\mu$ (as in \ref{lemma4:3}), we get a contradiction.
\hfill$\QED_{\ref{observation5:2}}$  

\begin{conclusion}
\label{conclusion6:3}
Under the hypothesis and the assumptions (a)-(c) of \ref{observation5:2}, for
the  $f^\ast\in {}^\kappa\ord$ defined by $f^\ast (i)=:\mu^+$, for every
$D\in E$ we have: $\rk^2_D (f,E)=\rk^3_D (f, E)=\mu^+$.
\end{conclusion}

\begin{theorem}
\label{theorem5:3}
[$DC + AC_\kappa + \oplus_\kappa(E) + \aleph_0<\kappa=\cf(\kappa)$, $E$
normal or just weakly normal as witnessed by $\i$]\ \ \ Assume (a)--(c) of
\ref{observation5:2}  and
\begin{description}
\item[(d)\ ] $S=\{i<\kappa:\cf(\mu_i^+)=\lambda\}$ is stationary (see
\ref{definition3:?}\ (2)), moreover
$\i^{-1}(S)\in D_0\in E$ for some $D_0$.
\end{description}
\underline{Then} $\cf(\mu^+_\kappa)=\lambda$
\end{theorem}

\Proof By \ref{observation5:2} + \ref{claim3:12} we know
$\rk^3_D(\langle\mu_i:i< 
\kappa\rangle, E)\geq\mu^+$ so by \ref{claim3:3}+\ref{claim5:1} there are
$f\in\prod\limits_{i<\kappa}(\mu^+_i +1)$ and $D\in E$ such that $S\in D$,
$\rk^2_D(f)=\rk^3_D(f)=\mu^+$. Without loss of generality each $\mu_i$
is singular; apply 
\ref{observation5:2}. By $AC_{\kappa}$ there is $\langle e_i:i<\kappa\rangle$,
$e_i\subseteq\mu^+_i$ cofinal, $\otp(e_i)=\cf(\mu_i^+)$. Now as in the
proof of \ref{lemma4:3} we know $\{\rk^3_D(h):D\in E, h\in
\prod\limits_{i<\kappa}e_i\}$ is cofinal in $\mu^+$, hence 
\[\{\rk^3_D(h):D\in E, h(i) \mbox{ the $\gamma$-th member of $e_i$ for $i\in
S$, 0 otherwise}\}\]
is a cofinal subset of $\mu^+$.\hfill$\QED_{\ref{theorem5:3}}$

\begin{conclusion}
\label{conclusion5.4} [as \ref{theorem5:3}+(a), (b), (c) of
\ref{observation5:2}] 
If $\langle\mu_i:i\leq\kappa\rangle$ is a strictly increasing  continuous
sequence of alephs \underline{then} for at most one $\lambda>\kappa$ the set
$\{i<\kappa:\cf(\mu^+_i)=\lambda\}$ is stationary 
\end{conclusion}
\bigskip

\par \noindent
{\sc Remark:} We can prove that for all $D\in E$,
$\prod\limits_{i<\kappa}\cf (\mu^+_i)/D$ has the same cofinality.
\bigskip

\noindent
{\sc Question 5.6:}
Under the hypothesis and assumptions of \ref{theorem5:3} can we in addition
have a stationary $S\subseteq \kappa$ such that $\langle\cf(\mu^+_i): i\in
S\rangle$ is with no repetition?
\bigskip

\section{ Nice $E$ exists}
\label{par8}
\begin{observation}
\label{obs8:1}
For an aleph $\lambda$ and a set $A^*$, $\lambda^{|A^*|}$ is equal to 
$$
\sum_{\beta<\theta(A^*)} \lambda^{|\beta|}\times |\beta|^{|A^*|}
$$
provided that $\aleph_0\leq |A^*|$ (otherwise replace $|\beta|^{|A^*|}$ by
$|\{f\in {}^{A^*}\beta:\rang(f)=\beta\}|$). 
\end{observation}

\begin{convention}
$\kappa,\sigma$ are alephs, $\sigma$ an uncountable cardinal such that
$DC_{<\sigma}$, $f$, $g$, $h$ will be functions from $\kappa$ to $\ord$.
\end{convention}

\begin{definition}
Assume $\cf(\kappa)\geq\theta$.
\begin{enumerate}
\item $P^{\sigma}_{\kappa,\alpha}=\{ f: f \mbox{ a function from some
ordinal }\gamma<\sigma \mbox{ into } \kappa\cup{}^{\kappa}\alpha\}$
partially ordered by $\subseteq$ is regarded as a forcing notion and
$G_{P^\sigma_{\kappa, \alpha}}$ or simply $G$ denote the generic.
\item $\name h =\bigcup G_{P^{\sigma}_{\kappa,\alpha}}$ (is forced to
be a function from $\sigma$ onto $\kappa\cup {}^{\kappa}\alpha$).
\item $\name{\bar \alpha}[G_{P^\sigma_{\kappa, \alpha}}]= \langle
\name \alpha_i: i<\sigma \rangle$ lists in increasing order 
$$
\{\sup(\kappa\cap\rang (\name h\restriction\gamma)):\gamma<\sigma\}.
$$
\item $\name Y=\{\langle\gamma , \name h (\gamma)\rangle : \gamma<\sigma,
\name h (\gamma)\in\kappa\} \cup\{\langle\gamma, i, \beta\rangle :
\gamma<\sigma, i<\kappa \mbox{ and } \name h(\gamma)\in
{}^{\kappa}\alpha , (\name h (\gamma))(i)=\beta\}$.
\item ${\name D}_{\kappa , \theta}^{\sigma}$ is the club filter on the
ordinal $\kappa$ in $K[{\name Y}]$ ($K$ - the core model, see Dodd and
Jensen [DJ]).
\item We can replace $\alpha$, ${}^{\kappa}\alpha$ by $f\in
{}^{\kappa}\ord$, $\prod\limits_{i<\kappa}(f(i)+1)$ respectively.
\end{enumerate}
\end{definition}
\medskip

\par \noindent
{\sc Remark:} On forcing for models of $ZF$ only see [Bu].
\medskip

\begin{observation}
\label{obs8:4}
\begin{enumerate}
\item $D^{\sigma}_{\kappa , \alpha}={\name D}^{\sigma}_{\kappa,\alpha}
\cap {\cal P}(\kappa)^{\v}$ does not depend on the generic subset $G$ of
$P_{\kappa ,\alpha}^{\sigma}$. 
\item If $\kappa$ is regular in $\v^{P^{\sigma}_{\kappa, \alpha}}$ (e.g.
$\kappa$ successor in \v ) \underline{then} $ D^{\sigma}_{\kappa , \alpha}$ is
normal, minimal among the normal filters on $\kappa$.
\end{enumerate}

\end{observation}

\begin{claim}
\label{claim8:5}
Assume that
\begin{description}
\item[$\boxtimes_{\sigma, \kappa, \alpha}$] for some/every
$G\subseteq P^{\sigma}_{\kappa,\alpha}$ generic over $\v$, in $K[{\name Y}]$
there are arbitrarily large Ramsey cardinals\footnote{or considerably less}.
\end{description}
\underline{Then} $\rk^2_{D^{\sigma}_{\kappa, \lambda}}(\alpha)<\infty$
for $E=E_0$ 
the family of filters on $\kappa$, or $E=E_1$ the family of normal filters on
$\kappa$ in the case of \ref{obs8:4}(2), or even 
$$
\begin{array}{ll}
E=E_2=
=\{D\colon\forces_{P^{\sigma}_{\kappa,\alpha}}& \mbox{ ``there is a normal
filter } D' \mbox{ on } \kappa\\
\ & \quad \mbox{ such that }D'\cap {\cal P}(\kappa)^{\v}=
D\mbox{''}\}.
\end{array}
$$
\end{claim}

\Proof As $P^{\sigma}_{\kappa, \alpha}$ is a homogeneous forcing, ``some $G$''
and ``every $G$'' are equivalent. The proof is as in [Sh-g], Ch.V 1.6, 2.9.
\hfill$\QED_{\ref{claim8:5}}$ 
\medskip

\par \noindent
{\sc Remark:} In [Sh-g] Ch.V we almost get away with $K(Y)$, $Y\subseteq
(2^{\aleph_1})^+$. We shall return to this later. 
\medskip

\begin{claim}
\label{claim8:6}
Assume that $\boxtimes_{\sigma, \kappa, \alpha}$ fails.
\begin{enumerate}
\item If in \v, $\lambda=\mu^+$, $\mu$ is singular $>\theta({}^{\sigma
>}({}^{\kappa}\alpha))$ \underline{then} $\lambda$ is regular but not
measurable. 
\item Bounds to cardinal arithmetic : if $\mu\geq\kappa$ \underline{then}
$\lambda^{\mu} \leq 2^{\mu}\times\lambda^+$
\end{enumerate}
\end{claim}

\Proof Let $G\subseteq P^{\theta}_{\kappa,\alpha}$ be generic over \v. First 
we shall prove that
\begin{description}
\item[$(*)_1$] \ \ \ \ \ \ $\lambda$ is still a cardinal in $\v[G]$.
\end{description}
Assume not. So $\v[G]\models |\lambda|=|\mu|$, and hence for some $P$-name
$\name H$ and $p\in P$, 
$$
p\forces_P\mbox{``}\name H\mbox{ is a one-to-one function from $\mu$ onto
$\lambda$''.}
$$
For $\alpha<\mu$ let 
$$
A_{\alpha}=\{\beta<\lambda : p\nforces_P {\name H}(\alpha)\neq\beta\}.
$$ 
So $\langle A_{\alpha}: \alpha<\mu\rangle\in\v$ and easily $\otp(A_{\alpha})<
\theta(P^{\sigma}_{\kappa,\alpha})$ (map $q\in P$ to $\beta$ if $p\leq q$,
$q\forces {\name H}(\alpha)=\beta$, and map $q\in P$ to $\min (A_\alpha)$
otherwise). But $\theta(P^{\sigma}_{\kappa,\alpha})\leq\mu$, so (as the sets
are well ordered) this implies 
$$
\v\models \lambda=|\bigcup\limits_{\alpha<\mu} A_{\alpha}|\leq
|\sum_{\alpha<\mu}\otp(A_\alpha)|\leq |\sum_{\alpha<\mu}\theta
(P^\sigma_{\kappa,\alpha})|\leq\mu\times\mu \leq\mu;
$$ 
a contradiction.  

Really we have proved 
\begin{description}
\item[$(*)_1^+$] \ \ \ every cardinal $\geq\theta({}^{\sigma}({}^{\kappa}
\alpha))$ in $\v$ is a cardinal in $\v[G]$.
\end{description}

Next we prove
\begin{description}
\item[$(*)_2$] $\lambda$ is regular in $\v$.
\end{description}
Assume not. So there are $A_2\subseteq\lambda =\sup (A_2)$, $\otp(A_2)
=\cf(\lambda)$, $A_0\subseteq\mu =\sup (A_0)$, $\otp (A_0)=\cf(\mu)$.  Let
$\mu^*= (\mu^+)^{K[{\name Y}[G]]}$ and let in \v\ we have $A_1\subseteq\mu^*
=\sup (A_1)$, $\otp (A_1)=\cf(\mu^*)$.  Let now $\bar A =\langle A_0, A_1,
A_2\rangle$. In $\v[G]$, $\mu$, $\lambda$ are still cardinals. Hence in
$K[{\name Y}[G]]$ too $\mu$, $\lambda$ are cardinals and even in $K[{\name
Y}[G], \bar A]$ we have $\lambda$ is a cardinal. Also $\mu^*\in
(\mu,\lambda)$, hence $\cf(\mu^*)<\mu$ hence
$K[{\name Y}[G], \bar A]\models \cf(\mu^*)<\mu$. Now $K[\name Y[G]]$, $K[\name
Y[G], \bar{A}]$ are models of $ZFC$, and in the latter
$\cf(\mu^*)\leq\otp(A_2)<\mu$. 

As $K[{\name Y}[G]]$ does not have unboundedly many Ramsey cardinals and
$\mu > \aleph_1^{\v}=\aleph_1^{K[{\name Y}[G]]}$, by Dodd and Jensen [DJ]
there is $B\in K[{\name Y}[G]]$ such that $K[{\name
Y}[G]]\models\mbox{ `` }|B|<(|\mu|^+)^{K[{\name Y}[G], \bar{A}]}$'' and
$A_1\subseteq B$, and so we get contradiction to $K[{\name
Y}[G]]\models\mbox{ ``}\mu^*\mbox{ a successor cardinal }>\mu\mbox{''}.$ 

Next:
\begin{description}
\item[$(*)_3$] if $\mu$ is singular, $\lambda$ is not measurable;
\end{description}
[as by the proof of $(*)_2$, $K[{\name Y}[G]]\models\mbox{``}\lambda=\mu^+
\mbox{''}$ so there is $e$ as in section \ref{par4}].
\hfill$\QED_{\ref{claim8:6}}$ 

\begin{claim}
\label{claim8:7}
[$AC_{{\cal P}({\cal P}(\kappa))}+DC$]\ \ \ \ \ Assume $E$ is a family
of filters on $\kappa$, $D_*=\min (E)$ (that is $D_*\in E$, and $[D\in
E\Rightarrow D_*\subseteq D]$) and for
some $f^*\in {}^{\kappa}\ord$, $\rk^2_{D_*}(f^*)=\infty$. \underline{Then}
$\rk^2_{D_*}(\alpha)=\infty$ for some $\alpha, \alpha<\theta({\cal P}({\cal
P}(\kappa))$ (and essentially $\alpha<\theta({\cal P}(\kappa))$) (i.e.
$\alpha$ stands 
for the function  from $\kappa$ to $\ord$ which is constantly $\alpha$).
\end{claim}

\Proof We first prove that there is such $\alpha<\theta({\cal P}({\cal
P}(\kappa)))$. If $D\in E$, $\rk^2_{D}(f)=\infty$ then for every ordinal
$\alpha$ there are $A= A_{f,\alpha}\in D^+$ and
$g=g_{f,\alpha}<_{D+A_{f,\alpha}} f$ such that $\rk^3_{D+A}(g)>\alpha$.
Without loss of generality $g_{f,\alpha}\leq f$; so we have only a set of
possible $(A_{f,\alpha}, g_{f,\alpha})$ (for given $f$).

\noindent By the ``$F$'' of $ZF$ there are $A_f$, $g_f<_{D+A_{f}} f$,
$g_f\leq f$ such that $\rk^3_{D+A_f}(g_f)=\infty$, so
$\rk^2_{D_1}(g_f)=\infty$ for every $D_1$, $D+A_f\subseteq D_1\in E$. Note: 
$$
\kappa+\kappa=\kappa,\quad |{\cal P}(\kappa)|^2=|{\cal P}(\kappa)|,\quad
|{\cal P}({\cal P}(\kappa))|^2=|{\cal P}({\cal P}(\kappa))|
$$ 
(really $|A|+1=|A|\Rightarrow |{\cal P}(A)|^2=|{\cal P}(A)|$). By $AC_{{\cal
P}({\cal P}(\kappa))}$: 
\begin{description}
\item[$(*)$] \ \ \ if $F\subseteq {}^{\kappa}\ord$ has cardinality $<|{\cal
P}({\cal P}(\kappa))|$ then we can find 
\[\{\langle A_{f, D}, g_{f,D}\rangle : f\in F, D\in E, \rk^2_D
(f)=\infty\}\]
such that $\rk^3_{D+A_{f, D}} (g_{f, D})=\infty$, $g_{f, D}<_{D+A_{f, D}}
f$, $g_{f,D}\leq f$. 
\end{description}
By $DC$ we can find $\langle F_n:n<\omega\rangle$, $F_n\subseteq {}^{\kappa}
\ord$, such that in $(*)$ for $F=F_n$ we have $F_{n+1}= F_n\cup\{g_{f, D}:
f\in F_n, D\in E\}$, and $F_0=\{f^*\}$. Let $F=\bigcup\limits_{n<\omega} F_n$,
and let $\bar A=\langle A_i: i<\kappa\rangle$ be defined by $A_i=\{f(i): f\in
F\}$. Clearly $|F_n|\leq |E_n|\leq|{\cal P}({\cal P}(\kappa))|$ and hence
$|A_i|\leq |{\cal P}({\cal P}(\kappa))|$. Hence $h^*(i)=:\otp
(A_i)<\theta({\cal P}({\cal P}(\kappa)))$ and $\sup\rang(h^*)<\theta({\cal
P}({\cal P}(\kappa)))$.  Now for every $f\in F$ we define $h_f \in
\prod\limits_{i<\kappa} g(i)$ by $h_f (i)=\otp(A_i\cap f(i))< g(i)$. It is now
easy to check that $\rk^2_{D_*} (h_{f^*})=\infty$ as required.

Next we prove that for some $\alpha<\theta({\cal P}(\kappa))$ we have
$\rk^2_{D_*} (\alpha)=\infty$.  As $AC_{{\cal P}({\cal P}(\kappa))}$, clearly
${\cal P}(\kappa)$ can be well ordered, so let $h_0 : {\cal
P}(\kappa)\rightarrow |{\cal P}(\kappa)|$ be one-to-one onto $|{\cal
P}(\kappa)|$, an aleph.  By the first part there is $\alpha^* < \theta({\cal
P}({\cal P}(\kappa)))$ with $\rk^2_{D_*} (\alpha^*)=\infty$. Hence we can find
$h_1^*: {\cal P}({\cal P}(\kappa)) \stackrel{\rm onto}{\longrightarrow}
\alpha^* +1$.
\medskip

\noindent{\sc Observation \ref{claim8:7}A}: $[AC_\kappa]$\ \ \ \ If $\alpha
<\theta({\cal P}(A^*))$, $|A^*|\times \kappa=|A^*|$ then ${}^{\kappa}
(\alpha+1)$ is a set of cardinality $\leq^* {\cal P}(A^*)$.
\smallskip

\noindent [Why? Let $h_1: {\cal P}(A^*)\stackrel{\rm onto}{\longrightarrow}
\alpha+1$ and let $h_2: A^*\times \kappa \rightarrow A^*$ be one to
one. For each $f\in {}^{\kappa} (\alpha+1)$ we can choose $\bar B=
\langle B_i: i<\kappa\rangle$, $h_1 (B_i)=f(i)$, so $B_i\subseteq
{\cal P} (A^*)$, and $\langle B_i: i<\kappa\rangle$ encodes $f$ and it
in turn can be coded by $\{ (i,x): i<\kappa\mbox{ and } x\in B_i\}\subseteq
\kappa\times A^*$, but it is not clearly if we have also the function
$f \mapsto \bar B^f$. However we can define a function $H$,
$\dom(H)={\cal P}(\A)$ and $\rang(H)\subseteq {}^\kappa(\alpha+1)$ by 
$$
(H(A))(i)= h_1(\{a\in \A: h_2(a, i)\in A\}).
$$
Clearly $H$ is a function from ${\cal P}(A)$ onto
${}^\kappa(\alpha+1)$, hence $|{}^\kappa(\alpha+1)|\leq^* |{\cal
P}(\A)|$.]\\
\null\hfill\qed$_{\ref{claim8:7}A}$ 
\medskip

Of course, apply \ref{claim8:7}A to $A^*={\cal P}(\kappa)$.

Also $|\fil(\kappa)|\leq |{\cal P}({\cal P}(\kappa))|$ and $|{\cal P}({\cal
P}(\kappa))|^2=|{\cal P}({\cal P}(\kappa))|$, so by $AC_{{\cal P}({\cal
P}(\kappa))}$ we can find $Y_0=\langle (A_{f,D}, g_{f,D}): f\in {}^{\kappa}
(\alpha^* +1), D\in E\rangle$ such that:
\begin{quotation}
\noindent$(*)$\ \ \ if $\rk^2_D (f)=\infty$ then $\rk^3_{D+A_{f,D}}
(g_{f,D})=\infty$, $A_{f,D} \in D^+$,
 
$g_{f,D} \in {}^{\kappa}(\alpha^* +1)$, $g_{f,D}<_{D+A_{f,D}} f$
\end{quotation}
Similarly by $AC_{{\cal P}({\cal P}(\kappa))}$ we can find $Y_1=\langle
d_{f,D}:  f\in {}^{\kappa} (\alpha^* +1), D\in E\rangle$ such that:
\begin{quotation}
$D\subseteq d_{f,D}\in E$ and $\rk^3_D (f)=\rk^2_{d_D} (f)$.
\end{quotation}
We now define the model $\gc$ with the universe ${\cal P}({\cal
P}(\kappa))$ - just put all the information needed below.

Clearly $|\gc|=|{\cal P}({\cal P}(\kappa))|$, so we have a choice function
$H^*$ on the family of definable (in $\gc$) non empty subsets of $\gc$. So for
$A\subseteq \gc$ we have the Skolem hull. Now we define by induction on
$\alpha\leq\kappa^+$ (an aleph) submodels $M_{\alpha}$ of $\gc$ increasing
continuously in $\alpha$ and of cardinality $\leq |{\cal P}(\kappa)|$.
\begin{quotation}

\noindent \underline{For $\alpha=0$}: the Skolem hull of
$\{\gamma:\gamma\leq |{\cal P}(\kappa)|\}$ 

\noindent \underline{For $\alpha$ limit}:
$\bigcup\limits_{\beta<\alpha} M_{\beta}$ 

\noindent \underline{For $\alpha=\beta+1$}: the Skolem hull of
$|M_{\alpha}|\cup\{x:x\in {}^{\kappa}|M_{\alpha}|\}$ 

\end{quotation}
Now ${}^{\kappa}(M_{\kappa^+})\subseteq M_{\kappa^+}$ (as $\kappa^+$ is
regular as $|{\cal P}(\kappa)|$ is an aleph). Let $H:M_{\kappa^+}\cap\ord \stackrel{\rm onto}{\longrightarrow}
\beta^*< |{\cal P}(\kappa)|^+$ be order preserving. So as we are assuming that
the conclusion fails, 
for every $f\in {}^{\kappa}(\alpha^*+1)\cap
M_{\kappa^+}$,  
$$
H\circ f^* \in {}^{\kappa}\beta^*\mbox{ satisfies }\rk^2_{D_*}
(H\circ f^*)<\infty.
$$  
Whereas $\rk^2_{D_*}(f^*)=\infty$. We can conclude as
in [Sh 386, 1.13]. \hfill$\QED_{\ref{claim8:7}}$ 

We can generalize \ref{obs8:1}-\ref{claim8:7}
\begin{claim}
\label{claim8:8}
[$DC_{<\sigma}$, $\sigma=\cf(\sigma)>\aleph_0$]\ \ \ \ Assume:
\begin{description}
\item[(a)] $D_*$ is a $\sigma$-complete filter on $A^*$
\item[(b)] $A^*=A_*\times\sigma$, and $\i:A^*\rightarrow\sigma$ is
$\i(a,\alpha)=\alpha$ 
\item[(c)] $D^*$ is the following filter on $A^*$; for $A\subseteq
A^*$:
\begin{quotation}
\noindent $A\in D^* \iff \{i<\sigma :\{a\in A_*: (i,a)\in A^*\}\in D_*\}\in 
D_{\sigma}$
where $D_{\sigma}$ is the minimal normal filter on $\sigma$.
\end{quotation}
\item[(d)] $E$= the family of $\i$-normal filters on $A^*$.
\end{description}
Then for every $f\in {}^{A^*}\alpha(*)$, $\rk^2_D (f)<\infty$,
provided that $\otimes^\sigma_{|A^*|,\alpha(*)}$.
\end{claim}
\medskip

\par \noindent
{\sc Concluding Remark}: We can deal with ${\rm fc}
({}^{\omega>}\omega)$ as in [Sh 386], [Sh 420]. Also concerning
\ref{claim8:7} change $|A^*|$. 

Also as in [Sh 386] we can define ${}^1\rk^2 (f,{\cal E})$ and
so improve \ref{conclusion3:??} (as in [Sh 386,\S 2]).
\eject

\centerline{\sc REFERENCES}
\bigskip
\bigskip
\bigskip
\begin{description}
\item[\mbox{[ApMg]}] A.~Apter, M.~Magidor, {\em Instances of Dependent
Choice and the Measurability of $\aleph_{\omega + 1}$}, Annals of Pure
and Applied Logic, submitted.
\item[\mbox{[Bu]}] E.~Bull, {\em Consecutive Large Cardinals}, Annals
of Mathematical Logic, 15(1978):161-191.
\item[\mbox{[DJ]}] A.~Dodd, R.~B.~Jensen, {\em The core model}, Annals
of Mathematical Logic, 20(1981):43-75.
\item[\mbox{[Gi]}] M.~Gitik, {\em All uncountable cardinals can be
singular}, Israel J. of Mathematics, 35(1980):61-88.
\item[\mbox{[GH]}] F.~Galvin, A.~Hajnal, {\em Inequalities for
cardinal powers}, Annals of Mathematics, 101(1975):491-498.
\item[\mbox{[J]}] T.~Jech, {\bf Set Theory}, Academic Press, New York
1978.
\item[\mbox{[J1]}] T.~Jech, {\bf Axiom of Choice}.
\item[\mbox{[Kf]}] G.~Kafkoulis, {\em Homogeneous Sequences of
Cardinals for OD Partition Relations}, Doctoral Dissertation, Cal.
Tech., 1989. 
\item[\mbox{[Sh-g]}] S.~Shelah, {\bf Cardinal Arithmetic}, Oxford Logic
Guides vol. 29(1994), General Editors: Dov M. Gabbai, Angus Macintyre, Dana
Scott, Oxford University Press.
\item[\mbox{[Sh 176]}] S.~Shelah, {\em Can you take Solovay inaccessible
away?}, Israel J. of Mathematics, 48(1984):1-47.
 \item[\mbox{[Sh 386]}] S.~Shelah, {\em Bounding $pp(\mu)$ when cf$(\mu)>
\aleph_0$ using ranks and normal ideals}, Chapter V of {\bf Cardinal
Arithmetic}, Oxford Logic Guides vol. 29(1994), General Editors: Dov M.
Gabbai, Angus Macintyre, Dana Scott, Oxford University Press.
\item[\mbox{[Sh 420]}] S.~Shelah, {\em Advances in Cardinal Arithmetic},
{\bf Finite and Infinite Combinatorics in Sets and Logic}, N.W. Sauer et al
(eds.), Kluwer Academic Publishers, 1993:355-383.
\item[\mbox{[Sh 506]}] S.~Shelah, {\em The pcf-theorem revisited}, {\bf
Mathematics of Paul Erd\H{o}s}, vol 2, edited by R. Graham and J. Nesetril,
in print. 
\end{description}
\shlhetal
\end{document}